\documentclass[11pt, a4paper]{article} 
\usepackage[utf8]{inputenc}

\usepackage{float}
\usepackage{tikz-cd}
\usepackage{amsmath}
\usepackage{amssymb}
\usepackage{amsfonts}
\usepackage{amsthm}
\usepackage{amscd}
\usepackage{extarrows}
\usepackage{mathtools}
\usepackage[all,2cell]{xy}
\usepackage[shortlabels]{enumitem}
\usepackage{hyperref}
\usepackage{euscript}
\usepackage[margin=3cm]{geometry}
\usepackage[bbgreekl]{mathbbol}
\usepackage{graphicx}
\graphicspath{{pics/}}

\counterwithin{figure}{subsection}

\usepackage{xstring} 

\usepackage{cleveref}%
\usepackage{autonum}

\usepackage{macros/utility-macros}
\usepackage{macros/diagram-macros}
\usepackage{macros/macros}
\usepackage{macros/environments}
\usepackage{macros/lax-add-macros}

\numberwithin{equation}{section}

\usepackage{xcolor}%
\hypersetup{%
    colorlinks,%
    linkcolor={red!50!black},%
    citecolor={blue!50!black},%
    urlcolor={blue!80!black}%
}

\setlist[enumerate,1]{label=(\arabic{*})}
\setlist[enumerate,2]{label=(\alph{*})}
\setlist[enumerate,3]{label=(\roman{*})}

\usepackage{subfiles}

\title{Lax Additivity} 
\author{
  Merlin Christ\footnote{Universit\"at Bonn, Mathematisches Institut, Endenicher Allee 60, 53115 Bonn, Germany, email: {\tt christ@math.uni-bonn.de}}
  , Tobias Dyckerhoff\footnote{Universität Hamburg, Fachbereich Mathematik, Bundesstraße 55, 20146 Hamburg, Germany, email: {\tt tobias.dyckerhoff@uni-hamburg.de}}
  , Tashi Walde\footnote{Universität Regensburg, Fakultät für Mathematik,  93040 Regensburg, Germany, email: {\tt tashi.walde@ur.de}}
}

\begin{document}

\maketitle

\begin{abstract}
  We introduce notions of lax semiadditive and lax additive \((\infty,2)\)-categories,
  categorifying the classical notions of semiadditive and additive $1$-categories.
  To establish a well-behaved axiomatic framework, we develop a calculus of lax
matrices and use it to prove that in locally cocomplete
\((\infty,2)\)-categories lax limits and lax colimits agree and are absolute.
  In the lax additive setting,
  we categorify fundamental constructions from homological algebra
  such as mapping complexes and mapping cones
  and establish their basic properties.
\end{abstract}

\tableofcontents

\noindent\\

\section{Introduction}

In this article, we propose an axiomatic framework of \emph{lax additive}
\((\infty,2)\)-categories, intended as a natural context to develop
foundational aspects of categorified homological algebra (analogously to the
familiar development of classical homological algebra building on additive
categories). 

Our motivation stems from several recent developments, some of the most
directly relevant ones being:

\begin{itemize}
	\item Categorified analogs of classical homological techniques have
	been very successful in the study of Fukaya--type categories in homological
	mirror symmetry. The categorical Picard-Lefschetz theory developed in
	\cite{Sei08} is a particularly well-proven example.

	\item Kapranov and Schechtman have proposed to study categorified
	analogs of perverse sheaves, termed perverse schobers \cite{KS14}. While this
	beautiful circle of ideas has already created substantial impact, the theory is
	still somewhat experimental and as of now there does not seem to exist a
	satisfying rigorous definition of perverse schobers in some natural generality. 

	\item Various foundational results from classical homological algebra
	have been shown to admit categorified variants replacing abelian groups by stable
	$\infty$-categories. An illustrative example is the categorified Dold-Kan
	correspondence (cf. \cite{dyck:dk, heine:dk}) which can be regarded as a 
	``proof of concept'' for the feasibility of categorifying some of the
	foundations of homological algebra. 

	\item Several examples of stable $\infty$-categories of algebraic or
	geometric origin have been shown to admit natural upgrades to complexes of
	stable $\infty$-categories (cf. \cite{CDW23}).
\end{itemize}

We see this work as a first step towards capturing the basic $2$-categorical
principles underlying these perspectives, with the final goal of creating a
unified picture of mutual benefit.
That being said, we recommend to read this paper as a companion to \cite{CDW23},
where the abstract axiomatic theory developed here appears in a very hands-on
way, illustrated by many examples and explicit constructions. 

Beyond these concrete applications, we feel that the $2$-categorical theory of
lax additivity developed in this work does have some intrinsic category
theoretic appeal, justifying its documentation in a standalone contribution. 
For example, we systematically introduce various types of lax matrices along
with categorified matrix multiplication rules. Based on this calculus, we prove
natural categorified variants of classical foundational results on (semi--)additive categories such as:
\begin{thm*}
  [\Cref{cor:lax-co-limits-are-bilimits} and \Cref{thm:lax-lims-absolute}]
  In locally cocomplete\footnote{
    An \((\infty,2)\)-category is locally cocomplete
    if all its hom-categories have colimits
    and if composition of \(1\)-arrows preserves colimits in each variable.
  } \((\infty,2)\)-categories
  \begin{itemize}
  \item
    lax limits and lax colimits coincide (when they exist) and
  \item
    all lax limits and lax colimits are absolute, i.e.\
    preserved by locally cocontinuous\footnote{
      A functor of \((\infty,2)\)-categories is locally cocontinuous
      if it induces a colimit-preserving functor on hom-categories.
    } functors.
  \end{itemize}
\end{thm*}

As a further illustration of the theory, we categorify basic additive
constructions from homological algebra such as mapping complexes and mapping
cones and establish their basic properties. 

\subsection{Rules of categorification}
\label{subsec:rules-of-categorification}

We begin with an overview of our preferred type of ``categorification'': It
arises from the insight that in some important respects stable
$\infty$-categories behave like categorified abelian groups, leading to the
``categorification rules'' in Table \ref{table:1}.
\begin{table}[h]
  \renewcommand{\arraystretch}{2}
  \begin{center}
    \begin{tabular}{lll}
      & classical & categorified\\ \hline
      1) & abelian group $A$ & stable $\infty$-category $\A$
      \\
      2) & element $x \in A$ & object $X \in \A$
      \\
      3) &  $y - x$ & $\cone(X \overset{f}{\to} Y)$
      \\
      4) &  $\sum (-1)^i x_i$
                  & $\tot(
                    X_0 \xrightarrow{d}  X_1 \xrightarrow{d} X_2\xrightarrow{d}
                    \cdots\xrightarrow{d} X_n)$
      \\
      5) & \parbox{6cm}{
           direct sum decomposition
      \\
      $C \cong A\oplus B$
      }
      &
        \parbox{6cm}{
        semiorthogonal decomposition
      \\
      $\C \simeq \langle \A, \B \rangle$
      }\\
      6) & external direct sum $A\oplus B$
                  & gluing along a functor/lax sum $\laxsum{\A}{F}{\B}$
    \end{tabular}
  \end{center}
 \caption{Categorification rules}
\label{table:1}
\end{table}

While the first two rules should be apparent, we start commenting on rule 3).
This is a first crucial difference between the classical and the categorified context:
In order to take a ``difference'' between objects $X,Y$ of a stable $\infty$-category $\A$,
we need to be given the additional datum of a morphism $f\colon X \to Y$%
---the difference will then be the cone of $f$.
Compliance with this rule will force us to include certain
$2$-categorical data which becomes invisible upon passing to the Grothendieck group $K_0$.
This typically results in rather natural \emph{lax variants}
of $1$-categorical constructions.  

Rule 4) is a natural generalization of Rule 3):
An alternating sum over $n$ elements will be categorified
by the totalization of an $n$-term complex in $\A$.
Here we do not only need to specify the differentials of this complex,
but also a coherent system of null homotopies%
---this is necessary to make sense of the totalization in the $\infty$-categorical context.

Rule 5) is almost evident after having accepted Rule 3):
While in a direct sum $A \oplus B$,
every element is uniquely the sum of elements from the components $A$ and $B$,
respectively, in a semiorthogonal decomposition
$\langle \A, \B \rangle$,
every object is uniquely an \emph{extension}
of an object $A \in \A$ by an object $B \in \B$.
Put differently, by shifting the exact triangle of the extension,
every object is uniquely the cone of a morphism $A[-1] \to B$,
thus connecting back to Rule 3). 

Conceptually distinct to a direct sum \emph{decomposition}
of a given abelian group are the universal properties satisfied by the
\emph{external} direct sum of a pair of abelian groups.
For its categorificaiton, it is not sufficient to just provide a pair of stable $\infty$-categories. As an additional datum, we need to specify a functor
$F\colon \A \rightarrow \B$
(similar to the additional choice of a morphism $f\colon X \to Y$ needed in Rule 3).
The categorified ``direct sum'' is then the \emph{lax sum} 
\begin{equation}
	\laxsum{\A}{F}{\B}
\end{equation}
of the diagram of stable $\infty$-categories described by $F$ (see Rule 6),
which is given by the commonly known construction of gluing along a functor.
The fact that this sum categorifies both the product and coproduct will be explained below in the context of lax additivity.
The lax sum admits a semiorthogonal decomposition with components $\A$ and $\B$, and vice versa,
any stable $\infty$-category with a semiorthogonal decomposition
can be described as a lax sum if and only if the semiorthogonal decomposition admits a \emph{gluing functor}. 

\subsection{Lax additivity}\label{intro:laxadd}

Of course abelian groups are rarely studied in isolation;
rather we consider the category \(\Ab\) of abelian groups.
This category has many important features,
but most importantly for us it is \emph{additive}.
Continuing our train of thought,
Table \ref{table:1} has a natural continuation in Table \ref{table:2} which explains what it means to say that
the \((\infty,2)\)-category of stable \(\infty\)-category,
or more generally any \((\infty,2)\)-category \(\exaddtwocat\), is \emph{lax additive}.
\begin{table}[h]
\renewcommand{\arraystretch}{2.5}
\begin{center}
  \begin{tabular}{lll}
    & classical & categorified\\ \hline
    7) & additive (\(\infty\)-)category \(\exaddcat\) & lax additive \((\infty,2)\)-category \(\exaddtwocat\)\\
    8) & \parbox{5.5cm}{hom-sets \(\exaddcat(X,Y)\) have addition} & \parbox{5.5cm}{hom-categories \(\exaddtwocat(\X,\Y)\) have colimits}\\
    9) & \parbox{5.5cm}{hom-sets \(\exaddcat(X,Y)\) are abelian groups} & \parbox{5.5cm}{hom-categories \(\exaddtwocat(\X,\Y)\) are stable}\\
     10) & \parbox{5cm}{finite direct sums\\\(\bigoplus\limits_{s=1}^kx_s=\prod\limits_{s=1}^kx_s=\coprod\limits_{s=1}^kx_s\)} &
         \parbox{6cm}{general lax bilimits\\\(\laxbilim\limits_{s:S}\X_s=\laxlim\limits_{s:S}\X_s=\laxcolim\limits_{s:S}\X_s\)}\\   
         11) & \parbox{5.5cm}{binary direct sums\\\(X\oplus Y =X\times Y=X\amalg Y\)} &%
         \parbox{6cm}{lax/oplax \(\Delta^1\)-bilimits\\
    \(
    \laxlima{\X}{F}{\Y}=
    \laxcolima{\X}{F}{\Y}=
    \oplaxlima{\X}{F}{\Y}=
    \oplaxcolima{\X}{F}{\Y}
    \)}\\
    12) & \parbox{5cm}{matrices
         \(
         \begin{psmallmatrix}
           m_{11}& m_{12}\\
           m_{21}& m_{22}
         \end{psmallmatrix}
         \)} &
               \parbox{6cm}{lax matrices
               \(
               \Donematrix{M_{11}}{M_{12}}{M_{21}}{M_{22}}
               \)
               }
    \\
    13) &
         \parbox{5cm}{matrix multiplication\\
         \((nm)_{us}=\sum_{t=t'}(n_{ut'}\circ m_{ts})\)
         }
                &
                  \parbox{6.5cm}{lax matrix multiplication\\
                  \((NM)_{us}=\colim\limits_{\gamma\colon t\to t'}(N_{ut'}\circ \gamma\circ M_{ts})\)
                  }
    \\
    14) &
         \parbox{6cm}{matrix multiplication, reparameterized\\
    \((nm)_{us}=\sum_{t}(-1)^t(n_{ut}\circ m_{ts})\)
    }
    &
      \parbox{6cm}{lax-oplax matrix multiplication\\
    \((NM)_{us}=\tot_{t}(N_{ut}\circ M_{ts})\)
    }
  \end{tabular}
\end{center}
 \caption{Categorification rules (cont.)}
\label{table:2}
\end{table}

Accepting our basic premise that abelian groups are to be categorified by stable
\(\infty\)-categories, Rule 9) requires no further comment.
Rule 8) is a convenient intermediate step,
categorifying the situation where the addition on hom-sets does not necessarily have inverses;
just like the uncategorified case,
many basic lemmas are most naturally expressed in this generality
leading to the notion of lax \emph{semi}-additive \((\infty,2)\)-category.

The direct sum of abelian groups is both a categorical product and a categorical coproduct,
a universal property that is taken as the definition in general additive categories.
Rule 10) states that the same definition can be categorified, by replacing finite products and coproducts,
i.e., limits and colimits indexed by finite discrete categories
\(S=\set{1,\dots,k}\),
with \emph{lax} limits and colimits indexed by arbitrary \(\infty\)-categories \(S\).
Apart from this change, the theory is exactly analogous:
if the hom-categories have colimits (categorifying addition)
then such lax limits and colimits always agree if they exits,
yielding the notion of lax bilimits.
Thus we obtain the main concept of this article:

\begin{defi*}[\Cref{defi:lax-additive}]
  An \((\infty,2)\)-category \(\exaddtwocat\) is lax additive if
  \begin{itemize}
  \item
    it is enriched in stable \(\infty\)-categories with colimits and
  \item
    it has all lax limits and lax colimits, which then automatically agree.
  \end{itemize}
\end{defi*}

In a lax additive \((\infty,2)\)-category,
we can say even more for the special choice of \(S=\Delta^1\),
in which case the constructions categorify the binary direct sum.
For this choice of $S$, all four possible universal 2-categorical constructions
(lax/oplax, limit/colimit)
associated to an \(S\)-diagram \(\X\xrightarrow{F}\Y\) agree with each other.
This is the content of Rule 11).
Note that this further explains the statement of Rule 6)
and this (op)lax bilimit is also called the lax sum $\laxsum{\X}{F}{\Y}$.

A convenient feature of additive categories is that maps
\(m\colon x_1\oplus x_2\to y_1\oplus y_2\)
between direct sums can be represented as matrices of the form
\begin{equation}
  {
    m =
    \begin{pmatrix}
      m_{11}\colon x_1\to y_1, & m_{12}\colon x_2\to y_1\\
      m_{21}\colon x_1\to y_2, & m_{22}\colon x_2\to y_2
    \end{pmatrix}.
  }
\end{equation}
Composing such maps then just amounts to the usual matrix multiplication.
Rule 13) shows how the usual matrix multiplication formula can be categorified,
yielding an analogous theory of matrices indexed in each coordinate
not by a finite set but by arbitrary \(\infty\)-categories.
These matrices are just a dependent version of bimodules,
which by Morita theory encode functors between module categories.

In \S\ref{subsec:rules-of-categorification} we have already seen how it is
conceptually easier to categorify subtraction rather than addition and more
generally alternating rather than ordinary sums.
In the lax additive setting we see a similar feature,
expressed in Rule 14),
where a suitable ``coordinate change'' yields
a more convenient formula for the categorified matrix multiplication
when we reparameterize it to use alternating sums rather than ordinary sums.
Categorically speaking this reparameterization
involves the identification of lax and oplax limits
and is therefore only available for certain special indexing categories
such as \(\Delta^1\).

\subsection{Mapping cones and mapping complexes}

Building on the notions introduced above, we explain how to categorify two
fundamental constructions within the lax additive framework: the mapping
complex between two chain complexes (see \Cref{con:lax-mapping-complex}) and
the mapping cone of a chain map (see \Cref{con:oplax-mapping-cone}).

We summarize the categorified formulas for their respective differentials in
Table \ref{table:3}.  We note that upon passing to \(K_0\), the signs in the
right hand side and left hand side differ only in a non-essential way, the
chain complexes are isomorphic.

\begin{table}[h]
\renewcommand{\arraystretch}{2.5}
\begin{center}
  \begin{tabular}{lll}
    & classical & categorified\\ \hline
    15)
    & \parbox{6cm}
      {
      \(\delta g_\bullet \coloneqq (dg + (-1)^\bullet gd)\)
      }
    &
      \parbox{6cm}
      {
      \(
      \delta^\lax g_\bullet \coloneqq \fib(dg \to gd)[\bullet]
      \)
      }
    \\
    16)
    & \parbox{6cm}
      {%
      \(d_{\operatorname{Cone}(f)}\coloneqq
      \begin{psmallmatrix}
        -d& 0\\
        -f& d
      \end{psmallmatrix}
      \)
      }
    &
      \parbox{6cm}
      {%
      \(
      d_{\operatorname{Cone}(f)}\coloneqq
      \Donematrixmixed{d}{0}{f}{d}
      \)
      }
  \end{tabular}
\end{center}
\caption{Categorification rules (cont.)}
\label{table:3}
\end{table}

We further note that unlike in the classical setting,
we do not require any signs in the differential of \(\operatorname{Cone}(f)\);
the appropriate signs guaranteeing \(\delta^2=0\)
are inbuilt into the alternating sums
defining the lax-oplax matrix multiplication.

The categorified mapping cone yields a natural notion
of null-homotopy \(H\) of a chain map \(f\),
given by \(\delta^\lax(H)=f\).
The categorified mapping cone interacts with this notion as expected:

\begin{thm*}[\Cref{cor:univ-prop-lax-cof-fib}]
  For every chain map \(f\colon \A\to \B\)
  for which sufficient adjoints exist
  and for each chain complex \(\C\),
  there is an equivalence between the (stable) \(\infty\)-categories of
  \begin{itemize}
  \item
    chain maps \(\operatorname{Cone}(f)\to \C\) and
  \item
    chain maps \(g\colon \B\to \C\)
    together with a lax null-homotopy \(H\) of \(gf\).
  \end{itemize}
\end{thm*}

The above mapping cone construction is only possible if certain adjoints of some of the involved functors exist. However, one can further analyze which data is precisely necessary to construct mapping cones, and give a more general construction of the categorified mapping cone, which takes as input extra data containing a degree $1$ map $\A\to \B$. We also prove a more general version of the above theorem, classifying maps out of this generalized mapping cone, see \Cref{thm:univ-prop-of-lh-cone}. In the case that sufficient adjoints exist, there exists a canonical choice for this extra data and \Cref{thm:univ-prop-of-lh-cone} specializes to the above theorem.

\subsection{Acknowledgements}

T.W.\ thanks Claudia Scheimbauer for conversations in the context of their joint work
that helped shape many of the key ideas regarding higher categorical additivity.
M.C.\ was funded by the Deutsche Forschungsgemeinschaft (DFG,
German Research Foundation) under Germany's Excellence Strategy -- EXC-2047/1 -- 390685813.
M.C.\ has received funding from the European Union’s Horizon 2020 research and innovation programme under the Marie Skłodowska-Curie grant agreement No 101034255.
M.C.\ and T.D.\ acknowledge support by the Deutsche Forschungsgemeinschaft under Germany’s Excellence
Strategy -- EXC 2121 “Quantum Universe” -- 390833306. 
T.D.\ acknowledges support of the VolkswagenStiftung through the Lichtenberg Professorship Programme.
T.W. acknowledges support by the Deutsche Forschungsgemeinschaft (DFG, German
Research Foundation) -- SFB 1085 -- ``Higher Invariants'' -- 224262486.
T.D. acknowledges support by the Deutsche Forschungsgemeinschaft (DFG,
German Research Foundation) – SFB 1624 – ``Higher structures, moduli spaces and
integrability'' – 506632645.


\section{Additive 1-categories}

To explain our philosophy,
let us first remind the reader of the classical story for ordinary additive categories.

We start by recalling the definition.

\begin{defi}
  A category \(\exaddcat\) is called \emph{additive} if:
  \begin{enumerate}
  \item
    \label{it:add:cmon-enrich}
    The category \(\exaddcat\) is enriched in abelian monoids;
    i.e.,
    each hom-set \(\exaddcat(x,y)\)
    has an associative, commutative addition \(+\) with neutral element \(0\)
    such that composition
    \begin{equation}
      \exaddcat(x,y)\times \exaddcat(y,z)\to\exaddcat (x,z)
    \end{equation}
    preserves \(+\) and \(0\) in each argument.
  \item
    \label{it:add:cab-enrich}
    Each commutative monoid \((\exaddcat(x,y),+,0)\)
    admits negatives, hence is an abelian group.
  \item
    \label{it:add:prod+coprod}
    The category \(\exaddcat\) admits finite products and coproducts
    (including empty ones).
  \item
    \label{it:add:prod=coprod}
    For each finite set of objects \(x_1,\dots, x_n\in\exaddcat\),
    the natural map
    \begin{equation}
      \label{eq:canmap-prod-coprod}
      \coprod_{s=1}^nx_s\longrightarrow\prod_{t=1}^nx_t,
    \end{equation}
    whose components \(x_s\to x_t\) are
    \begin{equation}
      \label{eq:cmap-coprod-prod}
      \begin{cases}
        1\colon x_s\to x_s&\text{, if } s=t\\
        0\in\exaddcat(x_s,x_t)&\text{, otherwise}
      \end{cases},
    \end{equation}
    is an isomorphism.
  \end{enumerate}
  A category only satisfying \ref{it:add:cmon-enrich}, \ref{it:add:prod+coprod}
  and \ref{it:add:prod=coprod} is called semiadditive.
\end{defi}

One typically identifies finite products and coproducts via the canonical map
\eqref{eq:canmap-prod-coprod}
and uses the symbol \(\oplus\) (called direct sum or biproduct) for both.

The use of the phrase ``is called additive if'' implies that being additive
is a property of the category \(\exaddcat\) rather than extra structure.
This is justified by the fact that the addition on the hom-sets of
an additive category is uniquely determined.
Explicitly, it is given by the following formula:
Given two maps \(f,g\colon x\to y\),
their sum is the composite
\begin{equation}
  x\to x\oplus x\xrightarrow{f\oplus g} y\oplus y\to y
\end{equation}
where the first map is the diagonal \(x\to x\times x\)
and the last map is the codiagonal \(y\amalg y\to y\).

In this sense, the biproduct structure \(\oplus\)
determines the addition structure \(+\) on the hom-sets.
The converse is also true, as explained by the following lemma:

\begin{lem}
  \label{lem:prod-coprod-equations}
  Let \(\exaddcat\) be a category enriched in abelian monoids.
  Let \(x_1,\dots, x_n\) be a finite set of objects in \(\exaddcat\).
  \begin{enumerate}
  \item
    \label{it:PI-prod-coprod}
    Let \(x\) be an object of \(\exaddcat\) equipped with a cone
    \(P=(p_s\colon x\to x_s)_{s=1}^n\)
    and a cocone
    \(I=(i_s\colon x_s\to x)_{s=1}^n\)
    satisfying the two equations
    \begin{enumerate}
    \item
      \label{it:IP=1}
      \begin{equation}
        \sum_{s=1}^n i_s\circ p_s = 1 \in \exaddcat(x,x)
      \end{equation}
    \item
      \label{it:PI=1}
      \begin{equation}
        p_t \circ i_s =
        \begin{cases}
          1 \in \exaddcat(x_s,x_t), \text{ if } s=t\\
          0, \text{ otherwise}
        \end{cases}
      \end{equation}
    \end{enumerate}
    Then \(P\) and \(I\) exhibit \(x\)
    as the product \(\prod_{s=1}^nx_n\)
    and as the coproduct \(\coprod_{s=1}^nx_n\),
    respectively.
    Morover, the canonical comparison map
    \eqref{eq:canmap-prod-coprod}
    is the identity \(1\colon x\to x\).
  \item
    \label{it:P-gives-I}
    Assume the product \(x=\prod_{s=1}^n\) exists
    and let
    \(P=(p_s\colon x\to x_s)_{s=1}^n\) be the product cone.
    Then there exists a unique cocone
    \(I=(i_s\colon x_s\to x)_{s=1}^n\)
    satisfying conditions
    \ref{it:IP=1} and \ref{it:PI=1} above.
  \item
    \label{it:I-gives-P}
    Dually, for every coproduct cocone
    \(I=(i_s\colon x_s\to x)_{s=1}^n\)
    there exists a unique cone
    \(P=(p_s\colon x\to x_s)_{s=1}^n\)
    satisfying \ref{it:IP=1} and \ref{it:PI=1}.
  \end{enumerate}
\end{lem}

Since \(\oplus\) and \(+\) determine each other, we have the following corollary:

\begin{cor}
  \label{cor:semiadd-homs-vs-biprod}
  \begin{enumerate}
  \item[]
  \item
    Let \(\exaddcat\) be a category enriched in abelian monoids.
    If \(\exaddcat\) admits finite products (equivalently, finite coproducts)
    then it is semiadditive.
  \item
    Let \(F\colon \exaddcat\to\exaddcat'\) be a functor between additive categories.
    The following are equivalent:
    \begin{enumerate}
    \item
      the functor \(F\) preserves finite products;
    \item
      the functor \(F\) preserves finite coproducts;
    \item
      the functor \(F\) preserves the addition on the hom-sets.
    \end{enumerate}
  \end{enumerate}
\end{cor}

\Cref{lem:prod-coprod-equations} is well known.
However, its proof will serve as a guide for its categorified counterpart,
so we shall explain it here:

\begin{proof}[Proof of \Cref{lem:prod-coprod-equations}]
  To prove \ref{it:PI-prod-coprod}, we assume
  \ref{it:IP=1} and \ref{it:PI=1} and show that \(P\) is a product cone;
  the statement about \(I\) is dual.
  We need to show that for each \(t\in\exaddcat\) the natural map
  \begin{equation}
    P_*\colon\exaddcat(t,x)\to \prod_{s=1}^n\exaddcat(t,x_s);
    \quad f\mapsto (p_s\circ f)_{s=1}^n
  \end{equation}
  is a bijection.
  Using \(I\) we can produce an explicit inverse via the formula
  \begin{equation}
    I_*\colon (f_s)_{s=1}^n\mapsto \sum_{s=1}^ni_s f_s.
  \end{equation}
  It satisfies
  \begin{align}
    (P_*\circ I_*) (f_s)_{s=1}^n
    &= P_* (\sum_{s=1}^n i_s f_s)\\
    &= (p_u\circ\sum_{s=1}^n i_s f_s)_{u=1}^n\\
    &= (\sum_{s=1}^n p_u i_s f_s)_{u=1}^n\\
    &= (f_u)_{u=1}^n
  \end{align}
  (using equation \ref{it:PI=1} in the last step)
  and
  \begin{align}
    (I_*\circ P_*)(f)
    &= I_*(p_s\circ f)_{s=1}^n\\
    &= \sum_{s=1}(i_sp_s f)\\
    &= (\sum_{s=1}i_sp_s)\circ f\\
    &= 1\circ f = f
  \end{align}
  (using equation \ref{it:IP=1} in the last step),
  as desired.
  Moreover, equation \ref{it:PI=1} says precisely that the identity
  \(1\colon x\to x\)
  satisfies the defining equation to be the map
  \eqref{eq:canmap-prod-coprod}.

  Next we prove \ref{it:P-gives-I}; the statement \ref{it:I-gives-P} is dual.
  By the univeral property of the product cone \(P\),
  there are unique maps \(i_s\colon x_s\to x\) satisfying
  equation \ref{it:PI=1}.
  These maps then assemble into the desired cocone \(I\).
  To verify equation \ref{it:IP=1} it suffices to postcompose
  with all the product projections \(p_u\) and compute
  \begin{align}
    p_u\circ \sum_{s=1}^ni_sp_s
    &= \sum_{s=1}^np_ui_sp_s\\
    &= p_u = 1\circ p_u
  \end{align}
  (using equation \ref{it:PI=1} in the second step).
\end{proof}

There is one further aspect of additive categories
whose categorification will be discussed here: matrix calculus. This is based
on the observation that in any category \(\exaddcat\) any map
\begin{equation}
  f\colon \coprod_{s=1}^nx_s\longrightarrow \prod_{t=1}^m y_t
\end{equation}
from a coproduct to a product can be encoded
through the bijection
\begin{equation}
  \exaddcat\left(\coprod_{s=1}^n x_s, \prod_{t=1}^m y_t\right)
  \cong \prod_{t=1}^m\prod_{s=1}^n\exaddcat(x_s,y_t)
\end{equation}
as an \(m\times n\)-matrix \((f_{ts})_{t=1, s=1}^{m,n}\)
whose entry \(f_{ts}\) is a map \(x_s\to y_t\).

The special feature of semiadditive categories
is that it makes sense to consider the composite
\begin{equation}
  h\colon \coprod_{s=1}^nx_s\xrightarrow{f} \prod_{t=1}^m y_t
  \cong
  \coprod_{t=1}^ny_t\xrightarrow{g} \prod_{u=1}^l z_u
\end{equation}
of two such maps, using the identification \eqref{eq:canmap-prod-coprod}.
This composite corresponds to a matrix
\begin{equation}
  (h_{us})_{u=1,s=1}^{l,n}\in\prod_{u=1}^l\prod_{s=1}^n\exaddcat(x_s,z_u).
\end{equation}
It is not hard to verify that
the matrix corresponding to the composite \(h\)
arises from the matrices of \(f\) and \(g\) by the usual rule for matrix multiplication:
\begin{equation}
  \label{eq:usual-atrix-multiplication}
  h_{us} = \sum_{t=1}^m g_{ut}f_{ts}
\end{equation}
From this perspecive, the identification
\eqref{eq:canmap-prod-coprod}
is just the identity matrix which has identities on the diagonal an zeroes everywhere else.

\section{Preliminaries}

Throughout this paper we use the notation ``\(x:A\)''
borrowed from homotopy type theory to say that \(x\) is a
term/inhabitant/element/object of the
(\(\infty\)-)groupoid, (\(\infty\)-)category,
or even \((\infty,2)\)-category \(A\).
When we construct an object ``\(F(x):B\) for each \(x:A\)'',
it is understood that \(F(x)\) is supposed to be functorial in \(x\)
in the appropriate sense.
This allows us to unambiguosly write formulas such as
\(\colim_{x:A}F(x)\) or \((F(x))_{x:A}\),
which of course only make sense with the additional functoriality in mind.

We reserve the notation \(x\in A\) for the case
when \(A\) is discrete, i.e.\ (equivalent to) a set;
in this case, the question of functoriality is vacuous.

\subsection{\texorpdfstring{\((\infty,2)\)}{(infinity,2)}-categories}
\label{subsec:laxinfty2}

In this paper, we think of \((\infty,2)\)-categories
as categories enriched in the \(\infty\)-category
\(\Catinfty\) of \(\infty\)-categories.
For a general treatment of enriched \(\infty\)-categories,
we refer to the work of Gepner and Haugseng~\cite{gepner-haugseng}.
For different approaches to \((\infty,2)\)-categories
we refer to \cite{lurie:2cat} and \cite{GR17I}.

Our goal is not to develop any \((\infty,2)\)-categorical
foundations but rather to develop the theory of lax additivity
while assuming that such foundations are already laid.
In practice, this means that none of our arguments and constructions
are performed explicitly in a model,
but only using the general high-level features which any
theory of \((\infty,2)\)-categories is expected to share.
We treat these ingredients axiomatically:

Let \(\extwocat\) be an \((\infty,2)\)-category.
\begin{itemize}
\item
  It has an underlying \(\infty\)-category \(\extwocatone\),
  and an underlying \(\infty\)-groupoid \(\extwocat^{\simeq}=(\extwocatone)^{\simeq}\).
\item
  It has a hom-functor
  \begin{equation}
    \extwocat(-,-)\colon\extwocat^\op\times \extwocat\to \CCatinfty,
  \end{equation}
  which takes values in the \((\infty,2)\)-category of \(\infty\)-categories.
  Occasionally, it is convenient to consider the hom-functor
  \begin{equation}
    \extwocat(-,-)\colon\extwocatone^\op\times\extwocatone \to \Catinfty
  \end{equation}
  as a functor of the underlying \(\infty\)-categories,
  and its associated Cartesian fibration
  \begin{equation}
    \Twcontra(\extwocat)=\contraGroth{\extwocat(-,-)}
    \to \extwocatone\times\extwocatone^\op.
  \end{equation}
\item
  There are composition functors
  \begin{equation}
    \label{eq:composition-map}
    \extwocat(X,Y)\times\extwocat(Y,Z)\to\extwocat(X,Z),
  \end{equation}
  functorial in \(X,Y,Z:\extwocat^{\simeq}\).
  Composition is coherently associative;
  this is formalized in \cite{gepner-haugseng} by encoding
  the \((\infty,2)\)-category \(\extwocat\)
  as an algebra in the monoidal \(\infty\)-category \((\Catinfty,\times)\)
  of a certain generalized nonsymmetric operad
  \(\Delta^\op_{\extwocat^{\simeq}}\to\Delta^\op\).
\item
  More generally, the composition map \eqref{eq:composition-map}
  is also natural in \(X:\extwocatone^\op\), \(Z:\extwocatone\)
  and dinatural in \(Y:\extwocatone\)
  (and not just in their groupoid cores).
  Thinking in terms of fibrations,
  this means that composition can be written as the dashed functor
  \begin{equation}
    \begin{tikzcd}
      \mixedGroth{X}{Y}{\extwocat(X,Y)}
      \times_{(Y:\extwocatone)}
      \mixedGroth{Y}{Z}{\extwocat(Y,Z)}
      \ar[dashed,r]
      \ar[d]
      &
      \mixedGroth{X}{Z}{\extwocat(X,Z)}
      \ar[d]
      \\
      (X:\extwocatone)
      \times
      (Y:\extwocatone)
      \times
      (Z:\extwocatone)
      \ar[r]
      &
      (X:\extwocatone)
      \times
      (Z:\extwocatone)
    \end{tikzcd}
  \end{equation}
  of mixed (Cartesian, coCartesian) fibrations.
\item
  It makes sense to talk about adjunctions
  \(f\dashv \ra{f}\colon X\to Y\)
  in \(\extwocat\).
  These are characterized by the fact that
  \begin{equation}
    ({f\circ})\dashv({\ra{f}\circ})
    \colon
    \extwocat(T,X)\to \extwocat(T,Y)
    \quad\text{and}\quad
    ({\circ \ra{f}})\dashv({\circ f})
    \colon
    \extwocat(X,T)\to\extwocat(Y,T)
  \end{equation}
  are adjunctions of \(\infty\)-categories for all \(T:\extwocat\).
\end{itemize}

For the purpose of developing the theory of lax additivity
we do not need the full coherent associativity of the composition law,
but only its incoherent shadow.
More precisely,
it suffices to postcompose the enrichment with the symmetric monoidal functor
\((\Catinfty,\times)\to(\ho{\Catinfty},\times)\),
and think of \(\extwocat\) as enriched in the homotopy category of
\(\infty\)-categories up to equivalence.

\subsection{Lax limits and colimits}

We start by recalling the definition of a lax limits and colimits
in a \((\infty,2)\)-category.
Let \(S\) be an \(\infty\)-category.

First, let \(\X\colon S\to \CCatinfty\) be a diagram of \(\infty\)-categories.
Let \(\laxlim\X\) be the \(\infty\)-category of sections
of the (covariant) Grothendieck construction
\(\covGroth{\X} \to S\)
associated to the functor \(\X\).
Informally, objects of \(\laxlim\X\) consist of
\begin{enumerate}
\item
  for each object \(s\) of \(S\), an object \(x_s\) in \(\X_s\),
\item
  for each edge \(f\colon s\to t\) in \(S\),
  a morphism \(x_f\colon\X_f(x_s)\to x_t\) in \(\X_t\),
\item
  for each \(2\)-simplex \(s\xrightarrow{f} t\xrightarrow{g} u\)
  (with the composite \(gf\) implicit) in \(S\)
  a \(2\)-simplex
  \begin{equation}
    \begin{tikzcd}
      &\X_{g}(x_t)\ar[rd,"x_g"]\\
      \X_{gf}(x_s)\ar[ru,"\X_g(x_f)"]\ar[rr,"x_{gf}"]&&x_u
    \end{tikzcd}
  \end{equation}
  in \(\X_u\),
\item
  and so on for higher simplices of \(S\).
\end{enumerate}
We will denote an object of \(\laxlim_{s:S}\X\)
as a tuple \((x_s)_{s:S}\).

Now, let \(\extwocat\) be an arbitrary \((\infty,2)\)-category
and \(\X\colon S\to \extwocat\) a diagram.
A lax cone over the diagram \(X\) with vertex \(L\) is
an object \((\phi_s)_{s:S}:\laxlim_{s:S}\extwocat(L,\X_s)\),
where
\(s\mapsto\extwocat(L,\X_s)\)
is the \(S\)-shaped diagram in \(\Catinfty\) obtained as the composite
\begin{equation}
  S\xrightarrow{\X}\extwocat\xrightarrow{\extwocat(L,-)}\Catinfty
\end{equation}
Unpacking the above, one sees that such a cone consists of
\begin{enumerate}
\item
  for each object \(s:S\), a structure map \(\phi_s\colon L\to \X_s\),
\item
  for each arrow \(f\colon s\to t\) in \(S\),
  a lax cone
  \begin{equation}
    \laxcone{\X_s}{L}{\X_t}{\phi_s}{\phi_t}{\X_f}
  \end{equation}
  i.e.\ a map \(\phi_f\colon \X_f\phi_s\to \phi_t\)
  in \(\extwocat(L,\X_t)\),
\item
  together with coherent pasting identifications,
  \(\phi_g\circ\X_g\phi_f\simeq \phi_{gf}\)
  for composable arrows \(s\xrightarrow{f}t\xrightarrow{g}u\) in \(S\).
\end{enumerate}

For each other object \(L':\extwocat\)
we have a canonical composition map
\begin{equation}
  \begin{tikzcd}
    \extwocat(L',L)\times\laxlim_{s:S}\extwocat(L,\X_s)
    \ar[r,"{-\circ -}"]
    \ar[d,"\Delta\times \id"]
    &
    \laxlim_{s:S}\extwocat(L',\X_s)
    \\
    \laxlim_{s:S}\extwocat(L',L)
    \times
    \laxlim_{s:S}\extwocat(L,\X_s)
    \ar[d,"\simeq"]
    \\
    \laxlim_{s:S}\extwocat(L',L)\times\extwocat(L,\X_s)
    \ar[uur,bend right]
  \end{tikzcd}
\end{equation}
which explicitly sends a cone \(\Phi=(\phi_s)_{s:S}\) with vertex \(L\)
and a morphism \(l\colon L'\to L\)
to the cone \(\Phi\circ l=(\phi_s\circ l)_{s:S}\).

Dually we can define the \(\infty\)-category
of lax cocones on \(\X\) with vertex \(L\)
as \(\laxlim_{s:S^\op}\extwocat(\X_s,L)\).
Explicitly, such a cocone \((\psi_s)_{s:S}\) has
structure maps \(\psi_s\colon \X_s\to L\)
and lax triangles
\begin{equation}
  \laxcocone[\phi_f]{\X_s}{L}{\X_t}{\psi_s}{\psi_t}{\X_f}
\end{equation}
over each arrow \(f\colon s\to t\) of \(S\).

\begin{defi}
  \label{defi:lax-limit}
  A cone \(P=(p_s)_{s:S}: \laxlim_{s:S}\extwocat(L,\X_s)\)
  is called a \emph{lax limit cone}
  if for each object \(L':\extwocat\) the functor
  \begin{equation}
    {P\circ-}\colon \extwocat(L',L)\longrightarrow\laxlim_{s:S}\extwocat(L',\X_s);
    \quad
    F\mapsto (p_s\circ F)_{s:S}
  \end{equation}
  is an equivalence of \(\infty\)-categories;
  in this case we call the object \(L\) a lax limit of the diagram
  \(\X\colon S\to \extwocat\)
  and write \(L=\laxlim_{s:S}\X_s\).

  Dually, we say that a cocone
  \(I=(i_s)_{s:S^\op}:\laxlim_{s:S^\op}\extwocat(\X_s,L)\)
  is a \emph{lax colimit cone}
  if for each \(L':\extwocat\) the functor
  \begin{equation}
    {-\circ I}\colon \extwocat(L,L')\longrightarrow\laxlim_{s:S^\op}\extwocat(\X_s,L');
    \quad
    F\mapsto (F\circ i_s)_{s:S^\op}
  \end{equation}
  is an equivalence;
  in this case we call \(L\) a lax colimit of \(\X\)
  and write \(L=\laxcolim_{s:S}\X_s\).
\end{defi}

\begin{rem}
  Our definition starts by \emph{defining} the lax limits of \(\infty\)-categories
  to be sections of the Grothendieck construction
  and then defining lax limits and lax colimits in arbitrary \((\infty,2)\)-categories
  by considering (co)representables.
  One can also define lax limits and colimits as a special case of
  \emph{weighted colimits},
  which can be defined directly in terms of ordinary limits/colimits.
  When using the latter definition, one can then \emph{compute}
  that lax limits of \(\infty\)-categories as sections of the Grothendieck construction,
  see \cite[Proposition~7.1 and Corollary~7.7]{ghn}.
\end{rem}

\begin{exa}
  \label{ex:lax-lim-cat}
  Let \(\extwocat=\CCatinfty\) be the \((\infty,2)\)-category
  of \(\infty\)-categories.
  Let \(\X\colon S\to\CCatinfty\) be a diagram.
  \begin{enumerate}
  \item
    As the notation suggests,
    the  lax limit of \(\X\) is the \(\infty\)-category
    \(L\coloneqq \laxlim_{s:S}\X_s\coloneqq \Fun_{S}(S,\covGroth[s:S]{\X_s})\)
    of sections of the corresponding Grothendieck construction.
    Indeed, naturally in  \(L':\CCatinfty\) we have the equivalence
    \begin{align}
      \CCatinfty(L',L)
      &=
      \Fun(L',\Fun_S(S,\covGroth{\X}))
      \\
      &\simeq
      \Fun_S(S\times L',\covGroth{\X})
      \\
      &\simeq
      \Fun_S(S,\covGroth[s:S]\Fun(L',\X_s))
      \\
      &=
      \laxlim_{s:S}(\CCatinfty(L',\X_s)),
    \end{align}
    which is induced by composition with the canonical cone
    \begin{equation}
      \label{eq:proj-functors-laxlim}
      P=
      \left(
        p_s\colon L=\laxlim_{S}\X \to \X_s
      \right)_{s:S}
    \end{equation}
    given by evaluation of sections.
  \item
    The lax colimit of the diagram \(\X\)
    is the contravariant Grothendieck construction
    \begin{equation}
      \laxcolim_{s:S}\X_s=
      \contraGroth[s:S]{\X_s},
    \end{equation}
    exhibited by the canonical cocone
    \begin{equation}
      I=
      \left(i_s\colon \X_s\to\contraGroth{\X}\right)_{s:S^\op}
    \end{equation}
    that includes the individual fibers.
  \end{enumerate}

  Assume now that the diagram \(\X\) takes values in \emph{stable}
  \(\infty\)-categories,
  \begin{enumerate}[resume]
  \item
    \label{it:ex:lax-limit-stable-cats}
    The \(\infty\)-category
    \(\laxlim_{s:S}\X_s=\Fun_S(S,\covGroth{\X})\) is again stable
    because limits and colimits of sections are computed pointwise.
    For the same reason, every functor
    \(F\colon L'\to\laxlim_{s:S}\X\)
    is exact if and only if each composite
    \(p_s\circ F\) is exact.
    It follows that the cone \(P\) exhibits the \(\infty\)-category
    \(\laxlim_{s:S}\X_s\) also as a lax limit in the \((\infty,2)\)-category
    of stable \(\infty\)-categories and exact functors.
  \item
    The \(\infty\)-category \(\contraGroth{\X}\),
    which is the lax colimit of \(\X\) in \(\CCatinfty\),
    is typically \emph{not} stable;
    to compute the lax colimit of \(\X\)
    in the \((\infty,2)\)-category of stable \(\infty\)-categories
    one therefore has to stabilize this \(\infty\)-category,
    which is a rather tricky operation.
    However, it will follow from the theory of lax matrices
    that---as long as the stable categories in question have enough colimits,
    for example because \(S\) is finite or because \(\X\) takes vales in \emph{presentable} stable \(\infty\)-categories---%
    this lax colimit indeed just agrees with the lax \emph{limit}
    which can be computed in \(\CCatinfty\);
    see \Cref{cor:lax-co-limits-are-bilimits} below.
  \end{enumerate}
\end{exa}

\begin{rem}
  \Cref{defi:lax-limit} can easily be modified to also define
  \emph{partially lax} limits and colimits (sometimes also called \emph{marked (co)limits}):
  If the indexing category \(S\) is equipped with some collection \(M\) of marked arrows,
  we can define the \(M\)-partially lax limit of a diagram
  \(\X\colon S \to \CCatinfty\)
  as the \(\infty\)-category of those sections of the Grothendieck construction
  \(\covGroth\X\to S\),
  whose value on the arrows in \(M\) is cocartesian.
  Then one defines partially lax limits and colimits in an arbitrary
  \((\infty,2)\)-category \(\extwocat\)
  analogously to \Cref{defi:lax-limit}.
  See also \cite{berman2020} and \cite{abellan-marked}
  for a general treatment of partially (op)lax (co)limits.

  In this paper, we mostly deal with (fully) lax limits or colimits
  (i.e. \(M\) only consists of the equivalences of \(S\)).
\end{rem}

\begin{exa}
  \label{exa:directed-pull-push}
  The only partially lax limits we need in this paper are the
  directed pullback and directed pushout,
  which we denote by
  \(\laxpull{\A}{\B}{\C}\) and \(\laxpush{\A}{\B}{\C}\).
  Abstractly, they are equipped with the universal squares
  \begin{equation}
    \oplaxsquare
    {\laxpull{\A}{\B}{\C}}
    {\A}{\C} {\B}
    {}{}
    {f}{g} {}
    \quad \text{and} \quad
    \oplaxsquare
    {\B}
    {\A}{\C}
    {\laxpush{\A}{\B}{\C}}
    {f}{g}{}{}{}
  \end{equation}
  inhabited by a (possibly noninvertible) \(2\)-morphism
  (for given \(\A,\B,\C\) and \(f,g\) which we omit from the notation).
  Concretely, they can be defined as partially lax limits/colimits
  with the arrow indexing \(g\) being marked
  or,
  equivalently,
  as partially oplax limits/colimits
  with the arrow indexing \(f\) being marked.
\end{exa}

\section{Lax matrices}
\label{subsec:lax-matrices}

Throughout this section,
let \(\extwocat\) be an \((\infty,2)\)-category
enriched in \(\infty\)-categories with colimits,
i.e.,
\begin{itemize}
\item
  each hom-category \(\extwocat(X,Y)\) has all colimits and
\item
  and each composition functor
  \(\extwocat(X,Y)\times\extwocat(Y,Z)\to\extwocat(X,Z)\)
  preserves colimits in each variable separately.
\end{itemize}

Analogously to the case of ordinary coproducts and products
(which corresponds to the case where the category \(S\) is just a set),
we can interpret maps from a lax colimit to a lax limit as a sort of matrices:
By the defining property we have
\begin{align}
  \extwocat(\laxcolim_{s:S}\X_s,\laxlim_{t:T}\Y_t)
  &\simeq\laxlim_{t:T}\extwocat(\laxcolim_{s:S}\X_s,\Y_t)\\
  &\simeq\laxlim_{t:T}\laxlim_{s:S^\op}\extwocat(\X_s,\Y_t)\\
  &\simeq\laxlim_{(t,s):T\times S^\op}\extwocat(\X_s,\Y_t)
\end{align}
so that we can interpret a map
\(\alpha\colon \laxcolim_{s:S}\X_s\to\laxlim_{t:T}\Y_t\)
as a tuple
\((\alpha_{t,s})_{(t,s):T\times S^\op}\)
which we think of as a matrix whose rows are indexed by \(T\)
and whose columns are indexed by \(S^\op\).
We define
\begin{equation}
  \laxmat{\extwocat}(\X,\Y)\coloneqq
  \laxlim_{(t,s):T\times S^\op}\extwocat(\X_s,\Y_t).
\end{equation}
Note that this is a well defined \(\infty\)-category
even when
\(\laxcolim\X\) and/or \(\laxlim\Y\) does not exist.
When \(\X\colon\singleton\to\extwocat\) is just an object \(X=\X_{*}\),
we still use the notation
\(\laxmat{\extwocat}(X,\Y)=\laxmat{\extwocat}(\X,\Y)\)
and observe that it is precisely the \(\infty\)-category
of lax cones on \(\Y\) with vertex \(X\);
and analogously for lax cocones

\begin{exa}
  \label{exa:2x2-matrices-intro}
  Let \(S=T=\Delta^1=\set{0\xrightarrow{f} 1}\) be the walking arrow
  and consider two diagrams
  \(\X\colon S\to \extwocat\) and \(\Y\colon T\to \extwocat\).
  Then we can compactly describe
  objects of
  \begin{equation}
    \laxmat{\extwocat}(\X,\Y)
    = \laxlim\left(
      \begin{tikzcd}
        \extwocat(\X_0,\Y_0)
        \ar[from=r]
        \ar[d]
        &
        \extwocat(\X_1,\Y_0)
        \ar[d]
        \\
        \extwocat(\X_0,\Y_1)
        \ar[from=r]
        &
        \extwocat(\X_1,\Y_1)
      \end{tikzcd}
    \right)
  \end{equation}
  as \(T\times S^\op\)-indexed diagrams in the Grothendieck construction,
  which we depict as follows:
  \begin{equation}
    \left(
      \begin{tikzcd}
        \alpha_{00}
        \ar[d]
        &
        \alpha_{01}
        \ar[l]
        \ar[d]
        \\
        \alpha_{10}
        &
        \alpha_{11}
        \ar[l]
      \end{tikzcd}
    \right)
  \end{equation}
  Explicitly unpacking this notation, such a matrix consists of:
  \begin{enumerate}
  \item
    four 1-morphisms 
    \begin{equation}
      \begin{matrix}
        \alpha_{00}\colon \X_0\to \Y_0
        &
        \alpha_{01}\colon \X_1\to \Y_0
        \\
        \alpha_{10}\colon \X_0\to \Y_1
        &
        \alpha_{11}\colon \X_1\to \Y_1
      \end{matrix}
    \end{equation}
  \item
    four 2-morphisms
    \begin{equation}
      \begin{tikzcd}
        \Y_f\circ\alpha_{00}\ar[d,"\alpha_{f0}"]
        &\alpha_{00}\ar[from=r,"\alpha_{0f}"]&\alpha_{01}\circ \X_f
        \\
        \alpha_{10}
        &&
        \Y_f\circ\alpha_{01}\ar[d,"\alpha_{f1}"]
        \\
        \alpha_{10}\ar[from=r,"\alpha_{1f}"]&\alpha_{11}\circ \X_f
        &
        \alpha_{11}
        \\
      \end{tikzcd}
    \end{equation}
  \item
    assembling into a commutative square
    \begin{equation}
      \begin{tikzcd}[column sep = large]
        \Y_f\circ\alpha_{00}
        \ar[d,"\alpha_{f1}"]
        &
        \Y_f\circ\alpha_{01}\circ \X_f
        \ar[l,"{\Y_f\circ \alpha_{0f}}"']
        \ar[d,"\alpha_{f1}\circ\X_f"]
        \\
        \alpha_{10}
        \ar[from=r,"\alpha_{1f}"']
        &
        \alpha_{11}\circ \X_f
      \end{tikzcd}
    \end{equation}
  \end{enumerate}
\end{exa}

We now introduce the lax matrix multiplication
which categorifies the classical formula
\eqref{eq:usual-atrix-multiplication}.
The classical formula involves
a finite sum of elements in some hom-set \(\exaddcat(x_s,z_u)\)
of the category \(\exaddcat\).
Our categorified analog of these sums will be categorical colimits.

\begin{con}
  \label{con:row-column-multiplication}
  Let \(S\) be an \(\infty\)-category
  and \(\X\colon S\to \extwocat\) a diagram.
  Passing to the Cartesian fibrations classifying
  the respecive hom-functors,
  we obtain a commutative square
  \begin{equation}
    \begin{tikzcd}
      \Twcontra(S)=
      \contraGroth {S(-,-)}
      \ar[d,"p"]
      \ar[r,"\alpha"]
      &
      \contraGroth {\extwocat(-,-)}=\Twcontra(\extwocat)
      \ar[d,"q"]
      \\
      S\times S^\op
      \ar[r,"\X\times \X^\op"]
      &
      \extwocatone\times \extwocatone^\op
    \end{tikzcd}
  \end{equation}
  which amounts to the dashed section
  \begin{equation}
    \begin{tikzcd}[column sep=large]
      &
      \contraGroth[(s,t)]{\extwocat(\X_s,\X_t)}
      \isCartesian
      \ar[r]
      \ar[d]
      &
      \Twcontra(\extwocat)
      \ar[d,"q"]
      \\
      \Twcontra(S)\ar[r,"p"]
      \ar[ur,dashed,"{\alpha}",bend left=20]
      &
      S\times S^\op
      \ar[r,"{\X\times\X^\op}"]
      &
      \extwocatone\times{\extwocatone^\op}
    \end{tikzcd}
  \end{equation}
  of the pullback-fibration \((\X\times\X^\op)^*(q)\)
  which informally sends an arrow
  \((f\colon s\to t):\Twcontra(S)\)
  to \(\X_f:\extwocat(\X_s,\X_t)\).

  We can now construct the composite functor
  \begin{align}
    &
      \laxlim_{s:S} \extwocat(L,\X_s)
      \times 
      \laxlim_{t:S^\op} \extwocat(\X_t,L')
    \\
    &
      =
      \Fun_{S\times S^\op}
      \left(S\times S^\op,
      \covGroth[(s,t)]
      {\extwocat(L,\X_s) \times\extwocat(\X_t,L')}
      \right)
    \\
    &
      \longrightarrow
      \Fun_{S\times S^\op}
      \left(\Twcontra(S),
      \covGroth[(s,t)]
      {\extwocat(L,\X_s) \times\extwocat(\X_t,L')}
      \right)
    \\
    &
      \longrightarrow
      \Fun_{S\times S^\op}
      \left(\Twcontra(S),
      \contraGroth[(s,t)]
      {\extwocat(\X_s,\X_t)}
      \times_{S\times S^\op}
      \covGroth[(s,t)]
      {\extwocat(L,\X_s) \times\extwocat(\X_t,L')}
      \right)
    \\
    &
      =
      \Fun_{S\times S^\op}
      \left(\Twcontra(S),
      \covGroth[s]
      {\extwocat(L,\X_s)}
      \times_{S}
      \contraGroth[(s,t)]
      {\extwocat(\X_s,\X_t)}
      \times_{S^\op}
      \covGroth[t]
      {\extwocat(\X_t,L')}
      \right)
    \\
    &
      \longrightarrow
      \Fun\left(\Twcontra(S), \extwocat(L,L')\right)
      \xrightarrow{\colim}
      \extwocat(L,L'),
  \end{align} 
  where
  \begin{itemize}
  \item
    the first arrow is pullback of sections along
    \(p\colon\Twcontra(S)\to S\times S^\op\),
  \item
    the second arrow adds the section \(\alpha\)
    in the first component of the fiber product,
  \item
    the third arrow is given by composition with the composition map
    \begin{equation}
      \covGroth[s]
      {\extwocat(L,\X_s)}
      \times_{S}
      \contraGroth[(s,t)]
      {\extwocat(\X_s,\X_t)}
      \times_{S^\op}
      \covGroth[t]
      {\extwocat(\X_t,L')}
      \longrightarrow
      \extwocat(L,L'),
    \end{equation}
  \item
    the last arrow is just the colimit functor
    in the \(\infty\)-category \(\extwocat(L,L')\).
  \end{itemize}
  On objects, this functor
  takes a lax cone and a lax cocone on \(\X\),
  \begin{equation}
    \Phi: \laxlim_{s:S}\extwocat(L,\X_s)
    \quad\text{and}\quad
    \Psi: \laxlim_{t:S^\op}\extwocat(\X_t,L'),
  \end{equation}
  and sends them to the map \(\Psi\circ_S \Phi\colon L\to L'\)
  defined by the formula
  \begin{equation}
    \label{eq:row-column-multiplication}
    {\Psi\circ_S \Phi}\coloneqq \colim_{(f:s\to t):\Twcontra(S)}
    (\Psi_{t}\circ\X_f\circ\Phi_s).
  \end{equation}
\end{con}

\begin{rem}
  \label{rem:set-matrix-mult}
  When \(S=\set{1,2,\dots, n}\) is a finite set
  the (Cartesian) twisted arrow category \(\Twcontra(S)\to S\times S\)
  can be canonically identified with the
  diagonal \(\Delta\colon S\to S\times S\).
  Under this identification the formula
  \eqref{eq:row-column-multiplication}
  simplifies to
  \begin{equation}
    (\Psi\circ_S\Phi)
    =\coprod_{\substack{s,t\in S\\s=t}} \Psi_t\circ\id\circ\Phi_s 
    =\coprod_{s=1}^n\Psi_s\circ\Phi_s
  \end{equation}
  which is just the usual multiplication
  \eqref{eq:usual-atrix-multiplication}
  of the row vector \(\Psi\)
  with the column vector \(\Phi\).

  When \(S=\singleton\) is even a singleton,
  this formula just returns the original composition in the \((\infty,2)\)-category
  \(\extwocat\).
  For this reason we drop the subscript \(S\)
  and just write \(-{\circ}-\) instead of \(-{\circ_{\singleton}}-\).
\end{rem}

We assemble our categorified analog of row-column multiplication
to the lax version of matrix multiplication:

\begin{con}
  \label{con:lax-matrix-multiplication}
  Let \(S,U\) be \(\infty\)-categories,
  and \(\X\colon S\to \extwocat\)
  and \(\Z\colon U\to \extwocat\) two diagrams.
  For each object \(Y:\extwocat\) we consider the functor
  \begin{align}
    &
      \laxlim_{s:S^\op}\extwocat(\X_s,Y)
      \times
      \laxlim_{u:U}\extwocat(Y,\Z_u)
    \\
    &
      =
      \laxlim_{(u,s):U\times S^\op}
      \extwocat(\X_s,Y)\times \extwocat(Y,\Z_u)
    \\
    &
      \longrightarrow
      \laxlim_{(u,s):U\times S^\op}
      \extwocat(\X_s,\Z_u),
  \end{align}
  induced by composition of \(\extwocat\).
  On objects it takes a lax cocone and a lax cone,
  \begin{equation}
    \Psi: \laxlim_{s:S^\op}\extwocat(\X_s,Y)
    \quad\text{and}\quad
    \Phi: \laxlim_{u:S}\extwocat(Y,\Z_u),
  \end{equation}
  and sends them to the matrix
  \(\Phi\circ\Psi : \laxlim_{(u,s):U\times S^\op}\)
  described by the formula
  \begin{equation}
    \label{eq:matrix-mult-over-point}
    (\Phi\circ\Psi)_{us}=\Phi_u\circ\Phi_s.
  \end{equation}
  More generally,
  we can replace the object \(Y:\extwocat\) by a diagram
  \(\Y\colon T\to\extwocat\)
  and consider the functor
  \begin{align}
    &\laxmat{\extwocat}(\X,\Y)\times\laxmat{\extwocat}(\Y,\Z)
    \\
    &
      =
      \laxlim_{(u,s):U\times S^\op}
      \left(
      \laxlim_{t:T} \extwocat(\X_s,\Y_t)
      \times
      \laxlim_{t':T^\op}\extwocat(\Y_{t'},\Z_u)
      \right)
    \\
    &
      \xrightarrow{\laxlim_{u,s}(-\circ_T-)}
      \laxlim_{(u,s):U\times S^\op}
      \extwocat(\X_s,\Z_u)
     =
    \laxmat{\extwocat}(\X,\Z)
  \end{align}
  which is given in componentwise in \(u,s\)
  by the composition functor
  from \Cref{con:row-column-multiplication}
  (applied to lax cones and cocones on \(\Y\)).
  Explicitly, this functor is given by the formula
  \begin{equation}
    \label{eq:matrix-mult}
    \Phi\circ_T\Psi
    =
    \left(
      \colim_{(f:t\to t'):\Twcontra(T)}
      (\Phi_{ut'}\circ\Y_{f}\circ\Psi_{ts})
    \right)_{(u,s):U\times S^\op}.
  \end{equation}
  This is what we call the lax matrix multiplication.

  Finally, we can assemble all lax matrices of different shapes
  into a category \(\holaxmat{\extwocat}\):
  \begin{itemize}
  \item
    Objects are equivalence classes of diagrams
    \(\X\colon S\to \extwocat\),
    where \(S\) is any small \(\infty\)-category.
  \item
    Morphisms from \(\X\colon S\to \extwocat\) to \(\Y\colon T\to\extwocat\)
    are equivalence classes of matrices
    \(\Phi : \laxmat{\extwocat}(\X,\Y)\).
  \item
    Composition is given by lax matrix multiplication
    of \Cref{con:lax-matrix-multiplication}.
  \end{itemize}
\end{con}

\begin{rem}
  Similarly to \Cref{rem:set-matrix-mult},
  we drop the subscript \(T\) in the case where \(T=\singleton\) is just a point.
\end{rem}

Note that the lax matrix multiplication is
functorial by construction,
making it in particular well defined on equivalence classes.
To prove that \(\holaxmat{\extwocat}\) is indeed a category,
we will thus only need to construct the identity matrix
and prove that lax matrix multiplication is associative up to equivalence.

In fact, we shall prove a slightly stronger statement:
\begin{lem}
  \label{lem:lax-matrix-associative}
  The lax matrix multiplication of \Cref{con:lax-matrix-multiplication} is
  \begin{enumerate}
  \item
    \label{it:lax-matrix-associative}
    associative up to natural equivalence,
    i.e.
    for diagrams
    \(\W\colon R\to\extwocat\)
    and \(\X\colon S\to\extwocat\)
    and \(\Y\colon T\to\extwocat\)
    and \(\Z\colon U\to \extwocat\)
    we have
    \begin{equation}
      {(-\circ_T-)\circ_S -}
      \simeq
      {-\circ_T(-\circ_S-)}
    \end{equation}
    as functors
    \begin{equation}
      \laxmat{\extwocat}(\W,\X)\times\laxmat{\extwocat}(\X,\Y)\times\laxmat{\extwocat}(\Y,\Z)
      \longrightarrow\laxmat{\extwocat}(\W,\Z)
    \end{equation}
  \item
    \label{it:lax-matrix-unital}
    unital up to natural equivalence,
    i.e.,
    for each diagram \(\X\colon S\to\extwocat\)
    there is a matrix
    \(\unitmat^\X:\laxmat{\extwocat}(\X,\X)\)
    with components
    \begin{equation}
      \label{eq:pointwise-formula-I-matrix}
      \unitmat^\X_{ts}=\colim_{f:S(s,t)}\X_f : \extwocat(\X_s,\X_t)
    \end{equation}
    such that
    \begin{equation}
      {\unitmat^\X\circ_S-}\simeq \id
      \quad\text{and}\quad
      {-\circ_S\unitmat^\X}\simeq \id
    \end{equation}
    as endofunctors of
    \(\laxmat{\extwocat}(\Y,\X) \)
    and
    \(\laxmat{\extwocat}(\X,\Y), \)
    respectively
    (for each other diagram \(\Y\colon U\to\extwocat\)).
  \end{enumerate}
\end{lem}

\begin{rem}
  The category
  \(\holaxmat{\extwocat}\) is of course only
  the truncation of an \((\infty,2)\)-category of lax matrices,
  which one could construct with more effort.
  Even \Cref{lem:lax-matrix-associative}
  only shows that lax matrix multiplication
  is associative up to natural equivalence,
  but does not exhibit any sort of coherence
  such as the pentagon.
  We will not be needing this additional layer of coherence in this article
  so \Cref{lem:lax-matrix-associative} will suffice.
\end{rem}

The proof of associativity is relatively straightforward:
\begin{proof}
  [Proof of \Cref{lem:lax-matrix-associative}, part~\ref{it:lax-matrix-associative}]
  For matrices
  \begin{equation}
    F:\laxlim_{S\times R^\op}(\W_r,\X_s),
    \quad
    G:\laxlim_{T\times S^\op}(\X_s,\Y_t),
    \quad
    H:\laxlim_{U\times T^\op}(Y_t,Z_u)
  \end{equation}
  we compute
  \begin{align}
    ((H\circ_TG)\circ_SF)_{ur}
    &\simeq\colim_{f\colon s\to s'}(H\circ_TG)_{us'}\X_f F_{sr}
    \\
    &\simeq\colim_{f\colon s\to s'}(\colim_{g\colon t\to t'}H_{ut'}\Y_gG_{ts'})\X_f F_{sr}
    \\
    &=\colim_{f\colon s\to s'}\colim_{g\colon t\to t'}(H_{ut'}\Y_gG_{ts'}\X_f F_{sr})
    \\
    &=\colim_{(f,g):\Twcontra(S)\times\Twcontra(T)}(H_{ut'}\Y_gG_{ts'}\X_f F_{sr})
    \\
    &\simeq\cdots\simeq (H\circ_T(G\circ_SF))_{ur}
  \end{align}
  naturally in \(F,G,H\) and \(u:U, r:R^\op\);
  where the third step uses that composition in \(\extwocat\) preserves
  colimits in each variable.
\end{proof}

Before we can prove part \ref{it:lax-matrix-unital}
of \Cref{lem:lax-matrix-associative}
we need to construct the unit matrices
\(\unitmat^\X:\laxmat{\extwocat}(\X,\X)\).

\begin{con}
  Consider the commutative square
  \begin{equation}
    \begin{tikzcd}
      \Twcov(S)=
      \covGroth {S(-,-)}
      \ar[d,"p"]
      \ar[r,"\alpha"]
      &
      \covGroth {\extwocat(-,-)}=\Twcov(\extwocat)
      \ar[d,"q"]
      \\
      S^\op\times S
      \ar[r,"\X^\op\times \X"]
      &
      \extwocatone^\op\times \extwocatone
    \end{tikzcd}
  \end{equation}
  induced by a diagram \(\X\colon S\to \extwocat\).
  Here the vertical maps are the coCartesian fibrations classifying
  the respective hom-functors of \(S\) and \(\extwocat\).
  The fact that the \((\infty,2)\)-category \(\extwocat\)
  is enriched in \(\infty\)-categories with colimits
  means that there exists an (essentially unique) left \(q\)-Kan extension of \(\alpha\)
  along \(p\),
  giving rise to the dashed lift
  \begin{equation}
    \begin{tikzcd}
      \Twcov(S)
      \ar[d,"p"]
      \ar[r,"\alpha"]
      &
      \Twcov(\extwocat)
      \ar[d,"q"]
      \\
      S^\op\times S
      \ar[ur,"\unitmat'",dashed]
      \ar[r,"\X^\op\times \X"]
      &
      \extwocatone^\op\times \extwocatone
    \end{tikzcd}
  \end{equation}
  Since the pullback of the coCartesian fibration \(q\) along
  \(\X^\op\times \X\) is, by definition,
  the coCartesian fibration
  \(\covGroth[(s,t):S^\op\times S]{\extwocat(\X_s,\X_t)}\to S^\op\times S\),
  this lift \(\unitmat'\) corresponds to a section of this fibration,
  i.e.,
  an object
  \begin{equation}
    \unitmat^\X=\unitmat :
    \laxlim_{(s,t):S^\op\times S}\extwocat(\X_s,\X_t)
  \end{equation}
  By the pointwise formula,
  we can explicitly compute the value of \(\unitmat^\X\) at \((s,t):S^\op\times S\)
  as the colimit of the composite
  \begin{equation}
    \Twcov(S)/(s,t)\xrightarrow{\alpha} \Twcov(\extwocat)/(\X_s,\X_t)\to\extwocat(\X_s,\X_t),
  \end{equation}
  which is the functor that informally maps
  \begin{equation}
    (f'\colon s'\to t', g\colon s\to s', h\colon t'\to t)\mapsto \X_h\circ \X_{f'}\circ \X_s.
  \end{equation}
  Since the inclusion
  \(S(s,t)\simeq \Twcov(S)_{(s,t)}\hra \Twcov(S)/(s,t)\)
  has a left adjoint (because \(q\) is a coCartesian fibration),
  it is homotopy terminal;
  thus we can compute the components \(\unitmat^\X_{ts}\)
  via the desired explicit formula \eqref{eq:pointwise-formula-I-matrix}.
\end{con}

\begin{rem}
  Since all pointwise colimits 
  \ref{eq:pointwise-formula-I-matrix}
  are taken over \emph{spaces} \(S(s,t)\)
  (as opposed to arbitrary \(\infty\)-categories),
  we see that for the construction of the unit matrix 
  we could have relaxed our assumption on \(\extwocat\)
  and only required it to be enriched in \(\infty\)-categories
  with groupoidal colimits.
\end{rem}

\begin{rem}
  When \(S\) is a set this formula simplifies to
  \begin{equation}
    \unitmat_{ts}=
    \colim_{f:S(s,t)}(\X_f)=
    \begin{cases}
      \id_{\X_s}, \text{ if } s=t\\
      \initialobj, \text{ if } s\neq t
    \end{cases}
  \end{equation}
  which is the direct analog of formula \eqref{eq:cmap-coprod-prod},
  with the initial object \(\initialobj\) of \(\extwocat(x_s,x_t)\)
  taking the role of the zero object of a commutative monoid.
\end{rem}

\begin{exa}
  \label{exa:2x2-matrices-unit-matrix}
  Continuing \Cref{exa:2x2-matrices-intro}, we consider a diagram
  \(\X\colon \Delta^1=\set{0\xrightarrow{10}1}\to \extwocat\).
  The unit \(\Delta^1\times(\Delta^1)^\op\)-matrix then looks as follows:
  \begin{equation}
    \label{eq:2x2-identity-matrix}
    \Donematrix{\id_{\X_0}}{\initialobj}{\X_{10}}{\id_{\X_1}}
    \colon
    \laxmat{\extwocat}(\X,\X)
  \end{equation}
  since the indexing space of the colimit
  \(\colim_{f:S(s,t)}\X_{f} \)
  is either empty in the case \(s=1, t=0\)
  or a singleton otherwise.
\end{exa}

We now finish the proof of \Cref{lem:lax-matrix-associative}.

\begin{proof}
  [Proof of \Cref{lem:lax-matrix-associative}, part~\ref{it:lax-matrix-unital}]
  We only treat the case of postcomposition with \(\unitmat^\X\);
  the other statement is dual.
  We need to show that for every diagram \(\Y\colon R\to \extwocat\),
  the functor
  \begin{equation}
    \unitmat^\X\circ_S-\colon
    \laxmat{\extwocat}(\Y,\X)\to\laxmat{\extwocat}(\Y,\X)
  \end{equation}
  is naturally equivalent to the identity.
  Naturally in \(u:S\), \(r:R^\op\) and  \(F:\laxmat{\extwocat}(\Y,\X)\),
  we compute (in the \(\infty\)-category \(\extwocat(\Y_r,\X_s)\)):
  \begin{align}
    (\unitmat^\X\circ_SF)_{ur}
    &\simeq
    \colim_{(f\colon s\to t):\Twcontra(S)}\left(\colim_{g:S(t,u)}\X_g\right)\X_fF_{sr}
      \\
    &=
      \label{eq:double-colimit-Tw}
      \colim_{\substack{(f\colon s\to t):\Twcontra(S)\\g:S(t,u)}}\X_{gf}F_{sr}
  \end{align}
  where the shape of the second colimit is the category
  \begin{equation}
    \covGroth[(f\colon s\to t):\Twcontra(S)]{S(t,u)}
    = \Twcontra(S)\times_{S^\op}(S/u)^\op
  \end{equation}

  Note that the diagram \((f,g)\mapsto \X_{gf}F_{sr}\)
  over which we are taking the colimit
  arises as the pullback of the diagram
  \begin{equation}
    S/u \to \extwocat(\Y_r,\X_u),
    \quad
    (h\colon s\to u) \mapsto \X_hF_{sr}
  \end{equation}
  along the functor
  \begin{equation}
    \gamma\colon\Twcontra(S)\times_{S^\op}(S/u)^\op\to S/u,
    \quad
    (f\colon s\to t, g\colon t\to u) \mapsto (gf\colon s\to u).
  \end{equation}
  This functor \(\gamma\) has a left adjoint
  \begin{equation}
    S/u \to\Twcontra(S)\times_{S^\op}(S/u)^\op,
    \quad
    (f\colon s\to u) \mapsto (f\colon s\to u, \id_u\colon u\to u),
  \end{equation}
  and is therefore homotopy terminal.
  This allows us to finish the computation:
  \begin{align}
    (\unitmat^\X\circ_SF)_{ur}
    &\simeq
    \label{eq:double-colimit-Tw2}
    \colim_{\substack{(f\colon s\to t):\Twcontra(S)\\g:S(t,u)}}\X_{gf}F_{sr}
    \\
    &\simeq
    \colim_{(h\colon s\to u):S/u}\X_hF_{sr}
    \\
    &\simeq \X_{\id_u}F_{ur}=F_{ur},
  \end{align}
  using in the last step that \(\id_u\colon u\to u\)
  is a terminal object of the comma category \(S/u\).
\end{proof}

We can characterize lax limit and colimits purely
in terms of the matrix calculus encoded in the category
\(\holaxmat{\extwocat}\):

\begin{lem}
  \label{lem:lax-limits-isos-holaxmat}
  Let \(\X\colon S\to \extwocat\) be a diagram.
  \begin{enumerate}
  \item
    A lax cone \(P:\laxmat{\extwocat}(L,\X)\) is a lax limit cone
    if and only if it is an isomorphism in the category
    \(\holaxmat{\extwocat}\).
  \item
    A lax cocone \(I:\laxmat{\extwocat}(\X,L)\) is a lax colimit cone
    if and only if it is an isomorphism in the category
    \(\holaxmat{\extwocat}\).
  \end{enumerate}
\end{lem}
\begin{proof}
  We prove the statement about lax cones; the other one is dual.

  First assume that \(P\) is a lax limit cone.
  Let \(\Y\colon T\to \extwocat\) be a diagram in \(\extwocat\).
  Using the defining universal property of \Cref{defi:lax-limit}
  on \(L'\coloneqq \Y_t\) for each \(t:T^\op\), we see that the functor
  \begin{equation}
    {P\circ-}
    \colon
    \laxmat{\extwocat}(\Y,L)
    =
    \laxlim_{t:T^\op}\extwocat(\Y_t,L)
    \xrightarrow{\simeq}
    \laxlim_{(s,t):S\times T^\op}\extwocat(\Y_t,\X_s)
    =\laxmat{\extwocat}(\Y,\X)
  \end{equation}
  is an equivalence of \(\infty\)-categories.
  In particular, after passing to equivalence classes of matrices,
  the map
  \begin{equation}
    {P\circ-}\colon \holaxmat{\extwocat}(\Y,L)\xrightarrow{\cong}\holaxmat{\extwocat}(\Y,\X)
  \end{equation}
  is a bijection.
  Since \(\Y\) was arbitrary, it follows that \(P\colon L\to\X\)
  is an isomorphism in the category \(\holaxmat{\extwocat}\).
 
  Conversely, assume that \(P\) has an inverse in \(\holaxmat{\extwocat}\),
  i.e., a lax cocone \(I:\laxmat{\extwocat}(\X,L)\)
  satisfying \(P\circ I \simeq \unitmat^\X\)
  and \(I\circ_S P\simeq \id_L\).
  Then for each \(L':\extwocat\) we have equivalences
  \begin{equation}
    (P\circ-)\circ (I\circ_S-)
    = (P\circ(I\circ_S-))
    \simeq {(P\circ I)\circ_S-}
    \simeq {\unitmat^\X\circ_S-}
    \simeq \id
  \end{equation}
  and
  \begin{equation}
    (I\circ_S-)\circ (P\circ-)
    = {I\circ_S(P\circ -)}
    \simeq {(I\circ_S P)\circ -}
    \simeq {\id_L\circ -}
    = \id
  \end{equation}
  as endofunctors of
  \begin{equation}
    \laxmat{\extwocat}(L',\X)
    \quad\text{and}\quad
    \extwocat(L',L),
  \end{equation}
  respectively
  (using \Cref{lem:lax-matrix-associative}),
  showing that \(P\circ-\) is an equivalence of \(\infty\)-categories,
  as required.
\end{proof}

Recall the standing assumption of this section,
that \(\extwocat\) is enriched in \(\infty\)-categories with colimits.

\begin{cor}
  \label{cor:lax-co-limits-are-bilimits}
  A diagram \(\X\colon S\to\extwocat\)
  admits a lax limit if and only if it admits a lax colimit.
  When they exist, the unit matrix
  \(\unitmat^\X:\laxmat{\extwocat}(\X,\X)\)
  corresponds to an equivalence
  \begin{equation}
    \label{eq:laxcolim=laxlim}
    \unitmatmap^\X\colon \laxcolim_S\X\xrightarrow{\simeq}\laxlim_S\X.
  \end{equation}
\end{cor}

\begin{proof}
  The diagram \(\X\) admits a lax (co)limit
  if and only if it is isomorphic in \(\holaxmat{\extwocat}\)
  to an object \(L\colon \singleton\to\extwocat\).
  In this case \(L\) is both the lax limit and the lax colimit,
  exhibited by mutually inverse lax (co)cones
  \(I\colon\X\to L\) and \(P\colon L\to \X\).
  By definition,
  the map \(\unitmatmap^\X\colon L\to L\)
  corresponding to the matrix \(\unitmat^\X\)
  is determined (up to equivalence) by the property that
  \(P\circ\unitmatmap^\X\circ I\simeq \unitmat^\X\).
  Since the identity \(\id_L\) satisfies this property,
  we conclude that \(\unitmatmap^\X\simeq \id_L\);
  in particular this map is an equivalence.
\end{proof}

\begin{rem}
  The comparison map \eqref{eq:laxcolim=laxlim}
  does not just depend on the objects
  \begin{equation}
    L=\laxlim\X
    \quad\text{and}\quad
    L'=\laxcolim\X
  \end{equation}
  but on the implicit lax (co)limit cones
  \(P\colon \laxmat{\extwocat}(L,\X)\)
  and
  \(I\colon \laxmat{\extwocat}(\X,L')\).
  Specifically, the map \(\unitmatmap^\X\)
  is characterized up to equivalence by the relation
  \begin{equation}
    P\circ\unitmatmap^\X\circ I \simeq \unitmat^\X
  \end{equation}
  or equivalently
  \begin{equation}
    \label{eq:inverse-of-comparison-map}
    \unitmatmap^\X\simeq (I\circ_SP)^{-1}
  \end{equation}
  (since \(P\) and \(I\) are isomorphisms in \(\holaxmat{\extwocat}\)
  and \(\unitmat^X\) is the identity on \(\X:\holaxmat{\extwocat}\)).
\end{rem}

Having described lax (co)limits via lax matrix formulas in
\(\holaxmat{\extwocat}\),
we can immediately deduce that all lax (co)limits are absolute
with respect to the colimit-enrichment.

\begin{thm}
  \label{thm:lax-lims-absolute}
  Let \(\extwocat,\extwocat'\) be \((\infty,2)\)-categories
  enriched in \(\infty\)-categories with colimits.
  Let \(\extwofun\colon\extwocat\to\extwocat'\)
  be a functor which preserves colimits on hom-cateogories.
  Then \(\extwofun\) preserves all lax colimits and lax limits.
\end{thm}

\begin{proof}
  \newcommand\holaxfun{\mathrm{hoLaxMat}_\extwofun}
  Since \(\extwofun\) preserves colimits on hom-categories,
  it induces a well defined functor
  \begin{equation}
    \holaxfun\colon \holaxmat{\extwocat}\to\holaxmat{\extwocat'},
  \end{equation}
  given by applying \(\extwofun\) pointwise to diagram and matrices.
  Since this functor necessarily sends isomorphisms to isomorphisms,
  \Cref{lem:lax-limits-isos-holaxmat} implies
  that \(\extwofun\) sends lax (co)limit cones to lax (co)limit cones.
\end{proof}

Finally, we deduce that lax matrix multiplication
corresponds to composition of maps between lax colimits/limits
in the case where those lax (co)limits exist.

\begin{prop}
  \label{prop:matrix-multiplication}
  Let
  \(\X\colon S\to\extwocat\),
  \(\Y\colon T\to\extwocat\),
  \(\Z\colon U\to\extwocat\)
  be diagrams indexed by \(\infty\)-categories
  and admitting lax limits/colimits.
  Then there is a commutative square of \(\infty\)-categories
  \begin{equation}
    \begin{tikzcd}
      \substack{
        \extwocat(\laxcolim_S\X,\laxlim_T\Y)
        \times
        \extwocat(\laxcolim_T\Y,\laxlim_U\Z)
      }
      \ar[d,"{-\circ(\unitmatmap^\Y)^{-1}\circ-}"]
      \ar[r,"\simeq"]
      &
      \substack{
        \laxmat{\extwocat}(\X,\Y)
        \times
        \laxmat{\extwocat}(\Y,\Z)
      }
      \ar[d,"-\circ_T-"]
      \\
      \extwocat(\laxcolim_S\X,\laxlim_U\Z)
      \ar[r,"\simeq"]
      &
      \laxmat{\extwocat}(\X,\Z)
    \end{tikzcd}
  \end{equation}
  In other words,
  after identifying lax colimits and lax limits via the canonical unit matrix,
  lax matrix multiplication corresponds precisely to function composition.
\end{prop}

\begin{proof}
  Denote by
  \(I_\X\colon \laxmat{\extwocat}(\X,\laxcolim\X)\) and
  \(P_\X\colon \laxmat{\extwocat}(\laxlim\X,\X)\)
  two lax (co)limits cones for the diagram \(\X\)
  (and similarly for \(\Y\) and \(\Z\)).
  The implicit identification
  \begin{equation}
    \extwocat(\laxcolim\X,\laxlim\Y)
    \xrightarrow{\simeq}
    \laxmat{\extwocat}(\X,\Y)
  \end{equation}
  is given explicitly as
  \(P_\Y\circ - \circ I_\X\),
  and similarly for the other horizontal maps.
  Therefore the desired commutative square is just the natural equivalence 
  \begin{align}
    (P_\Z\circ-\circ I_\Y)\circ_T(P_\Y\circ-\circ I_\X)
    &\simeq
    P_\Z\circ -\circ (I_\Y\circ_T P_\Y) \circ - \circ I_\X
      \\
    &\simeq P_\Z\circ (- \circ (\unitmatmap^\Y)^{-1} \circ -)\circ I_\X.
  \end{align} 
  using the equivalence \eqref{eq:inverse-of-comparison-map}
  in the last step.
\end{proof}

\begin{exa}
  \label{exa:2x2-matrices-matrix-comp}
  Continuing \Cref{exa:2x2-matrices-unit-matrix},
  we consider a \(\Delta^1\)-diagram \(\Y=(\Y_0\xrightarrow{\Y_{10}}\Y_1)\)
  in \(\extwocat\)
  and two lax matrices
  \begin{equation}
    F=
    \Donevector{F_0}{F_1}
    : \laxmat{\extwocat}(X,\Y)
    \quad\text{and}\quad
    G=\Donecovector{G_0}{G_1}
    : \laxmat{\extwocat}(\Y,Z)
  \end{equation}
  (just a lax cone and a lax cocone, really).
  The Cartesian twisted arrow category \(\Twcontra(\Delta^1)\)
  is the poset
  \begin{equation}
    \begin{tikzcd}
      \id_0\ar[from=r]&(0\xrightarrow{10} 1)\ar[d]\\
      &\id_1
    \end{tikzcd}
  \end{equation}
  The matrix product \(G\circ_{\Delta^1} F : \laxmat{\extwocat}(X,Z)=\extwocat(X,Z)\)
  is therefore the pushout of the diagram
  \begin{equation}
    \begin{tikzcd}
      G_0F_0\ar[from=r,"g F_0"']&G_1\Y_{10} F_0\ar[d,"G_1 f"]\\
      &G_1 F_1
    \end{tikzcd}
  \end{equation}
  computed in the \(\infty\)-category \(\extwocat(X,Y)\),
  where \(g\colon G_1\Y_{10}\to G_0\) and \(f\colon \Y_{10}F_0\to F_1\)
  are the \(2\)-cells encoded in \(G\) and \(F\), respectively.
  More general \(\Delta^1\times (\Delta^1)^\op\) matrices can then be multiplied
  in the usual row-by-column fashion since each entry \((GF)_{us}\)
  only depends on row \(u\) of \(G\) and column \(s\) of \(F\).
  For example,
  we can compute (with \(\X=\Y\) and \(F=\unitmat^\Y\))
  \begin{gather}
    G\circ_{\Delta^1} \unitmat=
    \Donematrix{G_{00}}{G_{01}}{G_{10}}{G_{11}}
    \Donematrix{\id_{\Y_0}}{\initialobj}{\Y_{10}}{\id_{\Y_1}}
    \\
    \simeq
    \Donematrix
    {\colim\left(G_{00}\leftarrow{G_{01}\Y_{10}\xrightarrow{=} G_{01}\Y_{10}}\right)}
    {\colim\left(\initialobj\xleftarrow{=}\initialobj\to G_{01}\right)}
    {\colim\left(G_{10}\leftarrow{G_{11}\Y_{10}\xrightarrow{=} G_{11}\Y_{10}}\right)}
    {\colim\left(\initialobj\xleftarrow{=}\initialobj\to G_{11}\right)}
    \simeq G
  \end{gather}
\end{exa}

\section{Lax additivity}\label{sec10.4}

Classical semi-additivity of a category \(\A\) manifests itself on two levels:
\begin{enumerate}
\item
  Each hom-set of \(A\) has a commutative monoid structure
  which allows to take sums \(\sum_{s\in S}f_s\) indexed by arbitary finite sets \(S\).
\item
  The category \(A\) allows
  direct sums \(\bigoplus_{s\in S}x_s\) indexed by finite sets \(S\)
  which are both products and coproducts.
\end{enumerate}

We categorify these notions by replacing (discrete) addition
\(\Sigma_{s\in S}\) on the hom-sets
by colimits \(\colim_{s: S}\) on the hom-categories
and (discrete) coproducts/products \(\coprod_{s\in S}\simeq \prod_{s\in S}\)
by lax bilimits \(\laxcolim_{s:S}\simeq\laxlim_{s:S}\),
which are now indexed by arbitary small
\(\infty\)-categories \(S\) rather than finite sets.

\begin{defi}
  Let \(\extwocat\) be an \((\infty,2)\)-category enriched
  in \(\infty\)-categories with groupoidal colimits.
  (This means that each hom-category \(\extwocat(X,Y)\) has all colimits indexed by \(\infty\)-groupoids
  and that composition preserves such colimits in each variable.)
  Let \(\X\colon S\to \extwocat\) a diagram
  indexed by an \(\infty\)-category \(S\).
  A \emph{lax bilimit} of \(\X\) consists of
  a lax colimit cone \(I:\laxmat{\extwocat}(\X,L')\)
  and
  a lax limit cone \(P:\laxmat{\extwocat}(L,\X)\)
  such that the canonical map
  \begin{equation}
    \unitmatmap^\X\colon L'\to L
  \end{equation}
  corresponding to the unit matrix
  \(\unitmat^\X:\laxmat{\extwocat}(\X,\X)\)
  is an equivalence.
  We identify \(L\) and \(L'\) via \(\unitmatmap^\X\)
  and write
  \begin{equation}
    \laxbilim_S\X
    \quad\text{or}\quad\laxbilim_{s:S}\X_s
  \end{equation}
  for both/either of them.

  When \(\X\colon S\to\singleton\xrightarrow{X}\extwocat\)
  is a constant diagram,
  we write \(\laxtensor{S}{X}\) or \(\laxpower{S}{X}\)
  for the constant lax bilimit
  \(\laxbilim_{s:S}X= \laxbilim_S\X\).
  In the special case where \(S=\Delta^1\) so that
  \(\X=(\X_0\xrightarrow{H} \X_1)\) is just an arrow in \(\extwocat\),
  we also write
  \(\X_0\laxplus[H]\X_1\) instead of  \(\laxbilim_{\Delta^1}\X\).
\end{defi}

\begin{rem}
  When convenient we drop the typographical distinction between a matrix
  \(F:\laxmat{\extwocat}(\X,\Y)\) and the associated map
  \(\mattomap{F}:\laxbilim\X\to\laxbilim\Y\).
  More generally, \Cref{prop:matrix-multiplication} justifies
  dropping the typographical distinction between
  matrix multiplication \(G\circ_T F\)
  and composition \(\mattomap{G}\circ \mattomap{F}\) of the associated maps 
  \begin{equation}
    \begin{tikzcd}
      \laxbilim_S \X
      \ar[r,"\mattomap{F}"]
      &
      \laxbilim_T \Y
      \ar[r,"\mattomap{G}"]
      &
      \laxbilim_U \Z
    \end{tikzcd}
  \end{equation}
  between lax bilimits in \(\extwocat\).
\end{rem}

We can now finally define the notion of lax semiadditivity.

\begin{defi}
  An \((\infty,2)\)-category \(\extwoaddcat\)
  is called \emph{(finitely) lax semiadditive} if
  \begin{enumerate}
  \item
    it is enriched in \(\infty\)-categories with (finite) colimits
    (with functors preserving them),
  \item
    each diagram \(S\to \extwoaddcat\)
    indexed by a (finite) small \(\infty\)-category \(S\)
    admits a lax bilimit.
  \end{enumerate}
\end{defi}

\begin{rem}
  We have seen in \Cref{cor:lax-co-limits-are-bilimits} that
  in the presence of sufficently many colimits in the hom-categories,
  \emph{every} lax limit or colimit is automatically a lax bilimit.
  Thus the second condition could be weakened to just require the existence
  of lax limits \emph{or} lax colimits.
\end{rem}

The final step is to categorify the passage from
semiadditve to additive categories
which amounts to requiring the hom-monoids to be abelian groups.
Following our philosophy of \S\ref{subsec:rules-of-categorification},
abelian groups should be replaced by stable \(\infty\)-categories
leading to the following easy definition:

\begin{defi}
  \label{defi:lax-additive}
  A (finitely) lax semiadditive \((\infty,2)\)-category \(\extwoaddcat\)
  is called \emph{(finitely) lax additive}
  if every hom-\(\infty\)-category \(\extwoaddcat(X,Y)\) is stable.
\end{defi}

\begin{rem}
  Denote by \(\extwocat^\op\) the \((\infty,2)\)-category
  obtained from \(\extwocat\)
  by reversing the directions of the \(1\)-morphisms,
  i.e., \(\extwocat^\op(X,Y)=\extwocat(Y,X)\).
  If \(\extwocat\) is enriched in (stable)
  \(\infty\)-categories with colimits,
  then so is \(\extwocat^\op\).
  Moreover lax limits/colimits/bilimits in \(\extwocat^\op\)
  correspond to lax colimits/limits/bilimits in \(\extwocat\).
  Thus an \((\infty,2)\)-category \(\extwoaddcat\)
  is (finitely) lax (semi)additive if and only if \(\extwoaddcat^\op\)
  is (finitely) lax (semi)additive.
\end{rem}

\begin{exa}
  Lax limits in the \((\infty,2)\)-category \(\PrL\) of presentable
  \(\infty\)-categories exist
  and are computed as underlying \(\infty\)-categories,
  i.e., as sections of the Grothendieck construction as in \Cref{ex:lax-lim-cat};
  indeed, for any small diagram \(\X\colon S\to \PrL\)
  the \(\infty\)-category \(L\coloneqq\Fun_S(S,\covGroth \X)\) is again presentable
  (see \cite[Proposition~5.5.3.17 and Proposition~5.5.3.3]{lurie:htt})
  and a functor into it preserves colimits
  if and only if it does so after postcomposing
  with each pointwise evaluation map \(L\to \X_s\).
  Since \(\PrL\) is enriched in \(\infty\)-categories with colimits
  it follows that it is a lax semiadditive \((\infty,2)\)-category.
  The full \((\infty,2)\)-category \(\StL\)
  of presentable \emph{stable} \(\infty\)-categories
  is closed under lax colimits and enriched in stable \(\infty\)-categories,
  thus it is lax additive.
  The \((\infty,2)\)-category \(\StR\)
  of presentable stable \(\infty\)-categories and right adjoint functors
  is only finitely lax additive,
  since composition with a right adjoint functor is exact
  but does not preserve arbitrary colimits.
  The \((\infty,2)\)-category \(\st\) of stable \(\infty\)-categories
  and exact functors is finitely lax additive.

  Note that one can replace all presentability assumptions
  by just requiring the relevant \(\infty\)-categories
  to have colimits (or limits, in the case of \(\StR\))
  and for the functors between them to preserve them.
\end{exa}

\begin{rem}
  In the \((\infty,2)\)-category \(\StL\),
  the lax bilimit of a diagram \(F\colon \A\to \B\)
  comes with a semiorthogonal decompositions with components \(\A\) and \(\B\)
  and gluing functor \(F\) in the sense of \cite{DKSS:spherical},
  see also \Cref{sec:recollements}.
  The idea to use matrices to describe coordinate change
  for such semiorthogonal decompositions already appears in \cite{DKS:nspherical}.
  Enlarging $\StL$ by a suitably defined $(\infty,2)$-category of exact
  profunctors, we may even describe general semiorthogonal decompositions via
  lax bilimits -- in the context of dg categories, this corresponds to the gluing
  constructions for bimodules, as established in \cite{KL15}.
\end{rem}

\begin{con}
  \label{con:interpret-space-matrices}
  Let \(T,S\) be small \(\infty\)-categories and
  let
  \begin{equation}
    F\colon T\times S^\op\to \Spaces\simeq \PrL(\Spaces,\Spaces)
  \end{equation}
  be a matrix of spaces.
  Let \(\extwoaddcat\) be a lax semiadditive \((\infty,2)\)-category.
  For each \(X:\extwoaddcat\), denote by \(F^X\coloneqq F\otimes \id_X\)
  the matrix
  \begin{equation}
    F^X\colon T\times S^\op\xrightarrow{F}\Spaces
    \xrightarrow{-\otimes \id_X}\extwoaddcat(X,X),
  \end{equation}
  where the second functor arises from the tensoring by \(\Spaces\)
  on the \(\infty\)-category \(\extwoaddcat(X,X)\) with colimits.
  In this way, the space-valued matrix \(F\)
  gives rise to a map \(F^X\colon X^S\to X^T\),
  which we call the action of \(F\) on \(X\).
  Similarly, when the hom-categories are pointed/stable,
  we can interpret in \(\extwoaddcat\) every matrix of pointed spaces/spectra,
  by using the corresponding tensoring.
\end{con}

\begin{lem}
  Action of matrices is compatible with matrix multiplication,
  i.e., we have equivalences
  \(F^X\circ G^X \simeq (F\circ G)^X\)
  whenever \(F,G\) are composable matrices of spaces/pointed spaces/spectra
  and \(X\) is an object in a correspondingly enriched \((\infty,2)\)-category.
\end{lem}
\begin{proof}
  We have, naturally in \(u:U\), \(s:S^\op\):
  \begin{align}
    (F\circ G)^X_{us}
    &\simeq
      \left(\colim_{f\colon t\to t'}F_{ut'}\otimes G_{ts}\right)\otimes \id_X
      \\
    &\simeq
      \left(\colim_{f\colon t\to t'}(F_{ut'}\otimes \id_X)\circ (G_{ts}\otimes \id_X)\right)
    \\
     &\simeq (F^X\circ G^X)_{us}
  \end{align}
  using that the tensoring preserves colimits
  and that \(\id_X\circ \id_X\simeq \id_X\).
\end{proof}

\begin{exa}
  The universal identity \(\Delta^1\times (\Delta^1)^\op\)-matrix is
  \begin{equation}
    \unitmat=
    \Donematrix{\singleton}{\emptyset}{\singleton}{\singleton}
    \colon \Delta^1\times(\Delta^1)^\op\to \Spaces.
  \end{equation}
  Indeed, every \((\infty,2)\)-category \(\extwocat\) whose
  hom-categories have initial objects,
  the unit matrix \(\unitmat^\X\) for the constant diagram
  \(\X\colon \Delta^1\to\singleton\xrightarrow{X}\extwocat\)
  is given by \(\unitmat^\X=\unitmat^X\coloneqq\unitmat\otimes\id_X\).

  More generally, for any \(\infty\)-category \(S\)
  the universal idenity \(S\times S^\op\)-matrix
  is the just the transpose of the hom-functor
  \begin{equation}
    \matrans{\unitmat}=S(-,-)\colon S^\op\times S\to\Spaces.
  \end{equation}
\end{exa}

\begin{exa}
  \label{exa:matrices-fib-cof}
  Consider the matrices
  \begin{equation}
    \Cof\coloneqq\Donematrix{\sphere^0}{\sphere^0}{0}{\sphere^0}
    \colon \Delta^1\times (\Delta^1)^\op\to \SpacesP
  \end{equation}
  and
  \begin{equation}
    \Fib\coloneqq\Donematrix[bicart]{0}{{\sphere[-1]}}{\sphere}{0}
    \colon \Delta^1\times (\Delta^1)^\op\to \Sp,
  \end{equation}
  defined in pointed spaces and spectra, respectively.
  If \(\A\) is a pointed \(\infty\)-category with colimits,
  then the matrix \(\Cof\) acts as the cofiber map
  \begin{equation}
    \Cof\nolimits_\A=\Donematrix{\id_A}{\id_A}{0}{\id_A}
    \colon \A^{\Delta^1}\to\A^{\Delta^1},
  \end{equation}
  Indeed, for every
  \begin{equation}
    \Donevector{a}{b}\colon\B\to\A^{\Delta^1}
  \end{equation}
  we can compute the matrix product
  \begin{equation}
    \Donematrix{\id_A}{\id_A}{0}{\id_A} \circ \Donevector{a}{b}
    \simeq
    \Donevector
    {\colim(a\leftarrow a\to b)}
    {\colim(0\leftarrow a\to b)}
    \simeq
    \Donevector
    {b}{\cof(a\to b)} = \Cof(a\to b).
  \end{equation}
  If \(\A\) is also stable,
  then a similar calculation shows that the matrix \(\Fib\)
  is inverse to \(\Cof\) and acts as the fiber map.
\end{exa}

For every matrix \(F\colon T\times S^\op\to \exaddcat(X,X)\)
corresponding to an arrow \(\laxpower{S}{X}\to \laxpower{T}{X}\)
in \(\extwoaddcat\),
we denote by \(\matrans{F}\) the ``transposed'' matrix
\begin{equation}
  \matrans{F}\colon S^\op\times (T^\op)^\op \simeq T\times S^\op \xrightarrow{F}
  \exaddcat(X,X).
\end{equation}
describing the dual map \(\laxpower{T^\op}{X}\to\laxpower{S^\op}{X}\).

\begin{lem}
  \label{lem:matrix-action-by-homs}
  In the setting of \Cref{con:interpret-space-matrices}
  there are commutative diagrams
  \begin{equation}
    \begin{tikzcd}
      \extwoaddcat(X,\laxpower{S}{Y})
      \ar[r,"\simeq"]
      \ar[d,"{\extwoaddcat(X,F^Y})"']
      &
      \extwoaddcat(X,Y)^S
      \ar[d,"{F^{\extwoaddcat(X,Y)}}"]
      \\
      \extwoaddcat(X,\laxpower{T}{Y})
      \ar[r,"\simeq"]
      &
      \extwoaddcat(X,Y)^T
    \end{tikzcd}
    \quad\text{and}\quad
    \begin{tikzcd}
      \extwoaddcat(\laxtensor{S}X,Y)
      \ar[r,"\simeq"]
      \ar[from=d,"{\extwoaddcat(F^X,Y)}"]
      &
      \extwoaddcat(X,Y)^{S^\op}
      \ar[from=d,"{\matrans{F}^{\extwoaddcat(X,Y)}}"']
      \\
      \extwoaddcat(\laxtensor{T}X,Y)
      \ar[r,"\simeq"]
      &
      \extwoaddcat(X,Y)^{T^\op}
    \end{tikzcd}
  \end{equation}
\end{lem}

\begin{proof}
  We do the second computation; the first one is similar.
  Functorially in
  \(M : \extwoaddcat(X,Y)^{T^\op}\simeq \extwoaddcat(\laxtensor{T}X,Y)\)
  and \(s: S^\op\)
  we compute
  \begin{align}
    \extwoaddcat(F^X,Y)(M)_s= (M\circ F^X)_s
    &= \colim_{(f\colon t\to t'):\Twcontra(T)} M_{t'}\circ F^X_{ts}
    \\
    &= \colim_{f:\Twcontra(T)} M_{t'}\circ (F_{ts}\otimes \id_X)
    \\
    &\simeq  \colim_{(f\colon t\to t'):\Twcontra(T)} F_{ts}\otimes M_{t'}
    \\
    &\simeq \colim_{(g\colon t'\to t) : \Twcontra(T^\op)} F_{ts}\otimes M_{t'}
    \\
    &\simeq \colim_{(g\colon t'\to t) : \Twcontra(T^\op)}
      (F_{ts}\otimes \id_{\extwoaddcat(X,Y)})(M_{t'})
    \\
    &=\colim_{(g\colon t'\to t) : \Twcontra(T^\op)}
      ((\matrans{F})_{st}\otimes \id_{\extwoaddcat(X,Y)})(M_{t'})
    \\
    &= \matrans{F}^{\extwoaddcat(X,Y)}(M)_s.
  \end{align}
  Apart from expanding the various definitions, we have used
  \begin{itemize}
  \item
    that \(M_{t'}\circ -\colon \extwoaddcat(X,X)\to\extwoaddcat(X,Y)\)
    preserves colimits and hence the tensoring \(F_{ts}\otimes-\),
  \item
    the canonical identification \(\Twcontra(T)\simeq\Twcontra(T^\op)\)
    which reverses source and target,
  \item
    that colimits of functors \(\extwoaddcat(X,Y)\to\extwoaddcat(X,Y)\)
    are computed pointwise,
    hence also the tensoring \(F_{ts}\otimes -\).
    \qedhere
  \end{itemize}
\end{proof}

\begin{lem}
  Let \(\extwoaddcat\) be a lax semiadditive \((\infty,2)\)-category.
  Then \(\extwoaddcat\) is lax additive if and only if
  each hom-category \(\extwoaddcat(X,Y)\) is pointed
  and the matrix \(\Cof\) from \Cref{exa:matrices-fib-cof}
  acts invertibly on each \(X:\extwoaddcat\).
  If this is the case, then the inverse is given by the action of \(\Fib\).
\end{lem}

\begin{proof}
  Assuming that all hom-categories \(\extwoaddcat(X,Y)\)
  of the lax semiadditive \((\infty,2)\)-category \(\extwoaddcat\)
  are pointed,
  they are stable if and only if the cofiber functor
  \begin{equation}
    \Cof\nolimits^{\extwoaddcat(X,Y)}
    \colon \extwoaddcat(X,Y)^{\Delta^1}
    \to \extwoaddcat(X,Y)^{\Delta^1}
  \end{equation}
  is invertible.
  Using \Cref{lem:matrix-action-by-homs},
  we can identify this functor with
  \(\extwoaddcat(X,\Cof^Y)\).
  Thus
  \(\extwoaddcat\) is lax additive
  if and only if
  all \(\extwoaddcat(X,Y)\) are stable,
  if and only if
  all \(\extwoaddcat(X,\Cof^Y)\) are invertible,
  if and only if
  all \(\Cof^Y\) are invertible,
  as claimed.
\end{proof}

\subsection{Oplax additivity}

So far we have focused our discussion exclusively
on \emph{lax} limits and colimits, as opposed to \emph{oplax} ones.
We could have of course passed to the 2-morphism dual everywhere
(obtained from an \((\infty,2)\)-category \(\extwocat\)
by replacing each hom-category \(\extwocat(X,Y)\)
by its opposite)
and told an analogous story using \emph{oplax} colimits/limits/bilimits.
This would lead to what we might call
oplax (semi)additive \((\infty,2)\)-categories \(\extwoaddcat\),
which are enriched in (stable) \(\infty\)-categories with limits
and allow the formation of oplax bilimits
\begin{equation}
  \oplaxbilim_S\X \coloneqq \oplaxcolim_S\X\simeq \oplaxlim_S\X
\end{equation}
of any diagram \(\X\colon S\to\extwoaddcat\)
indexed by a small \(\infty\)-category.

For the convenience of the reader,
we summarize the main formulas of this dual theory;
all the constructions and proofs are dual to the ones we saw earlier.

\begin{enumerate}
\item
  For two diagrams \(\X\colon S\to \extwoaddcat\)
  and \(\Y\colon T\to\extwocat\),
  we define
  \begin{equation}
    \oplaxmat{\extwocat}(\X,\Y)
    \coloneqq\oplaxlim_{(t,s):T\times S^\op}\extwocat(\X_s,\Y_t)
  \end{equation}
  as the category of oplax matrices from \(\X\) to \(\Y\).
  Explicitly, such matrices are sections
  of the contravariant Grothendieck construction
  \begin{equation}
    \contraGroth[(t,s):T\times S^\op]{\extwoaddcat(\X_s,\Y_t)}
    \longrightarrow T^\op\times S.
  \end{equation}
\item
  The oplax matrix multiplication
  \begin{equation}
    \oplaxmat{\extwocat}(\X,\Y)
    \times
    \oplaxmat{\extwocat}(\Y,\Z)
    \to
    \oplaxmat{\extwocat}(\X,\Z)
  \end{equation}
  is given by the formula
  \begin{equation}
    (\Phi\circ\Psi)_{us}
    \coloneqq\lim_{(f\colon t\to t'):\Twcov(T)}\Phi_{ut'}\Y_{f}\Psi_{ts}.
  \end{equation}
\item
  The oplax unit matrix for a diagram \(\X\colon S\to\extwocat\) is
  \begin{equation}
    \unitmatop^\X = \left((t,s)\mapsto\lim_{f: S(s,t)}\X_f\right)
    : \oplaxmat{\extwocat}(\X,\X)
  \end{equation}
\end{enumerate}

\begin{exa}
  The \((\infty,2)\)-category \(\StL\) is
  enriched in stable \(\infty\)-categories and admits oplax limits;
  thus it is finitely oplax additive.
  It is \emph{not} oplax additive
  because composition of functors does not preserve arbitrary limits.
\end{exa}

\begin{exa}
  \label{exa:oplax-Delta-1-mult}
  As in \Cref{exa:2x2-matrices-matrix-comp},
  consider a \(\Delta^1\)-diagram \(\Y=(\Y_0\xrightarrow{\Y_{10}}\Y_1)\) in \(\extwocat\).
  The oplax matrix product over \(\Y\) of
  an oplax cocone (= \(\Delta^1\)-row vector) \(G\)
  with
  an oplax cone (= \((\Delta^1)^\op\)-column vector) \(F\)
  is given by the dual formula
  \begin{equation}
    \Donecovectorop{G_0}{G_1}
    \circ_{\Delta^1}
    \Donevectorop{F_0}{F_1}
    =
    \lim\left(
      \begin{tikzcd}
        G_0F_0\ar[r]& G_1 \Y_{10} F_0
        \\
        &
        G_1F_1\ar[u]
      \end{tikzcd}
    \right).
  \end{equation}
  General \({(\Delta^1)^\op\times \Delta}\)-matrices can then be multiplied
  in the usual row-by-column way.
  The oplax unit matrix on \(\Y\) is
  \begin{equation}
    \unitmatop^\Y
    =
    \Donematrixop
    {\id_{\Y_0}}{\singleton}{\Y_{10}}{\id_{\Y_1}},
  \end{equation}
  where \(\singleton\) is the terminal object of \(\extwocat(\Y_1,\Y_0)\).
\end{exa}

\section{Coordinate change for \texorpdfstring{\(\Delta^1\)}{Delta1}-matrices}

We have seen that any lax semiadditive \((\infty,2)\)-category
admits a nicely behaved calculus of lax matrices.
However, if we apply \(K_0\) componentwise
to the lax matrix multiplication for lax \(\Delta^1\)-bilimits
(see \Cref{exa:2x2-matrices-unit-matrix}
and \Cref{exa:2x2-matrices-matrix-comp})
we obtain the very unusual formula
\begin{equation}
  \ppsmallmatrix{a_0& a_1}\circ\ppsmallmatrix{b_0\\ b_1} = a_0b_0 - a_1b_0 + a_1b_1
\end{equation}
or, more generally \(A\circ B = A I^{-1} B \),
where \(I=K_0(\unitmat^{\Delta^1})=\ppsmallmatrix{1&0\\1&1}\)
is the new unit matrix.

The goal of this section is to introduce a convenient
``coordinate change'' in the lax additive
(as opposed to merely lax \emph{semi}additive) setting,
which up to a sign recovers the usual matrix multiplication on \(K_0\).

The key ingredient is the cofiber-fiber-equivalence
\begin{align}
  \Coff\colon \Fun(\Delta^1,\A) &\xleftrightarrow[\phantom{\longrightarrow}]{\simeq}
  \Fun(\Delta^1,\A)\noloc \Fibb
  \\
  (\fib(u') =b\xrightarrow{u}a)
  &\leftrightarrow (a\xrightarrow{u'} b' = \cof(u))
\end{align}
for every stable \(\infty\)-category \(\A\).

More precisely, we make use of the following dependent version
of the cofiber-fiber-equivalence
which identifies the oplax limit over an arrow with the lax limit.

\begin{lem}[Lemma~1.3 in \cite{DJW19-BGP}]
  \label{lem:dependent-fiber-cofiber}
  Let \(f\colon\A\to\B\) be a diagram of stable \(\infty\)-categories.
  Then there is a natural equivalence
  \begin{equation}
    \label{eq:dependent-cofiber-fiber-equiv}
    \Coff\colon
    \oplaxlim_{\Delta^1}(\B\xleftarrow{f}\A)
    \xleftrightarrow{\simeq}
    \laxlim_{\Delta^1}(\A\xrightarrow{f}\B)
    \noloc\Fibb
  \end{equation}
  described by the formula
  \begin{equation}
    (b=\fib(u'), a , b\xrightarrow{u} fa)
    \leftrightarrow
    (a,b'=\cof(u), fa\xrightarrow{v} b').
  \end{equation}
\end{lem}

While not strictly necessary,
it is convenient to implement this equivalence
by explicit matrices using a combination of the lax and oplax matrix calculus.

For the remainder of the section,
let \(\extwoaddcat\) be a lax additive \((\infty,2)\)-category.
Then \(\extwoaddcat\) is in particular enriched
in \(\infty\)-categories with finite limits,
so that we have available the finite oplax matrix calculus
(dual to the one in \S\ref{subsec:lax-matrices})
as long as we restrict to diagrams indexed by
\emph{finite} \(\infty\)-categories \(S\).

\begin{con}
  \label{con:oplax-fiber-cones}
  Let \(\X\colon \Delta^1\to \extwoaddcat\) be a diagram,
  \(\X=(\X_0\xrightarrow{F}\X_1)\)
  and let \(P=(p_s)\) and \(I=(i_s)\) (indexed by \(s:\Delta^1\))
  be the lax limit/colimit cone exhibiting
  the lax bilimit \(\X_0\laxplus[F]\X_1\).
  We construct the oplax cone and cocone
  \begin{equation}
    \Fibb\coloneqq
    \oplaxcone {\X_0}{\X_0\laxplus[F]\X_1}{\X_1}
    {p_0} {\fibb} {F}
    \quad\text{and}\quad
    \Fibdual\coloneqq
    \oplaxcocone {\X_0}{\X_0\laxplus[F]\X_1}{\X_1}
    {\fibdual} {i_1} {F}
  \end{equation}
  as follows:
  \begin{itemize}
  \item
    Recall that
    \begin{equation}
      p_0=\Donecovector{\id_{\X_0}}{0}
      \quad\text{and}\quad
      i_1=\Donevector{0}{\id_{\X_1}}
    \end{equation}
    are just the top row and right column of the unit matrix \(\unitmat^\X\),
    viewed as a map out of or into the lax limit of \(\X\), respectively.
  \item
    Additionally we define the row and column vectors
    \begin{equation}
      \fibb\coloneqq\Donecovector{0}{[-1]_{\X_1}}
      \quad\text{and}\quad
      \fibdual\coloneqq\Donevector{[-1]_{\X_0}}{0}
    \end{equation}
    obtained by passing to the vertical and horizontal fibers
    of the unit matrix \(\unitmat^\X\).
    \item
      By construction, \((p_0,\fibb)\) and \((\fibdual,i_1)\) fit into 
      a \((\Delta^1)^\op\times(\Delta^1)^\op\)-matrix
      and a \(\Delta^1\times\Delta^1\)-matrix,
      \begin{equation}
        \Fibb\coloneqq
        \Donevectorop{p_0}{\fibb}
        =
        \Donematrixmixedop[bicart]{\id_{\X_0}}{0}{0}{[-1]_{\X_1}}
        \quad\text{and}\quad
        \Fibdual\coloneqq
        \Donecovectorop{\fibdual}{i_1}
        =
        \Donematrixmixed[bicart]
        {[-1]_{\X_0}}
        {0}{0}
        {\id_{\X_1}},
      \end{equation}
      respectively.
      As indicated, we can view these matrices as a column vector of row vectors,
      or a row vector of column vectors, respectively,
      thus yielding the desired oplax cone and cocone.
    \end{itemize}
  \end{con}

The following lemma explains the name of the cones
\(\Fibb\) and \(\Fibdual\)
in terms of the maps they represent/corepresent.

\begin{lem}
  Let \(Y:\extwoaddcat\).
  \begin{enumerate}
  \item
    The induced map
    \begin{equation}
      \laxlim_{s:\Delta^1}\extwoaddcat\left(Y, \X_s\right)
      \xleftarrow{\simeq}
      \extwoaddcat\left(Y,\X_0\laxplus[F]\X_1\right)
      \xrightarrow{\Fibb\circ -}
      \oplaxlim_{s:\Delta^1}\extwoaddcat(Y, \X_s)
    \end{equation}
    is precisely the dependent fiber functor
    of \Cref{lem:dependent-fiber-cofiber}
    for the \(\Delta^1\)-diagram
    \(\extwoaddcat(Y,\X_0)\to\extwoaddcat(Y,\X_1)\).
  \item
    The induced map
    \begin{equation}
      \laxlim_{s:(\Delta^1)^\op}\extwoaddcat\left(\X_s,Y\right)
      \xleftarrow{\simeq}
      \extwoaddcat\left(\X_0\laxplus[F]\X_1,Y\right)
      \xrightarrow{-\circ \Fibdual}
      \oplaxlim_{s:(\Delta^1)^\op}\extwoaddcat(\X_s,Y)
    \end{equation}
    is precisely the dependent fiber functor
    of \Cref{lem:dependent-fiber-cofiber}
    for the \((\Delta^1)^\op=\Delta^1\)-diagram
    \(\extwoaddcat(\X_0,Y)\leftarrow\extwoaddcat(\X_1,Y)\).
  \end{enumerate}
\end{lem}

\begin{proof}
  A quick matrix computation for each
  \(x =\Donevector{x_0}{x_1} : \laxlim_{s}\extwoaddcat(Y,\X_s)\) shows
  \begin{equation}
    \Fibb \circ x
    = \Donevectorop
    {p_0\circ_{\Delta^1} x}
    {\fibb\circ_{\Delta^1} x}
    = \Donevectorop{x_0}{\fib (Fx_0\to x_1)}
    : \oplaxlim_{s}\extwoaddcat(\X_s,Y),
  \end{equation}
  as required.
  Similarly, for each
  \(\comap{x} = \Donecovector{\comap{x}_0}{\comap{x}_1}
  : \laxlim_{s:(\Delta^1)^\op}\extwoaddcat(\X_s,Y)\)
  we have
  \begin{equation}
    \comap{x}\circ \Fibdual
    =\Donecovectorop{\fib(\comap{x}_1F\rightarrow \comap{x}_0)}{\comap{x}_1}
    : \oplaxlim_{s}\extwoaddcat(\X_s, Y),
  \end{equation}
  as desired.
\end{proof}

As an immediate application of \Cref{lem:dependent-fiber-cofiber}
we therefore get that the oplax cone/cocones \(\Fibb\) and \(\Fibdual\)
exhibit the lax bilimit \(\X_0\laxplus[F]\X_1\)
also as an oplax limit and colimit.
The following lemma makes a more precise statement,
showing that \(\Fibb\) and \(\Fibdual\) are inverse up to a shift.

\begin{lem}
  \label{lem:lax=oplax-over-Delta1}
  The oplax cone \(\Fibb\) and cocone \(\Fibdual\) from 
  \Cref{con:oplax-fiber-cones}
  are, up to negative shift \([-1]\),
  mutually inverse with respect to the oplax matrix multiplication.
  In particular, \(\Fibb\) and \(\Fibdual[1]\) (or \(\Fibb[1]\) and \(\Fibdual\))
  exhibit the lax bilimit \(\X_0\laxplus[F]\X_1\)
  also as an oplax bilimit of the diagram \(\X\colon \Delta^1\to\extwoaddcat\).
\end{lem}
\begin{proof}
  An explicit computation with the oplax matrix multiplication over \(\Delta^1\)
  (recall \Cref{exa:oplax-Delta-1-mult})
  shows:
  \begin{equation}
    \Fibb\circ \Fibdual=
    \Donematrixop{p_0\fibdual}{p_0i_1}{\fibb\fibdual}{\fibb i_1}
    \simeq
    \Donematrixop{\id_{\X_0}[-1]}{0}{F[-1]}{\id_{\X_1}[-1]}
    =\unitmatop^\X[-1]
  \end{equation}
  and
  \begin{align}
    \Fibdual\circ_{\Delta^1} \Fibb
    &\simeq\lim(\fibdual p_0\to i_1Fp_0\leftarrow i_1\fibb)
    \\
    &\simeq
      \lim
      \left(
      \Donematrix{[-1]}{0}{0}{0}
      \to
      \Donematrix{0}{0}{F}{0}
      \leftarrow
      \Donematrix{0}{0}{0}{[-1]}
      \right)
    \\
    &\simeq
      \Donematrix{[-1]}{0}{F[-1]}{[-1]} \simeq\unitmat^\X[-1],
  \end{align}
  where in the second computation
  we omit the straightforward verification
  that the unnamed arrows appearing in the last matrix
  are indeed those of \(\unitmat^\X[-1]\).
\end{proof}

\begin{rem}
  While there is a distinguished choice for the
  cofiber-fiber equivalence \eqref{eq:dependent-cofiber-fiber-equiv},
  the \Cref{lem:lax=oplax-over-Delta1} provides two
  (but equally distinguished) ways
  to identify the lax bilimit \(\X_0\laxplus[F]\X_1\)
  and the oplax bilimit \(\X_0\oplaxplus[F]\X_1\),
  depending on whether we look at the represented map
  (using \(\Fibb\) and treating them as (op)lax limits)
  or the corepresented map
  (using \(\Fibdual\) and treating them as (op)lax colimits).
  These two ways are not equivalent:
  they differ precisely by a suspension.
\end{rem}

\begin{rem}
  We now have several different ways to represent maps
  \(\X_0\laxplus[F]\X_1 \xrightarrow{\alpha} \Y_0\laxplus[G]\Y_1\),
  with the passage between them implemented
  by applying the cofiber-fiber equivalence to rows and/or columns
  of a matrix.

  \begin{equation}
    \label{eq:fibers-of-2x2-matrices}
    \begin{CD}
      {\Donematrixmixedop{\alpha_{00}}{\alpha_{01}}{\alpha^v_{10}}{\alpha^v_{11}}}
      @>{-\circ\Fibdual}>>
      {\Donematrixop{\alpha^h_{00}}{\alpha_{01}}{\alpha^{hv}_{10}}{\alpha^v_{11}}}
      \\
      @A{\Fibb\circ -}AA @A{\Fibb\circ -}AA
      \\
      {\Donematrix{\alpha_{00}}{\alpha_{01}}{\alpha_{10}}{\alpha_{11}}}
      @>{-\circ\Fibdual}>>
      {\Donematrixmixed{\alpha^h_{00}}{\alpha_{01}}{\alpha^h_{01}}{\alpha_{11}}}
    \end{CD}
  \end{equation}
  Here, the notation is chosen as follows:
  \begin{itemize}
  \item
    Subscripts indicate the source and target of each entry, reading right to left.
    For example, \(\alpha_{10},\alpha^h_{10},\alpha^v_{10},\alpha^{hv}_{10}\)
    all live in \(\extwoaddcat(\X_0,\Y_1)\), etc.
  \item
    Superscripts record in which direction (horizontal and/or vertical) one has to take fibers
    to obtain the new object from the original lax-lax matrix (lower left).
    For example, \(\alpha^h_{10}\) is the fiber of the horizontal map
    \(\alpha_{11}F\to \alpha_{10}\)
    while \(\alpha^{hv}_{10}\)
    can be computed
    either as the fiber of the vertical map
    \(G\alpha^h_{00}\to \alpha^h_{10}\)
    or equivalently as the fiber of the horizontal map
    \(\alpha^v_{11}F\to\alpha^v_{10}\).
  \end{itemize}
\end{rem}

\begin{rem}
  \label{lem:mixed-martix-multiplication}
  Consider two composable maps
  \begin{equation}
    \X \xrightarrow{\beta} \Y_0\laxplus[G]\Y_1\xrightarrow{\alpha} \Z,
  \end{equation}
  Each of the two maps \(\beta\) and \(\alpha\)
  can be represented by a matrix in two ways,
  depending on whether we treat
  the middle term as a lax or oplax bilimit.
  The following table shows the four corresponding possible
  row-column-multiplications
  with the standard lax multiplication in the lower left.
  General \(2\times 2\) matrices describing maps between (op)lax
  bilimits over \(\Delta^1\) can then be multiplied in the usual
  row-by-column fashion.

  \begin{equation}
    \label{eq:lax-oplax-multiplications}
    \begin{array}{c|cc}
      \circ
      &
        \Donecovector{\alpha_0}{\alpha_1}
      &
        \Donecovectorop{\alpha^h_0}{\alpha_1}
      \\
      \hline
      \Donevectorop{\beta_0}{\beta^v_1}
      &
        \cof\left(\alpha_1\beta^v_1\to\alpha_1G\beta_0\to\alpha_0\beta_0\right)
      &
        \lim\left(
        \begin{tikzcd}[row sep = tiny, column sep = tiny]
          \alpha^h_0\beta_0
          \ar[r]
          &
          \alpha_1G\beta_0
          \\
          &
          \alpha_1\beta^v_1
          \ar[u]
        \end{tikzcd}
        \right) [1]
      \\
      \Donevector{\beta_0}{\beta_1}
      &
        \colim\left(
        \begin{tikzcd}[row sep = tiny, column sep = tiny]
          \alpha_0\beta_0\ar[from=r]& \alpha_1 G\beta_0\ar[d]\\
          &\alpha_1\beta_1
        \end{tikzcd}
        \right)
      &
        \cof\left(\alpha^h_0\beta_0\to\alpha_1G\beta_0\to \alpha_1\beta_1\right)
    \end{array}
  \end{equation}
  Observe how the entry in top right differs from
  the standard oplax multiplication
  (see \Cref{exa:oplax-Delta-1-mult}) by a shift \([1]\).
  The reason for this is that we used
  the canonical cofiber-fiber-equivalence
  \eqref{eq:dependent-cofiber-fiber-equiv}
  both horizontally and vertically,
  which amounts to using the identification
  \(\Fibb\colon \laxbilim \Y\to \oplaxbilim \Y\)
  when discussing maps \emph{into} the (op)lax bilimit
  but the identification
  \(\Fibdual\colon \oplaxbilim\Y\to\laxbilim\Y\)
  when discussing maps \emph{from} the (op)lax bilimit;
  we have seen in \Cref{lem:lax=oplax-over-Delta1},
  that these two identifications are only inverse up to shift.

  The following table depicts the unit matrix with respect
  to each of the four multiplications;
  they are just obtained from the standard lax unit matrix (lower left)
  by passing to horizontal and/or vertical fibers.

  \begin{equation}
    \begin{array}{c|cc}
      &\text{lax}&\text{oplax}
      \\
      \hline
      \text{oplax}
      &
        \Donematrixmixedop[bicart]{\id}{0}{0}{[-1]}
                 &
                   \Donematrixop{[-1]}{0}{[-1]}{[-1]}
      \\
      \text{lax}
      &
        \Donematrix{\id}{0}{\id}{\id}
                 &
                   \Donematrixmixed[bicart]{[-1]}{0}{0}{\id}
    \end{array}
  \end{equation}
\end{rem}

From the matrix multiplication formulas
of \eqref{eq:lax-oplax-multiplications}
we can immediately see the advantage of this change of coordinates.
By working with lax-oplax or oplax-lax matrices,
we obtain, on \(K_0\) the formulas
\begin{equation}
  \ppsmallmatrix{a_0&a_1}\ppsmallmatrix{b_0\\b_1} =
  \pm(a_0b_0 - a_1b_1).
\end{equation}
The fact that matrix multiplication now involves an alternating sum
rather than an ordinary sum is a feature, rather than a bug.
In the next section we will see, for example,
how we can express the differential of the mapping cone
of a chain map \(f\colon (A_\bullet, \alpha)\to(B_\bullet,\beta)\)
by directly categorifying the canonical matrix
\(\delta=\ppsmallmatrix{\alpha&0\\f &\beta}\)
without having to introduce any signs;
the signs are already part of the matrix multiplication.

Another convenient feature is that the identification
between lax-oplax and oplax-lax matrices
is compatible with the passage to adjoints in the following sense:

\begin{con}
  \label{con:matrix-adjoints}
  Assume that \(G\colon \Y_0\to\Y_1\) has a right adjoint \(G\dashv \ra{G}\).
  Then \Cref{cor:lax-oplax-adjoint}, applied to the adjunctions
  \((G\circ)\dashv(\ra{G}\circ)\) and \((\circ\ra{G})\dashv(\circ G)\)
  yields equivalences
  \begin{equation}
    \extwoaddcat(-,\Y_0\laxplus[G]\Y_1)
    \simeq
    \extwoaddcat(-,\Y_1\oplaxplus[\ra{G}]\Y_0)
    \quad\text{and}\quad
    \extwoaddcat(\Y_1\laxplus[\ra{G}]\Y_0,-)
    \simeq
    \extwoaddcat(\Y_0\oplaxplus[G]\Y_1,-)
  \end{equation}
  given explicitly by passing to vertical and horizontal transposes
  \begin{equation}
    \Donevector[u]{Gy_0}{y_1}
    \leftrightarrow
    \Donevectorop[\overline{u}]{\ra{G}y_1}{y_0}
    \quad\text{and}\quad
    \Donecovector[\overline{v}]{\comap{y}_1}{\comap{y}_0\ra{G}}
    \leftrightarrow
    \Donecovectorop[v]{\comap{y}_0}{\comap{y}_1G},
  \end{equation}
  where we have added the application of the gluing functor
  in the matrix to make the effect of the transposition more apparent
  (usually we would just write something like
  \((\comap{y}_0\to \comap{y}_1)\)).
\end{con}

\begin{lem}
  For each \(\X,\Z:\extwoaddcat\),
  we have a commutative square
  \begin{equation}
    \begin{tikzcd}
      \extwoaddcat(\X,\Y_0\laxplus[G]\Y_1)
      \times
      \extwoaddcat(\Y_0\oplaxplus[G]\Y_1,\Z)
      \ar[r]
      \ar[d,"{\simeq}"]
      &
      \extwoaddcat(\X,\Z)
      \ar[d,equal]
      \\
      \extwoaddcat(\X,\Y_1\oplaxplus[\ra{G}]\Y_0)
      \times
      \extwoaddcat(\Y_1\laxplus[\ra{G}]\Y_0,\Z)
      \ar[r]
      &
      \extwoaddcat(\X,\Z)
    \end{tikzcd}
  \end{equation}
  where the horizontal maps are the oplax-lax and lax-oplax matrix multiplication,
  respectively,
  and the left vertical map is the equivalence of
  \Cref{con:matrix-adjoints}.
\end{lem}

\begin{proof}
  For each
  \begin{equation}
    \Donecovectorop[v]{\comap{y}_0}{\comap{y}_1}
    \colon \Y_0\oplaxplus[G]\Y_1\to \Z
    \quad\text{and}\quad
    \Donevector[u]{y_0}{y_1}\colon\X\to\Y_0\laxplus[G]\Y_1,
  \end{equation}
  the two different row-column products
  are the cofiber in \(\extwoaddcat(\X,\Z)\) of the two composite maps
  \begin{equation}
    \comap{y}_0y_0\xrightarrow{vy_0}\comap{y}_1Gy_0
    \xrightarrow{\comap{y}_1u}
    \comap{y}_1y_1
    \quad\text{and}\quad
    \comap{y}_0y_0
    \xrightarrow{\comap{y}_0\overline{u}}
    \comap{y}_0\ra{G}y_1
    \xrightarrow{\overline{v}y_1}
    \comap{y}_1y_1.
  \end{equation}
  A straightforward computation
  using the triangle identities for \(G\dashv \ra{G}\)
  shows that these two maps are in   canonically identified;
  hence so are their cofibers.
\end{proof}

\section{Chain complexes and chain maps}
\label{subsec:chain-maps}

Throughout this section, let \(\extwoaddcat\)
be a finitely lax additive \((\infty,2)\)-category.

Let \(\ZZ=(\ZZ,\leq)\) be the standard poset of integers.
A chain complex in \(\extwoaddcat\) is a functor \(\ZZ^\op\to\extwoaddcat\),
depicted as
\begin{equation}
  \begin{tikzcd}
    \dots\ar[r,"\alpha"]
    &
    \A_2\ar[r,"\alpha"]
    &
    \A_1\ar[r,"\alpha"]
    &
    \A_0\ar[r,"\alpha"]
    &
    \dots,
  \end{tikzcd}
\end{equation}
with the conditions that each \(\alpha^k\)
is a zero object in \(\extwoaddcat(\A_n,\A_{n-1-k})\) for \(k\geq 2\).

There are various notions of chain maps,
corresponding to different notions of natural transformations
of diagrams \(\ZZ^\op\to\extwoaddcat\)
in the \(2\)-categorical context
(see also \S\ref{sec:adjoints-in-diagrams}).
\begin{itemize}
\item
  A \emph{chain map} (without further qualifier)
  is a commutative diagram of the form
  \begin{equation}
    \label{eq:a-chain-map}
    \begin{tikzcd}
      \dots\ar[r,"\alpha"]
      &
      \A_2\ar[r,"\alpha"]
      \ar[d,"f_2" description]
      &
      \A_1\ar[r,"\alpha"]
      \ar[d,"f_1" description]
      &
      \A_0\ar[r,"\alpha"]
      \ar[d,"f_0" description]
      &
      \dots
      \\
      \dots\ar[r,"\beta"']
      &
      \B_2\ar[r,"\beta"']
      \ar[ur,"\simeq",phantom]
      &
      \B_1\ar[r,"\beta"']
      \ar[ur,"\simeq",phantom]
      &
      \B_0\ar[r,"\beta"']
      &
      \dots
    \end{tikzcd}
  \end{equation}
  Chain complexes and chain maps in \(\extwoaddcat\)
  assemble into an \((\infty,2)\)-category \(\CChA\),
  defined as a full sub-2-category of \(\FUN(\ZZ^\op,\extwoaddcat)\).
\item
  A \emph{lax chain map} is a diagram of the form,
  commuting only up to possibly noninvertible \(2\)-cells.
  \begin{equation}
    \label{eq:a-lax-chain-map}
    \begin{tikzcd}
      \dots\ar[r,"\alpha"]
      &
      \A_2\ar[r,"\alpha"]
      \ar[d,"f_2" description]
      &
      \A_1\ar[r,"\alpha"]
      \ar[d,"f_1" description]
      &
      \A_0\ar[r,"\alpha"]
      \ar[d,"f_0" description]
      &
      \dots
      \\
      \dots\ar[r,"\beta"']
      &
      \B_2\ar[r,"\beta"']
      \ar[ur, shorten >= 5, shorten <= 5 , Rightarrow]
      &
      \B_1\ar[r,"\beta"']
      \ar[ur, shorten >= 5, shorten <= 5 , Rightarrow]
      &
      \B_0\ar[r,"\beta"']
      &
      \dots
    \end{tikzcd}
  \end{equation}
  Chain complexes and chain maps in \(\extwoaddcat\)
  assemble into an \((\infty,2)\)-category \(\CChAlax\),
  defined as a full sub-2-category of
  \(\FUNlax(\ZZ^\op,\extwoaddcat)\).
\item
  Dually, we define the full sub-2-category 
  \(\CChAoplax\subset\FUNoplax(\ZZ^\op,\extwoaddcat)\)
  of chain complexes and \emph{oplax chain maps},
  which explicitly look as follows:
  \begin{equation}
    \label{eq:an-oplax-chain-map}
    \begin{tikzcd}
      \dots\ar[r,"\alpha"]
      &
      \A_2\ar[r,"\alpha"]
      \ar[d,"f_2" description]
      &
      \A_1\ar[r,"\alpha"]
      \ar[d,"f_1" description]
      &
      \A_0\ar[r,"\alpha"]
      \ar[d,"f_0" description]
      &
      \dots
      \\
      \dots\ar[r,"\beta"']
      &
      \B_2\ar[r,"\beta"']
      \ar[ur, shorten >= 5, shorten <= 5 , Leftarrow]
      &
      \B_1\ar[r,"\beta"']
      \ar[ur, shorten >= 5, shorten <= 5 , Leftarrow]
      &
      \B_0\ar[r,"\beta"']
      &
      \dots
    \end{tikzcd}
  \end{equation}
\end{itemize}

Given two chain complexes \((\A_\bullet,\alpha)\) and \((\B_\bullet,\beta)\),
we write
\begin{equation}
  \Maplax[0]{\A}{\B}\hookleftarrow\MapAB\hookrightarrow\Mapoplax[0]{\A}{\B}
\end{equation}
for the three corresponding \(\infty\)-categories of
lax chain maps, chain maps and oplax chain maps \(\A\to \B\).
More generally, we write
\begin{equation}
  \Maplax[\loffset]{\A_\bullet}{\B_\bullet}\coloneqq
  \Maplax[0]{\A_\bullet}{\B_{\loffset+\bullet}}
  \quad\text{and}\quad
  \Mapoplax[\loffset]{\A_\bullet}{\B_\bullet}\coloneqq
  \Mapoplax[0]{\A_\bullet}{\B_{\loffset+\bullet}}
\end{equation}
for the (stable) \(\infty\)-category of 
lax/oplax degree-\(\loffset\)-maps from \(\A\) to \(\B\).
Abstractly, these various \(\infty\)-categories
are just the hom-categories in the \((\infty,2)\)-categories
\(\FUN(\ZZ^\op,\extwoaddcat)\),
\(\FUNlax(\ZZ^\op,\extwoaddcat)\) and
\(\FUNoplax(\ZZ^\op,\extwoaddcat)\).
For us, a more useful description/definition will be
as certain sections of certain tautological fibrations.

For the rest of this section, fix two chain complexes
\((\A_\bullet,\alpha)\) and \((\B_\bullet,\beta)\)
and an integer \(\loffset\in\ZZ\).

\begin{con}
  \label{con:cats-of-chain-maps}
  Consider the functor
  \begin{equation}
    \ZZ^\op\times\ZZ\xrightarrow{\B\times \A^\op} \extwoaddcat\times\extwoaddcat^\op
    \xrightarrow{\extwoaddcat^\op(-,-)}\st;
    \quad
    (m,n) \mapsto \extwoaddcat(\A_n,\B_m)
  \end{equation}
  and its two mixed Grothendieck constructions
  \begin{equation}
    \label{eq:ch-hom-fibration}
    q:\mixedGroth{m:\ZZ^\op}{n:\ZZ}\extwoaddcat(\A_n,\B_m)
    \to \ZZ\times \ZZ
    \quad\text{and}\quad
    q':\mixedGroth{n:\ZZ}{m:\ZZ^\op}\extwoaddcat(\A_n,\B_m)
    \to \ZZ^\op\times\ZZ^\op,
  \end{equation}
  (contravariant,covariant) and
  (covariant,contravariant), respectively.
  We can identify oplax and lax chain maps \(\A\to\B\)
  with sections of \(q\) and \(q'\) on the diagonal.
  More precisely, we define
  \begin{equation}
    \Mapoplax[\loffset]{\A}{\B}
    \coloneqq 
    \Fun_{\ZZ\times\ZZ}
    \left(
      \ZZ(\loffset),
      \mixedGroth{m:\ZZ^\op}{n:\ZZ}\extwoaddcat(\A_n,\B_m)
    \right)
  \end{equation}
  and
  \begin{equation}
    \Maplax[\loffset]{\A}{\B}
    \coloneqq
    \Fun_{\ZZ^\op\times\ZZ^\op}
    \left(
      \ZZ^\op(\loffset),
      \mixedGroth{n:\ZZ}{m:\ZZ^\op}\extwoaddcat(\A_n,\B_m)
    \right)
  \end{equation}
  where
  \begin{equation}
    \ZZ(\loffset)\coloneqq\setP{(n+\loffset,n)}{m=n+\loffset}\subset \ZZ\times\ZZ
  \end{equation}
  is the \(\loffset\)-shifted diagonal.
  Concretely, such a section consists of
  objects \(f_n:\extwoaddcat(\A_n,\B_{n+\loffset})\)
  and morphisms
  \begin{equation}
    f_n\to f_{n+1}
    \quad\text{or}\quad
    f_{n+1}\to f_{n}
  \end{equation}
  in the corresponding Grothendieck construction,
  amounting to morphisms
  \begin{equation}
    \label{eq:2-cells-in-section}
    f_{n-1}\alpha\to \beta f_{n}
    \quad\text{or}\quad
    \beta f_{n}\to f_{n-1}\alpha,
  \end{equation}
  in \(\extwoaddcat(\A_{n},\B_{n+\loffset-1})\)
  respectively.
  We say that \(f_\bullet\) is an \emph{oplax} or \emph{lax}
  \emph{chain map of degree \(\loffset\)}.

  The full subcategories of
  \(\Mapoplax[0]{\A}{\B}\) and \(\Maplax[0]{\A}{\B}\)
  spanned by those sections where the corresponding maps 
  \eqref{eq:2-cells-in-section}
  are equivalences
  are canonically equivalent to each other by passing to inverses;
  we define this common full subcategory to be
  \(\MapAB\);
  it consists of the \emph{degree}-\(\loffset\)-\emph{chain maps}.
\end{con}

\begin{rem}
  By the standing assumption that \(\exaddtwocat\) is finitely lax additive,
  the diagram \(\exaddtwocat(\A_\bullet,\B_\bullet)\colon \ZZ^\op\times\ZZ\to \st\)
  takes values in stable \(\infty\)-categories and exact functors,
  hence any of its associated categories of sections will again be stable
  and restriction functors between them will be exact
  (see also \Cref{ex:lax-lim-cat} \ref{it:ex:lax-limit-stable-cats}).
  
  This includes the \(\infty\)-categories
  \(\Mapoplax[\loffset]\A\B\)
  and
  \(\Maplax[\loffset]\A\B\)
  of shifted diagonal sections
  as well as other variants defined below,
  such as the shifted upper triangular sections of
  \Cref{con:mapoplaxgeq-ut-sections}.
  Moreover, the full subcategory \(\MapAB\hookrightarrow\Maplax[0]\A\B\)
  (defined by the condition that each of the maps \(\beta f_{n+1}\to f_n\alpha\)
  is an equivalence)
  manifestly contains the zero section and is closed under fibers and cofibers;
  hence it is a stable subcategory.
\end{rem}

\begin{rem}
  We shall not unravel the definition of \(\FUNlax\) and \(\FUNoplax\)
  and show that the mapping categories therein do indeed agree with
  the \(\infty\)-categories constructed in \Cref{con:cats-of-chain-maps}.
  For the purpose of this paper,
  the reader may take this construction as the definition.
\end{rem}

It will be useful to study more general sections of the
fibrations \eqref{eq:ch-hom-fibration}.

\begin{con}
  \label{con:mapoplaxgeq-ut-sections}
  Denote by \(\Zgeq[\loffset]\hookrightarrow\ZZ\times\ZZ\)
  the full shifted triangular subposet
  of those \((m,n)\) satisfying \(m\geq n+\loffset\).
  We write
  \begin{equation}
    \Mapoplaxgeq[\loffset]\A\B
    \coloneqq
    \Fun_{\ZZ\times\ZZ}\left(\Zgeq[\loffset],
      \mixedGroth{m:\ZZ^\op}{n:\ZZ}\extwoaddcat(\A_n,\B_m)
    \right).
  \end{equation}
  for the (stable) \(\infty\)-category of shifted upper triangular sections;
  see \Cref{fig:section-in-Map-geq} for a depiction of such sections.
  For each \(\loffset\) we have the obvious (exact) restriction functors
  \begin{equation}
    \MapoplaxAB[\loffset]
    \xleftarrow{\resk}\MapoplaxgeqAB[\loffset]
    \xrightarrow{\reskgeq}\MapoplaxgeqAB[\loffset+1].
  \end{equation}
\end{con}

\begin{rem}
  \label{fig:section-in-Map-geq}
  A section \(F=(F_{m,n})_{m\geq \loffset +n}\) as in
  \Cref{con:mapoplaxgeq-ut-sections}
  with \(f^r_n=F_{r+n,n}\)
  can be depicted as follows:
  \begin{equation}
    \label{eq:section-in-Map-geq}
    \begin{tikzcd}
      &\cdots\ar[r,"\alpha"]
      &\A_2\ar[r,"\alpha"]&\A_1\ar[r,"\alpha"]&\A_0\ar[r,"\alpha"]&\A_{-1}\ar[r,"\alpha"]
      &\cdots\\
      \vdots\ar[d,"\beta"]&\ddots\ar[from=dr,]\\
      \B_{ \loffset+2 }\ar[d,"\beta"]
      && f^{\loffset}_{2}\ar[from=dr,] & f^{\loffset+1}_{1} \ar[l]\ar[from=d]&
      f^{\loffset+2}_0\ar[l]\ar[from=d] & f^{\loffset+3}_{-1}\ar[l]\ar[from=d]\\
      \B_{\loffset+1 }\ar[d,"\beta"]
      &&& f^{\loffset}_{1}\ar[from=dr,] & f^{\loffset+1}_{0}\ar[l]\ar[from=d] & f^{\loffset+2}_{-1}\ar[l]\ar[from=d]\\
      \B_{\loffset }\ar[d,"\beta"]
      &&&& f^{\loffset}_{0}\ar[from=dr,] & f^{\loffset+1}_{{-1}}\ar[l]\ar[from=d]\\
      \B_{\loffset-1} \ar[d,"\beta"]
      &&&&& f^{\loffset}_{{-1}}\ar[from=rd,]\\
      \vdots
      &&&&&&\ddots
    \end{tikzcd}
  \end{equation}
  The complexes \(\A_\bullet\) and \(\B_\bullet\) are drawn for reference
  to indicate how the section \(F\) spreads across the fibers of the fibration.

  Note that the arrows and squares in the diagram
  \eqref{eq:section-in-Map-geq}
  take place in a Grothendieck construction,
  so that one needs to suitably postcompose with \(\beta\)
  or precompose with \(\alpha\) to obtain genuine arrows or squares in the hom-categories
  \(\extwoaddcat(\A_n,\B_{n+r})\).
  Explicitly unpacking the data of such a section \(F\), we see that
  for each \(r\geq \loffset\) and \(n\in\ZZ\) it contains
  \begin{itemize}
  \item
    an object \(f_n^r\) in \(\exaddtwocat(\A_n,\B_{n+r})\),
  \item
    arrows \(\eta^r_n\colon f_n^r\to \beta f_n^{r+1}\) (vertical) and
    \(\epsilon^r_n\colon f_{n-1}^{r+1}\alpha\to f_n^r\) (horizontal)
    in \(\exaddtwocat(\A_n,\B_{n+r})\),
  \item
    an arrow
    \(f_n^r\alpha \to \beta f_{n+1}^r\) (diagonal)
    in \(\exaddtwocat(\A_{n+1},\B_{n+r})\),
  \item
    a commutative square 
    \begin{equation}
      \begin{tikzcd}[column sep=huge]
        \beta f_{n}^r
        &
        \beta f_{n-1}^{r+1}\alpha
        \ar[l,"\beta \epsilon^r_{n}"']
        \\
        f_{n}^{r-1}
        \ar[u,"\eta^{r-1}_{n}"]
        &
        f_{n-1}^r\alpha
        \ar[u, "\eta^r_{n-1} \alpha"']
        \ar[l, "\epsilon^{r-1}_{n}"]
        \ar[ul]
      \end{tikzcd}
    \end{equation}
    in \(\exaddtwocat(\A_{n},\B_{n+r-1})\)
    when \(r>\loffset\),
    and just the upper right triangle when \(r=\loffset\)
    (because \(f^{\loffset-1}_n\) is not defined).
  \end{itemize}
  Observe how the diagonal part
  \(f_\bullet^\loffset = F\resk : \Mapoplax[\loffset]\A\B\)
  precisely encodes
  the datum of an oplax degree-\(\loffset\) chain map
  as in \Cref{con:cats-of-chain-maps}.
\end{rem}

The following lemma states that we can ``crop''
redundant zeroes in a section
\(f:\MapoplaxgeqAB[\loffset]\).

\begin{lem}
  \label{lem:more-zeros-are-redundant}
  Denote by
  \begin{equation}
    U_\loffset^r\coloneqq \setP{(m,n)}{\loffset\leq m-n \leq\loffset+r}
    \subset \Zgeq[\loffset]
  \end{equation}
  the \(\loffset\)-shifted diagonal strip of width \(r\).
  The canonical restriction functor along
  \(U_\loffset^r\hookrightarrow \Zgeq[\loffset]\)
  induces an equivalence
  \begin{equation}
    \label{eq:restriction-to-strip-0}
    \MapoplaxgeqAB[\loffset]_{\reskgeq[\loffset+r]=0}
    \xrightarrow{\simeq}
    \Fun_{\ZZ\times\ZZ}
    \left(
      U^r_\loffset, \mixedGroth{m:\ZZ^\op}{n:\ZZ}{\extwoaddcat(\A_n,\B_m)}
    \right)_{\reskgeq[\loffset+r]=0},
  \end{equation}
  where on both sides we are only considering those sections
  which are zero on the \(r\)-th off diagonal and beyond.
\end{lem}

\begin{proof}
  First of all, we claim that the restriction functor
  \eqref{eq:restriction-to-strip-0}
  admits a fully faithful left adjoint given by left \(q\)-Kan extension
  (\(q\) is the fibration \eqref{eq:ch-hom-fibration}).
  The pointwise \(q\)-Kan extension formula trivializes,
  since for each \(m\geq n+\loffset+r\)
  the overcategory \(U^r_k/(m,n)\)
  has a terminal object given by the vertical edge
  \((n+\loffset+r,n)\to (m,n)\).
  We thus only have to argue that there are sufficiently many coCartesian edges
  over these vertical edges \((n+\loffset+r,n)\to(m,n)\).
  Since we are, by definition,
  only considering sections whose value at \((n+\loffset+r,n)\) is zero,
  this is automatic;
  the resulting Kan extended diagram is zero on \(\Zgeq[\loffset{+}r]\).
  The result follows since by construction the essential image
  of this left \(q\)-Kan extension is precisely
  \(\MapoplaxgeqAB[\loffset]_{\reskgeq[\loffset{+}r]=0}\).
\end{proof}

\begin{rem}
  Even when the \(r\)-th off-diagonal is zero,
  we cannot crop the diagram any further without losing information.
  In other words, the restriction functor
  \begin{equation}
    \label{eq:restriction-to-strip-no-0}
    \MapoplaxgeqAB[\loffset]_{\reskgeq[\loffset+r]=0}
    \xrightarrow{\simeq}
    \Fun_{\ZZ\times\ZZ}
    \left(
      U^{r-1}_\loffset, \mixedGroth{m:\ZZ^\op}{n:\ZZ}{\extwoaddcat(\A_n,\B_m)}
    \right),
  \end{equation}
  is \emph{not} typically an equivalence
  because the commutative squares
  \begin{equation}
    \cdsquareOpt
    {\beta f^{r-1}_n}{0}
    {f^{r-2}_n}{f^{r-1}_{n-1}\alpha}
    {leftarrow}
    {leftarrow}
    {leftarrow}
    {leftarrow}
  \end{equation}
  at the edge of the strip \(U^{r}_\loffset\)
  carry more data than just the composable arrows
  \begin{equation}
    f^{r-1}_{n-1}\alpha\to f^{r-2}_n\to \beta f^{r-1}_n,
  \end{equation}
  namely a trivialization of their composite.
\end{rem}

\begin{rem}
  \label{rem:data-of-section-f-h-0}
  In the special case \(r=2\),
  \Cref{lem:more-zeros-are-redundant} says that a section
  \(f:\MapoplaxgeqAB[\loffset]\)
  satisfying \(f\reskgeq[\loffset{+}2]=0\)
  amounts to the following data:
  \begin{itemize}
  \item
    objects \(f_n\coloneqq f^\loffset_n:\extwoaddcat(\A_n,\B_{\loffset+n})\),
  \item
    objects \(h_n\coloneqq f^{\loffset+1}_n:\extwoaddcat(\A_n,\B_{\loffset+n+1})\),
  \item
    commutative squares
    \begin{equation}
      \label{eq:square-in-data-f-h-0}
      \cdsquareOpt
      {\beta h_n}{0}
      {f_n}{h_{n-1}\alpha}
      {leftarrow}
      {leftarrow}
      {leftarrow}
      {leftarrow}
    \end{equation}
    in \(\extwoaddcat(\A_n,\B_{\loffset+n})\).
  \end{itemize}
\end{rem}

\begin{lem}
  \label{lem:trivialized-chain-maps}
  The canonical evaluation maps
  at the individual \(\loffset\)-shifted diagonal entries
  assemble into an equivalence
  \begin{equation}
    \MapoplaxgeqAB[\loffset]_{\reskgeq[\loffset+1] = 0}
    \xrightarrow{\simeq}
    \prod_{n\in\ZZ}\extwoaddcat(\A_n,\B_{\loffset+n}),
    \quad
    f^\bullet_\bullet\mapsto (f^\loffset_n)_{n\in\ZZ}
  \end{equation}
  of (stable) \(\infty\)-categories.
  Here the left hand side denotes the kernel of the restriction functor
  \begin{equation}
    \MapoplaxgeqAB[\loffset]
    \xrightarrow{\reskgeq}
    \MapoplaxgeqAB[\loffset+1],
  \end{equation}
  i.e., the full subcategory of those shifted upper triangular sections
  that vanish on the first off-diagonal and beyond.
\end{lem}

\begin{proof}
  By \Cref{lem:more-zeros-are-redundant} (with \(r=1\)),
  we may restrict our sections to the strip
  \(U^{1}_\loffset\subset \Zgeq[\loffset]\),
  which, as a poset, is simply isomorphic to \(\ZZ\)
  via \((m,n)\mapsto m+n\).
  Therefore a diagram of shape \(U^{1}_\loffset\)
  just amounts to a sequence of objects and arrows.
  If such a diagram is zero on odd-indexed objects
  \(2n+1\triangleq(n+1,n)\),
  then all arrows are uniquely determined
  and the only relevant data are the values
  at the even-indexed objects \(2n\triangleq(n,n)\).
\end{proof}

\begin{lem}
  \label{lem:reskgeq-has-adjoints}
  The restriction functor
  \begin{equation}
    \MapoplaxgeqAB[\loffset]\to\MapoplaxgeqAB[\loffset+1]
  \end{equation}
  admits a left adjoint \(j\) and a right adjoint \(j'\), both fully faithful,
  given by left and right \(q\)-Kan extension, respectively.
  A section \eqref{eq:section-in-Map-geq}
  lies in the essential image
  of \(j\) / \(j'\)
  if and only each leftmost horizontal / bottommost vertical
  edge is coCartesian / Cartesian,
  i.e.\ induces an equivalence
  \begin{equation}
    f^{\loffset+1}_{n-1}\alpha
    \xrightarrow{\simeq}
    f^{\loffset}_{n}
  \quad\quad
    \text{/}
    \quad\quad
    f^{\loffset}_{n}
    \xrightarrow{\simeq}
    \beta f^{\loffset+1}_{n}.
  \end{equation}
\end{lem}

\begin{proof}
  The pointwise left \(q\)-Kan extension formula at
  \((n+\loffset,n)\) along the inclusion
  \(\Zgeq[\loffset{+}1]\hookrightarrow\Zgeq[\loffset]\)
  trivializes,
  since the overcategory
  \(\Zgeq[n{+}\loffset{+}1]/(n+\loffset,n)\)
  has a terminal object \((n+\loffset,n-1)\).
  Therefore the desired left \(q\)-Kan extension exists
  if and only if each horizontal edge
  \((n+\loffset,n)\leftarrow(n+\loffset,n-1)\)
  admits a coCartesian lift.
  Since the fibration
  \begin{equation}
    q\colon \mixedGroth{m:\ZZ^\op}{n:\ZZ}\extwoaddcat(\A_n,\B_m)
    \longrightarrow \ZZ\times \ZZ
  \end{equation}
  is (by construction) coCartesian in the second variable,
  this is always the case.

  The argument for the right adjoint is dual.
\end{proof}

Going forward, it will be convenient to reformulate such statements
using recollements of stable \(\infty\)-categories.
See \Cref{sec:recollements} for a brief summary of the theory
as we will use it without further explicit mention.

\begin{cor}
  \label{cor:reskgeq-recollement}
  The restriction functor \(\reskgeq\) is part of a recollement
  \begin{equation}
    \label{eq:recollement-kgeq}
    \cdrecollementcompl
    {
      \prod_{n\in\ZZ}\extwoaddcat(\A_n,\B_{\loffset+n})
    }
    {
      \MapoplaxgeqAB[\loffset]
    }
    {
      \MapoplaxgeqAB[\loffset{+}1]
    }
    {}
    {\reskgeq[\loffset{+}1]}
    {j}
    {}
    {j'}
    {}
  \end{equation}
  of stable \(\infty\)-categories
  with gluing functor
  \(
  f^\bullet_\bullet\mapsto
  \left(
    \fib(f_{n-1}^{\loffset+1}\alpha\to\beta f^{\loffset+1}_{n})
  \right)_n
  \).
\end{cor}

\begin{proof}
  \Cref{lem:reskgeq-has-adjoints}
  provides the two fully faithful adjoints \(j\) and \(j'\)
  of the functor \(\reskgeq[\loffset+1]\)
  as left and right \(q\)-Kan extension, respectively.
  Since the kernel of this functor is identified with
  \(\prod_{n\in \ZZ}\exaddtwocat(\A_n,\B_{k+n})\)
  by \Cref{lem:trivialized-chain-maps},
  this determines the desired recollement.

  It remains to compute the gluing functor.
  From the pointwise formulas for the relative Kan extensions
  we see that for each \(f^\bullet_\bullet:\MapoplaxgeqAB[\loffset{+}1]\)
  the canonical transformation \(j(f)\to j'(f)\) is given
  on the main diagonal by the structure map
  \begin{equation}
    j(f)^\loffset_n= f^{k+1}_{n-1}\alpha\to\beta f^{\loffset+1}_{n}=j'(f)^\loffset_n
  \end{equation}
  (for \(n\in \ZZ\)); passing to fibers yields the desired formula for the gluing functor.
\end{proof}

\begin{rem}
  Note that neither of the two adjoints in the left half
  of the recollement \eqref{eq:recollement-kgeq}
  are the tautological functor \(f^\bullet_\bullet\mapsto (f^\loffset_n)_{n\in \ZZ}\)
  that evaluates a section
  at the individual entries of the \(\loffset\)-shifted diagonal.
\end{rem}

\begin{rem}
  We can think of 
  \(\MapoplaxgeqAB[\loffset]_{\reskgeq=0}\)
  as the \(\infty\)-category of degree-\(\loffset\) chain maps
  \(f\colon \A_\bullet\to\B_{\loffset+\bullet}\)
  with trivialized structure map \(f\alpha\to\beta f\).
  Note that this is not a full subcategory of \(\MapoplaxAB[k]\).
  Indeed, the restriction functor
  \begin{equation}
    \MapoplaxgeqAB[\loffset]_{\reskgeq=0}
    \to \MapoplaxAB[\loffset]
  \end{equation}
  which forgets the trivialization is neither full nor faithful.
\end{rem}

The restriction functor to the diagonal does
not, in general, have analogous adjoints.
This does happen in the special case where
the differentials of the chain complexes
\((\A_\bullet,\alpha)\) and/or \((\B_\bullet,\beta)\)
have left adjoints:

\begin{lem}
  \label{lem:resk-has-adjoints}
  Consider the restriction functor
  \begin{equation}
    \label{eq:res-kgeq-diag}
    \resk\colon\MapoplaxgeqAB[\loffset]\longrightarrow
    \MapoplaxAB[\loffset]
  \end{equation}
  \begin{enumerate}
  \item
    Assume that each differential \(\beta\) has a left adjoint.
    Then this restriction functor has a fully faithful left adjoint
    given by relative left Kan extension.
    Explicitly it is given by
    \begin{equation}
      f^{\loffset+1}_n\coloneqq \la{\beta}f^{\loffset}_n
      \quad\text{and}\quad
      f^{r}_n=0 \text{ for } r\geq \loffset +2
    \end{equation}
    with the non-trivial vertical arrows amounting to the units
    \(f^\loffset_n\to\beta\la{\beta}f^\loffset_n\) of the adjunction.
  \item
    \label{it:alpha-left-adjoint}
    Assume that each differential \(\alpha\) has a left adjoint.
    Then this restriction functor has a fully faithful right adjoint
    given by relative right Kan extension.
    Explicitly it is given by
    \begin{equation}
      f^{\loffset+1}_n\coloneqq f^{\loffset}_{n+1}\la{\alpha}
      \quad\text{and}\quad
      f^{r}_n=0 \text{ for } r\geq \loffset +2
    \end{equation}
    with the non-trivial horizontal arrows amounting to the counits
    \(f^{\loffset}_{n+1}\la{\alpha}\alpha\to f^{\loffset}_{n+1} \)
    of the adjunction.
   \end{enumerate}
\end{lem}

\begin{proof}
  The two statements are dual; we focus on \ref{it:alpha-left-adjoint}.

  We observe that the relevant undercategories
  \((m,n)/\ZZ(\loffset)\)
  (for \((m,n):\Zgeq[\loffset]\))
  have an initial object \((m,m-\loffset)\).
  Therefore the desired pointwise right \(q\)-Kan extension exists
  if we can guarantee that each horizontal edge
  \((m,m-\loffset)\to (m,n)\) has a Cartesian lift.
  In general, the fibration
  \begin{equation}
    q\colon \mixedGroth{m:\ZZ^\op}{n:\ZZ}{\extwoaddcat(\A_n,\B_m)}
    \longrightarrow
    \ZZ\times \ZZ
  \end{equation}
  is only Cartesian in the first variable, not in the second.
  Being Cartesian in the second variable amounts to
  each \(\extwoaddcat(\alpha,\B_m)\) having a right adjoint
  which is guaranteed because each
  \(\alpha\colon \A_n\to\A_{n-1}\) has a left adjoint by assumption.
  The explicit formulas are an immediate consequence
  of this pointwise construction
  using \(\la{\beta}\la{\beta}=0\)
  to obtain the vanishing beyond the first off-diagonal.
\end{proof}

\begin{lem}
  \label{lem:kernel-of-both-restrictions}
  \label{lem:inclusion-of-oplax}
  There is an equivalence, canonical up to shift, between
  \begin{itemize}
  \item
    the full subcategory
    \begin{equation}
      \setP{f}{\forall r\neq \loffset+1 : f^r_\bullet=0}
      \subset
      \MapoplaxgeqAB[\loffset]
    \end{equation}
    of those sections \(f\) which are non-zero only on the first off-diagonal and
  \item
    the \(\infty\)-category \(\MapoplaxAB[\loffset{+}1]\)
    of oplax degree-\((\loffset{+}1)\) chain maps.
  \end{itemize}
  Explicitly it sends a section \(f^\bullet_\bullet\)
  to a chain map with components
  \(g_n\coloneqq f^{\loffset+1}_n[-n] : \extwoaddcat(\A_n,\B_{\loffset+n+1})\).
\end{lem}

\begin{rem}
  Note that the equivalence of \Cref{lem:kernel-of-both-restrictions}
  is \emph{not} induced by the obvious restriction functor
  \begin{equation}
    \MapoplaxgeqAB[\loffset]\xrightarrow{\resk[\loffset{+1}]} \MapoplaxAB[\loffset{+}1]
  \end{equation}
  which, when restricted to 
  \(\setP{f}{\forall r\neq \loffset+1 : f^r_\bullet=0}\)
  only hits oplax chain maps with trivial structure map.
\end{rem}

\begin{proof}
  According to \Cref{rem:data-of-section-f-h-0},
  the data of a section
  \(f : \MapoplaxAB[\loffset]\)
  with
  \(f\reskgeq[\loffset{+}2] = 0\)
  and
  \(f^\loffset_\bullet=0\)
  amounts to
  \begin{itemize}
  \item
    1-morphisms \(h_n\coloneqq f^{\loffset+1}_n : \extwoaddcat(\A_n,\B_{\loffset+n+1})\)
  \item
    and commutative squares
    \begin{equation}
      \begin{tikzcd}
        \beta h_n&0\ar[l]\\
        0\ar[u]&h_{n-1}\alpha\ar[l]\ar[u]
      \end{tikzcd}
    \end{equation}
    in \(\extwoaddcat(\A_{n},\B_{\loffset+n})\)
    which amount precisely to morphisms
    \(\phi_{n}\colon h_{n-1}[-n+1]\alpha\to\beta h_{n}[-n]\).
  \end{itemize}
  Thus setting \(g_n\coloneqq h_n[-n]\),
  this is precisely the data of an oplax degree-\((\loffset{+}1)\) map
  \(g=(g_\bullet,\phi_\bullet):\MapoplaxAB[\loffset{+}1]\).
\end{proof}

\begin{prop}
  \label{prop:recollement-a-b-adjoints}
  Assume that all differentials \(\alpha\) and \(\beta\)
  have left adjoints.
  Then the restriction functor
  \begin{equation}
    \label{eq:kgeq-restricted-k-restriction}
    \resk\colon \MapoplaxgeqAB[\loffset]_{\reskgeq[\loffset{+}2]=0}
    \longrightarrow\MapoplaxAB[\loffset]
  \end{equation}
  is part of a recollement
  \begin{equation}
    \label{eq:recollement-resk}
    \cdrecollementcompl
    {\MapoplaxAB[\loffset{+}1]}
    {\MapoplaxgeqAB[\loffset]_{\reskgeq[\loffset{+}2]=0}}
    {\MapoplaxAB[\loffset]}
    {i}{\resk}{j}{}{j'}{};
  \end{equation}
  with gluing functor
  \begin{equation}
    \label{eq:formula-gluing-functor-resk}
    \rho\colon
    f_\bullet \mapsto
    \left(
      \fib(\la{\beta}f_n\to f_{n+1}\la{\alpha})[-n]
    \right)_n.
  \end{equation}
  In particular, we have the dashed equivalence
  of (stable) \(\infty\)-categories.
  \begin{equation}
    \label{eq:Mapoplaxgeq-gluing}
    \begin{tikzcd}
      \MapoplaxgeqAB[\loffset]_{\reskgeq[\loffset{+}2]=0}
      \ar[r,dashed,leftrightarrow,"\simeq"]
      \ar[d,"{\resk}"]
      &
      \laxlima{\MapoplaxAB[\loffset]}{\rho}{\MapoplaxAB[\loffset{+}1]}
      \ar[dl,"p_0"]
      \\
      \MapoplaxAB[\loffset]
    \end{tikzcd}
  \end{equation}
\end{prop}
\begin{proof}
  By \Cref{lem:resk-has-adjoints},
  the restriction functor
  \eqref{eq:kgeq-restricted-k-restriction}
  has adjoints \(j\) and \(j'\)
  given by relative left and right Kan extension.
  Using \Cref{lem:kernel-of-both-restrictions}
  to identify the kernel then yields the recollement 
  \eqref{eq:recollement-resk}
  and the induced equivalence 
  \eqref{eq:Mapoplaxgeq-gluing}
  by the general theory.

  From the explicit construction in
  \Cref{lem:kernel-of-both-restrictions}
  it follows that the canonical transformation
  \(j\to j'\) between the two adjoints
  is given explicitly at \(f:\MapoplaxAB[\loffset]\) by
  the canonical mate
  \begin{equation}
    j(f)_n^{\loffset+1}=\la{\beta}f_n\to f_{n+1}\la{\alpha}=j'(f)_n^{\loffset+1}
  \end{equation}
  on the first off-diagonal;
  it is an equivalence
  (\(f_n\xrightarrow{=}f_n\) or \(0\xrightarrow{=} 0\)) everywhere else.
  The gluing functor
  \begin{equation}
    \MapoplaxAB[\loffset]\to\Ker(\resk)
  \end{equation}
  is given by the fiber of this transformation,
  therefore yields the desired formula \eqref{eq:formula-gluing-functor-resk}
  under the identification of \Cref{lem:inclusion-of-oplax}.
\end{proof}

\begin{defi}
  We denote by
  \begin{equation}
    \MapoplaxgeqABex[\loffset] \subset \MapoplaxgeqAB[\loffset]
  \end{equation}
  the full subcategory spanned by those sections \((f^r_n)\)
  such that all the induced squares
  \begin{equation}
    \cdsquareOpt[bicart]
    {\beta f_n^{r+1}}{\beta f_{n-1}^{r+2}\alpha}
    {f_n^{r}}{f_{n-1}^{r+1}\alpha}
    {leftarrow}
    {leftarrow}
    {leftarrow}
    {leftarrow}
  \end{equation}
  in \(\extwoaddcat(\A_n,\B_{r+n})\) are biCartesian
  (for all \(r\geq k\)) and call such sections \emph{exact}.
\end{defi}

\begin{lem}
  \label{lem:inclusion-of-lax}
  There is an equivalence of \(\infty\)-categories between
  \begin{itemize}
  \item
    the full subcategory
    \begin{equation}
      \MapoplaxgeqABex[\loffset]_{\reskgeq[\loffset{+}2]=0}
      \subset \MapoplaxgeqAB[\loffset]
    \end{equation}
    of those sections which are exact
    and vanish beyong the first off-diagonal and
  \item
    the \(\infty\)-category \(\MaplaxAB[\loffset{+}1], \)
    of \emph{lax} degree-\((\loffset{+}1)\) chain maps.
  \end{itemize}
  Explicitly it sends a section \((f^r_n)\) to a chain map with components
  \(g_n\coloneqq f^{\loffset+1}_n[-n]\).
\end{lem}

\begin{proof}
  In \Cref{rem:data-of-section-f-h-0},
  if we restrict to squares \eqref{eq:square-in-data-f-h-0}
  which are biCartesian,
  the data just amounts
  (by rotating the exact triangle forward and shifting by \([-n]\))
  to objects
  \(h_n=f^{\loffset+1}_n:\extwoaddcat(\A_n,\B_{\loffset+n+1})\)
  and maps
  \(\beta h_n[-n]\to h_{n-1}[-n+1]\alpha\) in
  \(\extwoaddcat(\A_n,\B_{\loffset+n})\).
  This is precisely the data of a lax
  degree-\(\loffset{+}1\) chain map \(g\) with components
  \(g_n\coloneqq h_n[-n]\), as desired.
\end{proof}

\begin{rem}
  \Cref{lem:inclusion-of-oplax} and \Cref{lem:inclusion-of-lax}
  explain how the \(\infty\)-category
  \(\MapoplaxgeqAB[\loffset]_{\reskgeq[\loffset{+}2]=0}\)
  contains both the oplax and the lax degree-\((\loffset{+}1)\)
  maps \(\A\to\B\).
  From the explicit constructions it is immediate that these two inclusions
  are compatible,
  in the sense that there is a commutative square
  \begin{equation}
    \cdsquareOpt
    {\Map(\A_\bullet,\B_{\loffset+\bullet+1})}
    {\MapoplaxAB[\loffset{+1}]}
    {\MaplaxAB[\loffset{+1}]}
    {\MapoplaxgeqAB[\loffset]_{\reskgeq[\loffset{+}2]=0}}
    {hookrightarrow}
    {hookrightarrow}
    {hookrightarrow}
    {hookrightarrow}
  \end{equation}
  and we have
  \begin{equation}
    \Map(\A_\bullet,\B_{\loffset+\bullet+1})
    =
    {\MapoplaxAB[\loffset{+1}]}
    \cap
    {\MaplaxAB[\loffset{+1}]}
  \end{equation}
  as full subcategories of 
  \(\MapoplaxgeqAB[\loffset]_{\reskgeq[\loffset{+}2]=0}\).
\end{rem}

\begin{con}[lax mapping complex]
  \label{con:lax-mapping-complex}
  Let \((f^r_n)=(f,h):\MapoplaxgeqAB[\loffset]\)
  be a section as in \Cref{rem:data-of-section-f-h-0}.
  If each square \eqref{eq:square-in-data-f-h-0} is biCartesian,
  then both of the maps
  \begin{equation}
    f_n\alpha\to \beta h_n \alpha
    \quad\text{and}\quad
    \beta h_n \alpha \to \beta f_{n+1}
  \end{equation}
  are equivalences,
  since their fibers/cofibers are 
  \begin{equation}
    h_{n-1}\alpha\alpha=0
    \quad\text{and}\quad
    \beta\beta h_{n+1}=0,
  \end{equation}
  respectively.
  Therefore the oplax degree-\(\loffset\) chain map \(f=f^\loffset\)
  is an actual chain map \(\A_\bullet\to\B_{\loffset+\bullet}\).
  Therefore the canonical restriction functor
  \begin{equation}
    \resk[\loffset]\colon
    \MapoplaxgeqAB[\loffset]
    \to
    \MapoplaxAB[\loffset]
  \end{equation}
  restricts to a functor
  \begin{equation}
    \delta\colon
    \MaplaxAB[\loffset{+}1]
    \simeq
    \MapoplaxgeqABex[\loffset]_{\reskgeq[\loffset{+}2]=0}
    \xrightarrow{\resk[\loffset]}
    \Map(\A_\bullet,\B_{\loffset+\bullet})
  \end{equation}
  whose kernel is precisely
  \(\Map(\A_\bullet,\B_{\loffset+\bullet+1})\).
  These differentials \(\delta\) assemble to what we call the
  \emph{lax mapping complex} \(\MaplaxAB\):
  \begin{equation}
    \label{eq:lax-mapping-complex}
    \begin{tikzcd}[column sep=0]
      \cdots
      \MaplaxAB[2]
      \ar[dr,"\delta"]
      \ar[rr,dashed,"\delta"]
      &&
      \MaplaxAB[1]
      \ar[rr,dashed,"\delta"]
      \ar[dr,"\delta"]
      &&
      \MaplaxAB[0]
      \cdots
      \\
      &\Map(\A_\bullet,\B_{1+\bullet})
      \ar[ur,hookrightarrow]
      &&\Map(\A_\bullet,\B_{\bullet})
      \ar[ur,hookrightarrow]
    \end{tikzcd}
  \end{equation}
  Unraveling, we get the explicit formula for the differential
  \begin{equation}
    \label{eq:formula-diff-lax}
    \delta(g_\bullet)_n=\fib(\beta g_n\to g_{n-1}\alpha)[n].
  \end{equation}
\end{con}

\begin{rem}
  Assume that the differentials \(\alpha\) and \(\beta\) have right adjoints
  \(\ra{\alpha}\) and \(\ra{\beta}\), respectively.
  Denote by \(\ra{\A}_\bullet\coloneqq(\A_{-\bullet},\ra{\alpha})\)
  and \(\ra{\B}_\bullet\coloneqq(\B_{-\bullet},\ra{\beta})\)
  the chain complex obtained from \(\A\) and \(\B\)
  by passing to right adjoints of the differentials.
  Note that for each \(n\in\NN\)
  there is a tautological equivalence of \(\infty\)-categories
  \begin{equation}
    \begin{tikzcd}[row sep=tiny]
      \Mapoplax[\loffset]{\ra{\A}}{\ra{\B}}
      \ar[r,leftrightarrow,"\simeq"]
      &
      \MaplaxAB[-\loffset]
      \\
      (f_{\bullet}, f\ra{\alpha}\to \ra{\beta}f),
      \ar[leftrightarrow,r]
      &
      (f_{-\bullet}, \beta f\to f\alpha)
    \end{tikzcd}
  \end{equation}
  by noting that both sides are sections
  \begin{equation}
    \begin{tikzcd}
      &
      \mixedGroth{m:\ZZ^\op}{n:\ZZ}{\extwoaddcat(\ra{\A}_n,\ra{\B}_m)}
      \ar[d,"q"]
      \ar[r,equals]
      &
      \mixedGroth{n:\ZZ}{m:\ZZ^\op}{\extwoaddcat(\A_n,\B_m)}
      \ar[d,"q'"]
      \\
      \ZZ(\loffset)
      \ar[r,hookrightarrow]
      \ar[ur,dashed,bend left=20]
      &
      \ZZ\times\ZZ
      \ar[r,"{\cong}",leftrightarrow,"(-1)\cdot"']
      &
      \ZZ^\op\times\ZZ^\op
    \end{tikzcd}
  \end{equation}
  of the same fibration.
  An explicit computation shows that under this equivalence,
  the differential
  \begin{equation}
    \delta\colon \MaplaxAB[\loffset+1]\to\MaplaxAB[\loffset]
  \end{equation}
  of the lax mapping complex \eqref{eq:lax-mapping-complex}
  is identified with the gluing functor
  \begin{equation}
    \rho\colon \Mapoplax[-\loffset-1]{\ra{\A}}{\ra{\B}}
    \to\Maplax[-\loffset]{\ra{\A}}{\ra{\B}}
  \end{equation}
  of \Cref{prop:recollement-a-b-adjoints}
  applied to the chain complexes \(\ra{\A}\) and \(\ra{\B}\).
\end{rem}

Once we have constructed the mapping complex,
we immediately get the corresponding notion of categorified chain homotopy.

\begin{defi}
  \label{defi:lax-null-homotopy}
  Let \(f\colon \A\to\B\) be a chain map.
  A \emph{lax null-homotopy} of \(f\)
  is a lax degree-\(1\) map \(h\colon \MaplaxAB[1]\)
  with \(\delta(h)=f\).
\end{defi}

\begin{rem}
  Clearly, one could dualize \Cref{con:lax-mapping-complex}
  and all the preceding lemmas
  to obtain the \emph{oplax} mapping complex and the resulting notion
  of an \emph{oplax} null-homotopy.
  This is the version which appears
  in \cite[Section~4.6]{CDW23},
  where this oplax mapping complex was constructed
  (in the special case \(\extwocat=\StL_k\)) via the
  product totalization of the canonical double complex
  \(\extwocat(\A_{-\bullet},\B_\bullet)\).
  We shall not give a detailled proof
  that these two different constructions agree;
  this is relatively straightforward by inspection
  of the terms of the complex and the explicit formulas
  for the differential.
\end{rem}

\begin{defi}
  \label{defi:adjointable-square}
  A commutative square
  \begin{equation}
    \cdsquare{\A}{\B}{\C}{\D}{g}{f}{f'}{g'}
  \end{equation}
  in an \((\infty,2)\)-category \(\exaddtwocat\)
  is called \emph{vertically left/right adjointable} if both \(f\) and \(f'\)
  have a left/right adjoint
  and the corresponding canonical mate
  \begin{equation}
    \la{f'}g'\to g\la{f}
    \quad
    \text{/}
    \quad
    g\ra{f}\to\ra{f'}g'
  \end{equation}
  is an equivalence.
  \emph{Horizontally left/right adjointable} is defined analogously
  but with \(g\) and \(g'\) having adjoints.
\end{defi}

For chain maps, we distinguish two types of adjointability conditions:
in the direction of the differentials and
in the direction of the chain map itself.
 
\begin{defi}
  Let \(f\colon (A_\bullet,\alpha)\to(\B_\bullet,\beta)\)
  be a chain map.
  \begin{itemize}
  \item
    We say that \(f\)
    is \emph{\ldiffBC{}} / \emph{\rdiffBC{}}
    if each square in the corresponding diagram \eqref{eq:a-chain-map}
    is horizontally left/right adjointable,
    i.e. all differentials \(\alpha\) and \(\beta\) admit left/right adjoints
    and the canonical mate
    \(\la{\beta}f \to f\la{\alpha}\) / \(f\ra{\alpha}\to\ra{\beta}f\)
    is an equivalence.
  \item
    We say that \(f\) is \emph{left/right adjointable} 
    if each square in the corresponding diagram \eqref{eq:a-chain-map}
    is vertically left/right adjointable,
    i.e., each component \(f_n\) has a left/right adjoint
    and the canonical mate
    \(\la{f}\beta\to \alpha\la{f}\) / \(\alpha\ra{f}\to\ra{f}\beta\)
    is an equivalence.
  \end{itemize}
\end{defi}

\section{The oplax mapping cone construction}
\label{sec:lax-mapping-cone}

Let \(\extwoaddcat\) be a finitely lax additive \((\infty,2)\)-category.
The goal of this section is to construct the oplax mapping cone
\(\oplaxmapcone{f}\)
of a chain map \eqref{eq:a-chain-map}
in \(\extwoaddcat\) by categorifying the usual formula
\begin{equation}
  \label{eq:classic-mapping-cone-diff}
  \mappingcone[n+1]{f} \coloneqq A_{n}\oplus B_{n+1}
  \xrightarrow{
    \ppsmallmatrix{-\alpha&0\\-f&\beta}
  }
  A_{n-1}\oplus B_{n}
  \eqqcolon
  \mappingcone[n]{f}
\end{equation}
for the differential. According to the philosophy outlined in
\S\ref{subsec:rules-of-categorification},
we need additional data to specify the mapping cone complex:
\begin{itemize}
\item
  To construct the terms of the mapping cone complex
  \(\oplaxmapcone{f}\coloneqq\A_{n-1}\laxplus\B_{n}\)
  as a lax bilimit,
  we need to specify 1-morphisms
  \(h\colon\A_{n-1}\to\B_{n}\)
  or \(k\colon \B_n\to \A_{n-1}\)
\item
  We need some suitable 2-categorical data to be able to write down
  the \(\Delta^1\times\Delta^1\)-analog of the differential matrix
  \eqref{eq:classic-mapping-cone-diff}.
\end{itemize}

We will also see that in the presence of sufficient compatible adjoints
to the differentials \(\alpha\), \(\beta\) and/or \(f\),
one can canonically construct such data using the various units/counits
and in this case we recover the fiber and cofiber of $f$ as in
\cite[Section~4.3]{CDW23}.

\begin{defi}
  \label{defi:enhanced-chain-map}
  We denote by
  \begin{equation}
    \MaplhAB\coloneqq \MapAB\times_{\MapoplaxAB[0]}
    \MapoplaxgeqAB[0]_{\reskgeq[2]=0}
  \end{equation}
  the \(\infty\)-category of those sections
  \eqref{eq:section-in-Map-geq} which
  are zero beyond the first off-diagonal
  and restrict to an honest chain map (as opposed to an oplax one)
  on the diagonal.
  Such sections are called
  \emph{\lhenhancedmor{}s}
  of chain complexes and are written
  \(F\colon (\A_\bullet,\alpha_\bullet)\lhenhmor(\B_\bullet,\beta_\bullet)\).
\end{defi}

The mnemonic ``lh'' stands for ``left-horizontal''
and is explained by \Cref{lem:universal-lh-enhancements},
where we construct canonical \lhenhancement{}s
in the presence of left adjoints in the horizontal (=differential) direction.

\begin{rem}
  \Cref{rem:data-of-section-f-h-0} tells us that an \lhenhancedmor{}
  \(F\colon (\A_\bullet,\alpha_\bullet)\lhenhmor(\B_\bullet,\beta_\bullet)\)
  consists of 1-morphisms
  \begin{equation}
    f_n\colon\A_n\to\B_n
    \quad\text{and}\quad
    h_n\colon \A_n\to\B_{n+1}
  \end{equation}
  together with (not necessarily biCartesian) commutative squares
  \begin{equation}
    \label{eq:eta-epsilon=0}
    \begin{tikzcd}
      h_{n-1}\alpha_n\ar[r]\ar[d,"\epsilon_n"]&0\ar[d]\\
      f_n\ar[r,"\eta_n"]&\beta_{n+1} h_{n}
    \end{tikzcd}
  \end{equation}
  in \(\extwoaddcat(\A_n,\B_n)\)
  such that each composite
  \begin{equation}
    \label{eq:epsilon-eta=1}
    f_{n-1}\alpha_n\xrightarrow{\eta_{n-1}\alpha_n}\beta_{n}h_{n-1}\alpha_n
    \xrightarrow{\beta_{n}\epsilon_n}\beta_{n}f_n
  \end{equation}
  is an equivalence
  (i.e., exhibits \(f\colon \A\to\B\) as a chain map).
  We say that \(F=(F,h,\epsilon,\eta)\) is an \emph{\lhenhancement{}}
  of the underlying chain map \(f\colon \A\to \B\).

  We can also depict such an \lhenhancedmor{} of chain complexes as follows
  \begin{equation}
    \label{eq:lhenhmor-as-pasting}
    \begin{tikzcd}
      \dots\ar[r,"\alpha"]
      &
      \A_2\ar[r,"\alpha"]
      \ar[dr,phantom,"\Leftarrow" very near start, "\Leftarrow" very near end]
      \ar[d,"f"{name=leftf} description]
      &
      \A_1\ar[r,"\alpha"]
      \ar[dr,phantom,"\Leftarrow" very near start, "\Leftarrow" very near end]
      \ar[d,"f"{name=middlef} description]
      \ar[dl,"h" description]
      &
      \A_0\ar[r,"\alpha"]
      \ar[d,"f"{name=rightf} description]
      \ar[dl,"h" description]
      &
      \dots
      \\
      \dots\ar[r,"\beta"']
      &
      \B_2\ar[r,"\beta"']
      &
      \B_1\ar[r,"\beta"']
      &
      \B_0\ar[r,"\beta"']
      &
      \dots
    \end{tikzcd}
  \end{equation}
  but note that this picture is not complete,
  since it does not depict the trivialization
  \(\eta\circ\epsilon\simeq 0\) encoded in the square \eqref{eq:eta-epsilon=0}.

  The forgetful functor
  \begin{equation}
    \MaplhAB\to \MapAB
  \end{equation}
  sends an \lhenhancedmor{} \(F=(f,h,\epsilon,\eta)\)
  to its underlying chain map
  by forgetting \(h\), \(\epsilon\) and \(\eta\)
  and only remembering the maps \(f\) and the equivalences
  \(f\alpha\simeq \beta f\);
  in the picture \eqref{eq:lhenhmor-as-pasting}
  this just amounts to pasting the triangular \(2\)-cells to form
  (commutative) squares.
  For each chain map \(f\colon \A_\bullet\to\B_\bullet\),
  we write \(\Maplh[f](A_\bullet,\B_\bullet)\)
  for the fiber of this forgetful functor over the object
  \(f:\Map(\A_\bullet,\B_\bullet)\);
  it is the (typically not stable) \(\infty\)-category
  of \lhenhancement{}s of the chain map \(f\).
\end{rem}

\begin{rem}
  An \lhenhancedmor{} is called exact
  if each square \eqref{eq:eta-epsilon=0} is biCartesian.
  We denote by
  \begin{equation}
    \MaplhexAB\coloneqq
    \MapoplaxgeqABex[0]_{\reskgeq[2]=0}
    \subset\MaplhAB
  \end{equation}
  the full subcategory of \exactlhenhancedmor{}s.
\end{rem}

\begin{rem}
  \label{rem:exact-lh-enhancements=null-homotopies}
  Note that under the identification of \Cref{lem:inclusion-of-lax},
  an \exactlhenhancement{} of a chain map is precisely
  a lax null-homotopy in the sense of
  \Cref{defi:lax-null-homotopy}.
\end{rem}

The following construction of the oplax mapping cone
is a tautological reformulation of what
the data of an \lhenhancedmor{} entails.

\begin{con}[Oplax mapping cone]
  \label{con:oplax-mapping-cone}
  Let \(F=(f,h,\epsilon,\eta)\colon (\A_\bullet,\alpha)\lhenhmor(\B_\bullet,\beta)\)
  be an \lhenhancedmor{} of chain complexes.
  We define the \emph{oplax mapping cone} of \(F\) to be the chain complex
  \begin{equation}
    \oplaxmapcone{F}\colon
    \dots\rightarrow
    \A_{n}\laxplus[h]\B_{n+1}
    \xrightarrow{\delta_{n+1}}
    \A_{n-1}\laxplus[h]\B_{n}
    \xrightarrow{\delta_{n}}
    \A_{n-2}\laxplus[h]\B_{n-1}
    \rightarrow\dots
  \end{equation}
  where the differential is the lax-oplax matrix
  \begin{equation}
    \delta_{n+1} \coloneqq\Donematrixmixed{\alpha_{n}}{0}{f_n}{\beta_{n+1}}
    \colon \A_n\oplaxplus[h]\B_{n+1} \to \A_{n-1}\laxplus[h]\B_{n}
  \end{equation}
  induced by the commutative square \eqref{eq:eta-epsilon=0}.
  Using the matrix multiplication formula from
  \Cref{lem:mixed-martix-multiplication}
  we compute the squared differential
  \begin{equation}
    \delta\circ\delta
    \simeq
    \Donematrixmixed{\cof(\alpha\alpha\to 0)}{\cof(0\to 0)}
    {\cof(f\alpha\to \beta h \alpha\to \beta f)}{\cof(0\to\beta\beta)}
  \end{equation}
  It is zero because \(\alpha^2\simeq 0\), \(\beta^2\simeq 0\)
  and the fact that the composite map \eqref{eq:epsilon-eta=1}
  is an equivalence.
\end{con}

Having constructed the mapping cone \(\oplaxmapcone{F}\)
with respect to the choice of the auxiliary \lhenhancement{}
of the underlying chain map \(f\),
it is natural to ask whether there are universal ways to produce such
\lhenhancement{}s.
These exists as long as the differentials \(\alpha\) and/or \(\beta\)
admit left adjoints:

\newcommand\inlhenh[1][(-)]{#1_\beta}
\newcommand\terlhenh[1][(-)]{#1_\alpha}
\begin{lem}
  \label{lem:universal-lh-enhancements}
  Let \((A_\bullet,\alpha)\) and \((B_\bullet,\beta)\)
  be two chain complexes in \(\extwoaddcat\).
  Consider the forgetful functor
  \begin{equation}
    p\colon \Maplh(\A_\bullet,\B_\bullet)\to \Map(\A_\bullet,\B_\bullet)
  \end{equation}
  \begin{enumerate}
  \item
    \label{it:beta-left-adjoint}
    If each differential \(\beta\) has a left adjoint,
    then \(p\) admits a fully faithful left adjoint \(\inlhenh\).
  \item
    \label{it:alpha-right-adjoint}
    If each differential \(\alpha\) has a left adjoint,
    then \(p\) admits a fully faithful right adjoint \(\terlhenh\).
  \item
    \label{it:zero-lh-enhancement-BC}
    Assume that both differentials \(\alpha\) and \(\beta\) admit left adjoints.
    The canonical transformation
    \(\inlhenh\to \terlhenh\)
    is an equivalence precisely on those chain maps
    \(f :\Map(\A_\bullet,\B_\bullet)\)
    which are \ldiffBC{}.
  \end{enumerate}
\end{lem}

\begin{proof}
  The first two statements are a direct consequence of
  \Cref{lem:resk-has-adjoints} (for \(\loffset=0\))
  by observing that both adjoints
  (if they exists)
  take values in
  \begin{equation}
    \MaplhAB\subset \MapoplaxgeqAB[0]
  \end{equation}
  when restricted to
  \begin{equation}
    \MapAB\subset\MapoplaxAB[0].
  \end{equation}

  To prove \ref{it:zero-lh-enhancement-BC} fix a chain map \(f:\Map(\A_\bullet,\B_\bullet)\)
  and consider the component \(\inlhenh[f]\to\terlhenh[f]\).
  The only place where it can possibly not be an equivalence is
  on the first off-diagonal.
  Unraveling the pointwise formula,
  one observes that the value at these off-diagonal places
  is given by the mates
  \(\la{\beta}f\to f\la{\alpha}\)
  of the equivalences \(f\alpha\to\beta f\);
  by definition \(f\) is \ldiffBC{} precisely
  if these mates are all equivalences.
\end{proof}

\begin{rem}
  Assume that all differentials \(\alpha\) and \(\beta\)
  admit left adjoints.
  Then the recollement \eqref{eq:recollement-resk}
  (for \(\loffset=0\))
  restricts to a recollement
  \begin{equation}
    \cdrecollement
    {\Mapoplax[1]{\A}{\B}}
    {\Maplh(\A,\B)}
    {\MapAB}
    {}{p}
  \end{equation}
  and therefore to an equivalence 
  \begin{equation}
    \begin{tikzcd}
      \Maplh(\A,\B)\ar[r,dashed,leftrightarrow,"\simeq"]
      \ar[d]
      &\laxlima{\MapAB}{\rho}{\Mapoplax[1]{\A}{\B}}
      \ar[d,"p_0"]
      \\
      \MapAB
      \ar[r,equal]
      &
      \MapAB
    \end{tikzcd}
  \end{equation}
  Pointwise over each \(f:\MapAB\)
  we thus have an equivalence
  \begin{equation}
    \Maplh[f](\A,\B)\simeq \rho(f)/\Mapoplax[1]{\A}{\B}.
  \end{equation}
  A glance at the explicit formula \eqref{eq:formula-gluing-functor-resk}
  shows that the gluing funtor \(\rho\) is zero precisely
  on those chain maps which are \ldiffBC{};
  in this case the \(\infty\)-category of \lhenhancement{}s
  is simply equivalent to \(\MapoplaxAB[1]\), hence in particular stable.
\end{rem}

The following corollary summarizes the situation over each chain map \(f\):
When the differentials of the chain complexes admit adjoints,
each chain map \(f\) can be canonically enhanced in two ways
yielding an initial or terminal object
in the category \(\Maplh[f](\A,\B)\) of \lhenhancement{}s of \(f\).
If the chain map is \ldiffBC{},
this \(\infty\)-category is stable
and these two canonical \lhenhancement{}s agree.

\begin{cor}
  \label{lem:unique-enhancement-diffadjoint}
  Let \(f\colon (\A_\bullet,\alpha)\to(\B_\bullet,\beta)\)
  be a chain map.
  \begin{enumerate}
  \item 
    \label{it:enhancement-beta}
    If each differential \(\beta\) admits a left adjoint \(\la{\beta}\)
    then \(f\) admits an initial \lhenhancement{} \(\inlhenh[f]\).
  \item
    \label{it:enhancement-alpha}
    Dually,
    if each differential \(\alpha\) admits a left adjoint \(\la{\alpha}\),
    then \(f\) admits a terminal \lhenhancement{} \(\terlhenh[f]\).
  \item
    If the chain map \(f\) is \ldiffBC{} 
    then the two \lhenhancement{}s \(\inlhenh[f]\) and \(\terlhenh[f]\) coincide.
    In this case we denote this \lhenhancement{} by \(\lhenh{f}\).
  \end{enumerate}
\end{cor}

\begin{proof}
  Follows from the adjunctions of \Cref{lem:universal-lh-enhancements}
  viewed pointwise over \(f:\Map(\A_\bullet,\B_\bullet)\).
\end{proof}

We now identify the mapping cones constructed
from the initial and terminal \lhenhancement{}
with those constructed in \cite[Construction~4.3.3]{CDW23}
using the \dpo{} and \dpb{}.

\begin{prop}
  \label{prop:mapping-cones-with-lax-squares}
  Let \(f\colon (\A_\bullet,\alpha)\to (\B_\bullet,\beta)\)
  be a chain map.
  \begin{enumerate}
  \item
    Assume that each \(\alpha\) admits a left adjoint and
    let \(\terlhenh[f]\)
    be its terminal \lhenhancement{} of \Cref{lem:unique-enhancement-diffadjoint}.
    Consider the oplax square
    \begin{equation}
      \begin{tikzcd}[column sep = huge]
        \A_{i}\ar[d,"f_i"']\ar[r,"\alpha"]
        \ar[dr,"\Downleftarrow" near start, phantom]
        &\A_{i-1}\ar[dl,"h_\alpha" description,""{name=M}]\ar[d]
        \\
        \B_{i}\ar[r]&\A_{i-1}\oplaxplus[h_\alpha]\B_{i}
        \ar[from=M,phantom,"\Downleftarrow" very near end]
      \end{tikzcd}
    \end{equation}
    obtained by pasting \(\epsilon\) with the oplax colimit cone.
    This square is a \dpo{},
    thus yields an identification
    \(
    \laxpush{\A_{i-1}}{\A_i}{\B_i}\xrightarrow{\simeq} \A_{i-1}\oplaxplus[h]\B_{i}.
    \)
    Under this identification
    the differential of \(\oplaxmapcone{\terlhenh[f]}\)
    corepresents the map
    \begin{equation}
      \label{eq:represented-by-diff-alpha}
      (\comap{a}_{i-1}\alpha\to \comap{b}_i f)
      \mapsto
      (\cof(\comap{a}_{i-1}\alpha\to \comap{b}_i f)\alpha\xrightarrow{\simeq} \comap{b}_i\beta f).
    \end{equation}
  \item
    Dually,
    if each \(\beta\) admits a left adjoint,
    then the terms of the cone \(\oplaxmapcone{\inlhenh[f]}\)
    are canonically identified with \(\laxpull{\A_{i-1}}{\B_{i-1}}{\B_i}\)
    and the differential represents the map
    \begin{equation}
      \label{eq:represented-by-diff-beta}
      (fa_i\to \beta b_{i+1})
      \mapsto
      (f\alpha a_i\xrightarrow{\simeq} \beta\fib(f a_i\to \beta b_{i+1}))[1]
    \end{equation}
  \end{enumerate}
\end{prop}

\begin{proof}
  Let \(\terlhenh[f]=(f,h_\alpha,\epsilon_\alpha,\eta_\alpha)\)
  be the terminal \lhenhancement{} of \(f\).
  We have to show that for each test object \(\C:\extwoaddcat\),
  the functor
  \begin{equation}
    \label{eq:inprf-on-representables}
    \oplaxlim_{\Delta^1}(\extwoaddcat(\A_{i-1},\C)
    \xleftarrow{\circ h_\alpha}\extwoaddcat(\B_{i},\C))
    =
    \extwoaddcat(\A_{i-1}\oplaxplus[h_\alpha]\B_i,\C)
    \to
    \laxpull{\extwoaddcat(\A_{i-1},\C)}
    {\extwoaddcat(\A_{i},\C)}
    {\extwoaddcat(\B_i,\C)}
  \end{equation}
  is an equivalence of (stable) \(\infty\)-categories.
  Explicitly, this functor sends a section 
  \(a_{i-1}'\xrightarrow{u}{\comap{b}_{i}h_\alpha}\)
  to the composite
  \begin{equation}
    a_{i-1}'\alpha\xrightarrow{u\alpha}
    {\comap{b}_{i}h_\alpha\alpha}
    =
    {\comap{b}_{i}f_i\la{\alpha}\alpha}
    \xrightarrow{\comap{b}_if_i\counit_\alpha}{\comap{b}_{i}f_{i}},
  \end{equation} 
  (where \(\counit_\alpha\colon \la{\alpha}\alpha\to \id\) is the counit of the adjunction
  \(\la{\alpha}\dashv\alpha\)),
  which is precisely its transpose under the adjunction
  \((\circ\alpha) \dashv (\circ\la{\alpha})\).
  Thus the functor
  \eqref{eq:inprf-on-representables}
  is an equivalence by
  \Cref{lem:lax-pullback-adjoint}\ref{it:lem:lax-pullback-right-adjoint}
  applied to
  \begin{equation}
    \extwoaddcat(\A_{i-1},\C)
    \xrightarrow{\circ\alpha}
    \extwoaddcat(\A_{i},\C)
    \xleftarrow{\circ f_i}
    \extwoaddcat(\B_{i},\C).
  \end{equation}
  To compute the map corepresented by the differential,
  we compute for each
  \begin{equation}
    (\comap{a}_{i-1}\xrightarrow{u} \comap{b}_{i}) : \extwoaddcat(\A_{i-1}\oplaxplus[{h_\alpha}]\B_{i},\C)
  \end{equation}
  the matrix product
  \begin{equation}
    \label{eq:calculation-corepresent-differential}
    (\comap{a}_{i-1}\xrightarrow{u} \comap{b}_{i})\circ
    \Donematrixmixed{\alpha}{0}{f_{i}}{\beta}
    =
    (\cof(\comap{a}_{i-1}\alpha\to \comap{b}_{i}f_{i}) \to \cof(0\to \comap{b}_i\beta)),
  \end{equation}
  where in the first entry we are taking the cofiber
  of the map
  \(\overline{u}\colon \comap{a}_{i-1}\alpha \xrightarrow{u\alpha} \comap{b}_{i}h_\alpha\alpha
  \xrightarrow{\comap{b}\epsilon_\alpha} \comap{b}_{i}f_{i} \).
  Note that in the matrix representation 
  \eqref{eq:calculation-corepresent-differential}
  we are omitting the application of the gluing functor as is customary.
  If we put this implicit application back in,
  we obtain the map
  \begin{equation}
    \cof(\comap{a}_{i-1}\alpha \to \comap{b}_{i}f_{i}) \to
    \comap{b}_{i}\beta h_\alpha = \comap{b}_{i}\beta f_{i+1}\la{\alpha} :
    \extwoaddcat(\A_{i},\C)
  \end{equation}
  which yields the desired map 
  \begin{equation}
    \cof(\comap{a}_{i-1}\alpha \to \comap{b}_{i}f_{i})\alpha \to \comap{b}_{i}\beta f_{i+1} :
    \extwoaddcat(\A_{i+1},\C)
  \end{equation}
  after transposing;
  it is just the equivalence
  \(\comap{b}_{i}f_{i}\alpha \simeq \comap{b}_{i}\beta f_{i+1}\)
  because \(\comap{a}_{i-1}\alpha\alpha=0\).

  The proof of the dual statement is analogous:
  We apply \(\extwoaddcat(\C,-)\) to reduce to the case of \(\infty\)-categories,
  where we apply \Cref{lem:lax-pullback-adjoint}\ref{it:lem:lax-pullback-left-adjoint}.
  Then we only have to perform the dual matrix computation
  \begin{equation}
    \Donematrixmixed{\alpha}{0}{f_i}{\beta}
    \circ
    \Donevector{a_i}{b_{i+1}}
    =
    \Donevector{\cof(\alpha a_i\to 0)}{\cof(f_ia_i\to \beta b_{i+1})}
    =\Donevector{\alpha a_i}{\fib(f_ia_i\to \beta b_{i+1})}[1]
  \end{equation}
  to obtain the desired formula.
\end{proof}

\begin{defi}
  The external shift of a chain complex \((\A_\bullet,\alpha)\)
  is defined as
  \begin{equation}
    \A[n]_\bullet\coloneqq (\A_{\bullet-n},\alpha[n]),
  \end{equation}
  where the terms are reindexed and the differentials
  are shifted internally in the stable \(\infty\)-categories
  \(\extwoaddcat(\A_i,\A_{i-1})\).
\end{defi}

\begin{con}
  Let \(f\colon (\A_\bullet,\alpha)\to (\B_\bullet,\beta)\) be a chain map.
  We define
  \begin{equation}
    \Cof(f)\coloneqq
    \oplaxmapcone{\terlhenh[f]}
    \quad\text{and}\quad
    \Fib(f)\coloneqq
    \oplaxmapcone{\inlhenh[f]}[-1],
  \end{equation}
  whenever these are defined,
  i.e., whenever \(\alpha\) or \(\beta\) has a left adjoint, respectively.
\end{con}

\begin{rem}
  \Cref{prop:mapping-cones-with-lax-squares} essentially states that
  this definition of \(\Fib(f)\) and \(\Cof(f)\)
  agrees with the one from \cite[Construction~4.3.3]{CDW23}
  in the case \(\exaddtwocat=\StL_k\).
\end{rem}

\begin{cor}
  \label{cor:cofib=fib[1]-ldiff}
  Let \(f\colon (\A_\bullet,\alpha)\to (\B_\bullet,\beta)\)
  be a \ldiffBC{} chain map.
  We have an equivalence
  \begin{equation}
    \Cof(f)\simeq \Fib(f)[1].
  \end{equation}
\end{cor}

\begin{proof}
  Since we assume that the chain map \(f\) is \ldiffBC{},
  \Cref{lem:unique-enhancement-diffadjoint} states
  that \(\terlhenh[f]\) and \(\inlhenh[f]\) are canonically equivalent
  as \lhenhancement{}s of the chain map \(f\).
  Therefore the chain complexes
  \(\Cof(f)=\oplaxmapcone{\terlhenh[f]}\) and \(\Fib(f)[1]=\oplaxmapcone{\inlhenh[f]}\)
  are also equivalent.
\end{proof}

So far we have used that one can express the
\dpb{}
\(\laxpull{\A_{i-1}}{\B_{i-1}}{\B_i}\)
and the \dpo{}
\(\laxpush{\A_{i-1}}{\A_i}{\B_i}\)
as a lax limit/colimit of a composite involving horizontal left adjoints.
To prove \cite[Proposition~4.3.12]{CDW23} we need an analogous discussion
using vertical right adjoints.
This change corresponts to changing the direction
of the gluing map between \(\B_{n}\) and \(\A_{n-1}\).
We start by defining the corresponding notion of enhancement.

\begin{defi}
  An \emph{\rvenhancedmor{}} 
  \(F\colon (\A_\bullet,\alpha_\bullet)\rvenhmor(\B_\bullet,\beta_\bullet)\)
  of chain complexes consists of 1-morphisms
  \begin{equation}
    f_n\colon\A_n\to\B_n
    \quad\text{and}\quad
    k_{n}\colon \B_{n}\to\A_{n-1}
  \end{equation}
  toghether with an oplax-lax matrix of the form
  \begin{equation}
    \label{eq:rv-matrix}
    \delta=\Donematrixmixedop{\beta}{f}{0}{\alpha}
    \colon \B_{n+1}\laxplus[k]\A_n\to\B_{n}\oplaxplus[k]\A_{n-1}
  \end{equation}
  such that the composite map
  \(f\alpha\to fkf\to \beta f\)
  is an equivalence (yielding the underlying chain map \(f\) of \(F\)).
  The resulting chain complex
  \((\oplaxmapcone[\bullet]{F}\coloneqq \B_\bullet\laxplus[k]\A_{\bullet-1},\delta)\)
  is called the oplax mapping cone of \(F\).
\end{defi}

The mnemonic ``rv'' stands for ``right-vertical''
and reflects the fact that there are canonical \rvenhancement{}s
in the presence of right adjoints in the vertical (=chain map) direction.

We shall now explain how such \rvenhancedmor{}s
assemble into an \(\infty\)-category.
For simplicity we will restrict to those,
where each \(f_n\) admits a right adjoint \(g_n\coloneqq\ra{f_n}\).

\begin{con}
  Define
  \begin{equation}
    \MaplaxleqBA[\loffset]\coloneqq \Fun_{\ZZ^\op\times\ZZ^\op}
    \left(
      \Zopleq[\loffset], \mixedGroth{n:\ZZ}{m:\ZZ^\op}\extwoaddcat(\B_n,\A_m)
    \right)
  \end{equation}
  to consist of sections defined on
  \begin{equation}
    \Zopleq[\loffset]\coloneqq \setP{(m,n)}{m\leq n+\loffset}
    \subset \ZZ^\op\times \ZZ^\op.
  \end{equation}
  Pictorially, such sections look as follows:
  \begin{equation}
    \begin{tikzcd}
      &\cdots\ar[r,"\beta"]
      &\B_2\ar[r,"\beta"]&\B_1\ar[r,"\beta"]&\B_0\ar[r,"\beta"]
      &\cdots\\
      \vdots\ar[d,"\alpha"]&\ddots\ar[dr]\\
      \A_{ \loffset+2 }\ar[d,"\alpha"]
      && g^{\loffset}_{2}\ar[d]\ar[dr]
      \\
      \A_{\loffset+1 }\ar[d,"\alpha"]
      &&g^{\loffset-1}_2\ar[d]\ar[r]& g^{\loffset}_{1}\ar[d]\ar[dr]\\
      \A_{\loffset }\ar[d,"\alpha"]
      &&g^{\loffset-2}_2\ar[r]&g^{\loffset-1}_1\ar[r]& g^{\loffset}_{0}\ar[dr]\\
      \vdots
      &&&&&\ddots
    \end{tikzcd}
  \end{equation}
  Denote by
  \begin{equation}
    \MaprvLAB\subset \MaplaxleqBA[0]
  \end{equation}
  the full subcategory of those sections \((g^r_n)\) satisfying:
  \begin{itemize}
  \item
    The lax chain map \(g_\bullet=g^0_\bullet\colon \B_\bullet\to\A_\bullet\)
    on the main diagonal is left adjointable,
    i.e.\
    each \(g_n\colon \B_n\to \A_n\) has a left adjoint \(\la{g_n}\)
    and the canonical mate
    \(\la{g_{n-1}}\alpha\to\beta\la{g_{n}}\) is an equivalence.
  \item
    The section is zero beyond the first off-diagonal,
    i.e. \(g^r_\bullet=0\) for \(r\leq -2\).
  \end{itemize}

  Using the dual of \Cref{rem:data-of-section-f-h-0}
  and by passing from an adjointable lax chain map
  \(g_\bullet=g^0_\bullet\colon \B_\bullet\to\A_\bullet\)
  to its adjoint
  \(f_\bullet\coloneqq \la{g_\bullet}\colon\A_\bullet\to \B_\bullet\)
  (which is an honest chain map),
  it is not hard to see that the data of such a section
  amounts precisely to that of an \rvenhancedmor{}
  whose underlying chain map \(f\) admits pointwise adjoints;
  the \(1\)-morphisms \(k_{n}\colon \B_{n}\to\A_{n-1}\)
  are the term \(g^{-1}_n\) on the first off-diagonal
  and the matrices \eqref{eq:rv-matrix}
  amount precisely to the squares
  \begin{equation}
    \cdsquareNA{k_{n+1}}{g_n}{0}{k_n}.
  \end{equation}

  Therefore we can view \(\MaprvLAB\) as the \(\infty\)-category of
  those \rvenhancedmor{}s,
  whose underlying chain map \(f\) admits pointwise adjoints.
  We have the canonical forgetful functor
  \begin{equation}
    \label{eq:forgetful-functor-rv}
    \begin{tikzcd}
      \MaplaxleqBA[0]_{\reskleq[-2]=0}\ar[r,"{\resk[0]}"]
      &\Maplax[0]{\B}{\A}\\
      \MaprvLAB \ar[u,hookrightarrow]\ar[r,"{\resk[0]}"]\ar[dr,dashed]
      & \set{\text{left adjointable }g_\bullet}\ar[u,hookrightarrow]
      \ar[d,leftrightarrow,"{\simeq}"']
      \\
      &\set{\text{pointwise right adjointable }f_\bullet}\ar[r,hookrightarrow]
      &\Map(\A,\B)
    \end{tikzcd}
  \end{equation}
  sending such an \rvenhancedmor{} to its underlying chain map.
  For each pointwise adjointable chain map
  \(f\colon \A\to\B\) we write
  \(\MaprvAB[f]=\MaprvLAB[f]\) for the fiber of
  this dashed functor;
  it is the \(\infty\)-category of \rvenhancement{}s of \(f\).
\end{con}

\begin{lem}
  The forgetful functor \eqref{eq:forgetful-functor-rv}
  is part of a recollement
  \begin{equation}
    \cdrecollementcompl
    {\Maplax[-1]{\B}{\A}}
    {\MaprvLAB}
    {\set{\substack{\text{pointwise right adjointable}\\ f_\bullet\colon \A\to \B}}}
    {}{}{j}{}{j'}{},
  \end{equation}
  whose gluing functor \(\rho\) computes the fiber of the canonical mate,
  i.e.,
  \begin{equation}
    \rho(f)=\left(\fib(\alpha \ra{f_n}\to \ra{f_{n-1}}\beta)[n]\right)_n.
  \end{equation}
\end{lem}

\begin{proof}
  Similarly to \Cref{lem:reskgeq-has-adjoints}
  and \Cref{cor:reskgeq-recollement} we compute that the relative
  left and right Kan extension along the diagonal
  \(\ZZ^\op\hookrightarrow\Zopleq\) always exist,
  yielding fully faithful left and right adjoints \(j\) and \(j'\)
  to the restriction functors \(\resk[0]\).
  Explicitly we have
  \begin{equation}
    j(g)^{-1}_n=\alpha g_n \to g_{n-1} \beta = j'(g)^{-1}_n
  \end{equation}
  (the structure map of \(g_\bullet\)) on the first off-diagonal
  and zero beyond it.
  Moreover, the kernel of the forgetful functor is the full subcategory
  \begin{equation}
    \setP{g}{\forall r\neq -1 : g^{r}_\bullet=0}
    \subset\MaplaxleqBA,
  \end{equation}
  which, similarly to \Cref{lem:inclusion-of-oplax}
  we can identify with \(\Maplax[-1]{\B}{\A}\)
  via the assignment
  \begin{equation}
    g^\bullet_\bullet\mapsto (g^{-1}_n[n])_n.
  \end{equation}
  The desired result follows by passing to adjoints, i.e.,
  \(g_\bullet\coloneqq \ra{f_\bullet}\).
\end{proof}

\begin{rem}
  The gluing functor for the recollement
  \begin{equation}
    \cdrecollement{\Maplax[-1]{\B}{\A}}
    {\MaplaxleqBA[0]_{\reskleq[-2]=0}}
    {\Maplax[0]{\B}{\A}}
    {}
    {\resk[0]}
  \end{equation}
  (before restricting the cokernel to the subcategory of
  left adjointable maps \(g\colon \B\to\A\))
  is nothing but the differential of the lax mapping complex
  \(\Maplax{\B}{\A}\).
\end{rem}

As a direct consequence we get the following result,
which provides the two canonical \rvenhancement{}s
of a pointwise right adjointable chain map.

\begin{cor}
  Let \(f\colon (\A_\bullet,\alpha)\to(\B_\bullet,\beta)\)
  be a pointwise right adjointable chain map.
  The \(\infty\)-category \(\MaprvAB[f]\) has
  \begin{enumerate}
  \item
    an initial object \(f^\alpha\coloneqq j(f)\)
    with \(k^\alpha=\alpha \ra{f}\) and where
    the vertical map \(\alpha\to kf=\alpha \ra{f} f\)
    in the matrix \eqref{eq:rv-matrix} is the unit;
  \item
    a terminal object \(f^\beta\coloneqq j'(f)\)
    with \(k^\beta=\ra{f}\beta\) and where
    the horizontal map \(\beta \leftarrow f k = f \ra{f} \beta \)
    in the matrix \eqref{eq:rv-matrix} is the counit.
  \item
    These two \rvenhancement{}s coincide
    if and only if the chain map \(f\) is right adjointable.
    In this case we denote this \rvenhancement{} by \(\rvenh{f}\).
  \end{enumerate}
\end{cor}

Analogously to \Cref{prop:mapping-cones-with-lax-squares},
we can exhibit the terms of the corresponding oplax mapping cones
\(\oplaxmapcone{f^\alpha}\) and \(\oplaxmapcone{f^\beta}\)
as a \dpo{} or \dpb{}, respectively.

\begin{prop}
  Let \(f\colon (\A_\bullet,\alpha)\to (\B_\bullet,\beta)\)
  be a chain map and assume that each \(f_i\) admits a right adjoint.
  \begin{enumerate}
  \item
    The oplax square
    \begin{equation}
      \begin{tikzcd}[column sep = huge]
        \A_i\ar[d,"f"']\ar[r,"\alpha"]
        \ar[dr,"\Downleftarrow" near start, phantom]
        &\A_{i-1}\ar[d,""]
        \\
        \B_i\ar[ur,"k^\alpha" description,""{name=M}]\ar[r]&
        \B_i\laxplus[k^\alpha]\A_{i-1}
        \ar[from=M,phantom,"\Downleftarrow" very near end]
      \end{tikzcd}
    \end{equation}
    yields an identification
    \(\laxpush{\A_{i-1}}{\A_i}{\B_{i}}\xrightarrow{\simeq}
    \B_i\laxplus[k^\alpha]\A_{i-1} \).
    Under this identification,
    the differential of \(\oplaxmapcone{f^\alpha}\)
    again corepresents the map \eqref{eq:represented-by-diff-alpha}.
  \item
    Dually, the terms of the cone \(\oplaxmapcone{f^\beta}\)
    are canonically identified with \(\laxpull{\A_{i-1}}{\B_{i-1}}{\B_i}\)
    and the differential
    again represents the map \eqref{eq:represented-by-diff-beta}.
  \end{enumerate} 
\end{prop}

\begin{proof}
  Similar to \Cref{prop:mapping-cones-with-lax-squares}; omitted.
\end{proof}

\begin{cor}
  \label{cor:cofib=fib[1]}
  The chain complexes \(\oplaxmapcone{f^\alpha}\)
  and \(\oplaxmapcone{f^\beta}\)
  also yield a construction for \(\Cof(f)\) and \(\Fib(f)[1]\), respectively.
  In particular, \(\Cof(f)\) and \(\Fib(f)[1]\) agree
  when the chain map \(f\) is right adjointable.
\end{cor}

\begin{cor}
  When \(f\) is both right adjointable and \ldiffBC{},
  the two canonical oplax mapping cones
  \(\oplaxmapcone{\lhenh{f}}\) and \(\oplaxmapcone{\rvenh{f}}\)
  agree.
\end{cor}

\begin{rem}
  Throughout this section there was a bias in our discussion,
  since we implicitly treated chain maps
  as being \emph{oplax},
  i.e.\ having directed squares of the form
  \begin{equation}
      \oplaxsquare{\A_{i+1}}{\A_i}{\B_{i+1}}{\B_i}{\alpha}{f_{i+1}}{f_i}{\beta}{}
  \end{equation}
  This was already apparent in the chosen direction for \dpo{}s and \dpb{}s
  in \Cref{exa:directed-pull-push}
  and accounts for the two possible choices we had
  when it came to adjointability conditions:
  having vertical right adjoints or horizontal left adjoints.
  We could rewrite this whole section with the opposite conventions
  and obtain the \emph{lax mapping cone}
  \(\laxmapcone{F}\) associated to a chain map \(f\)
  with suitable enhancements.
  In the case where \(f\)
  is left adjointable or \rdiffBC{}
  we could again construct a canonical lax mapping cone
  \(\laxmapcone{F}\)
  whose terms are identified both with
  \(\laxpush{\B_i}{\A_i}{\A_{i-1}}\) and with
  \(\laxpull{\B_{i}}{\B_{i-1}}{\A_{i-1}}\).
\end{rem}

\section{Universal property of the lax mapping cone}\label{sec:univpropcone}

The main reason for introducing the mapping cone of a chain map
\(f\colon (A_\bullet,\alpha)\to (B_\bullet,\beta)\)
between chain complexes in an additive category \(\A\)
is that it yields an explicit model for the cofiber of \(f\)
in the stable \(\infty\)-category \(\htpycat{\A}\)
of chain complexes up to chain homotopy.
In other words, it satisfies
\begin{align}
  \Map_{\htpycat{\A}}(\mappingcone{f},C)
  &\simeq
    \fib\left(\Map(B,C)\to\Map(A,C)\right)
  &=
    \set{(g\colon B\to C, h\colon gf\simeq 0)}
\end{align}
naturally in \(C:\htpycat{\A}\).

Already before passing to the stable \(\infty\)-category \(\htpycat{\A}\),
one can see a naive version of this universal property
characterizing the mapping cone up to isomorphism in \(\Ch(\A)\) via
\begin{equation}
  \label{eq:mapping-cone-aive-univ-prop}
  \Ch(\A)(\mappingcone{f},C) \cong
  \setP{(g, h)}{g\colon B\to C, h\colon gf\simeq 0}
\end{equation}
naturally in \(C:\Ch(\A)\).
In other words: maps out of \(\mappingcone{f}\) are chain maps
\(g\colon B\to C\) together with a null-homotopy of \(gf\).

Ultimately, we are of course
interested in understanding the categorified analog
of the homotopically meaningful universal property.
However, this is currently out of reach since we don't even know
what the correct analog of the stable \(\infty\)-category \(\htpycat{\A}\)
should be and in what sense we are supposed to view the mapping cone
as a cofiber.
Therefore, we now instead describe the categorified analog of
\eqref{eq:mapping-cone-aive-univ-prop}
in the hopes that it might lead to a better understanding
of the theory of categorified chain complexes up to homotopy.

\begin{thm}
  \label{thm:univ-prop-of-lh-cone}
  Let
  \(F\colon \A\lhenhmor\B\)
  be an \lhenhancedmor{} of chain complexes
  with underlying chain map \(f\).
  \begin{enumerate}
  \item
    For each chain complex \(\C:\Ch(\extwoaddcat)\)
    there is a natural equivalence of (stable) \(\infty\)-categories between
    \begin{itemize}
    \item
      chain maps \(\oplaxmapcone{F}\to\C\) and
    \item
      chain maps \(g\colon \B\to\C\)
      together with an \exactlhenhancement{} \(E\) of \(gf\)
      and a morphism \(E\to gF\) of \lhenhancement{}s of \(gf\).
    \end{itemize}
  \item
    For each chain complex \(\C:\Ch(\extwoaddcat)\)
    there is a natural equivalence of (stable) \(\infty\)-categories between
    \begin{itemize}
    \item
      chain maps \(\C\to\oplaxmapcone{F}[-1]\) and
    \item
      chain maps \(g\colon \C\to\A\)
      together with an \exactlhenhancement{} \(E\) of \(fg\)
      and a morphism \(Fg\to E\) of \lhenhancement{}s of \(fg\).
    \end{itemize}
  \end{enumerate}
\end{thm}

Before proving \Cref{thm:univ-prop-of-lh-cone},
we isolate the special case
where \(F\) is the initial or terminal \lhenhancement{} of \(f\).
\begin{cor}
  \label{cor:univ-prop-lax-cof-fib}
  Let \(f\colon (\A_\bullet,\alpha)\to(\B_\bullet,\beta)\) be a chain map.
  \begin{enumerate}
  \item
    Assume that all differentials \(\alpha\) have left adjoints.
    Then for each chain complex \(\C:\Ch(\extwoaddcat)\)
    there is an equivalence of (stable) \(\infty\)-categories
    between
    \begin{itemize}
    \item
      chain maps \(\Cof(f)\to \C\) and
    \item
      chain maps \(g\colon \B\to\C\)
      together with a lax null-homotopy \(E\) of \(gf\).
    \end{itemize}
  \item
    Assume that all differentials \(\beta\) have left adjoints.
    Then for each chain complex \(\C:\Ch(\extwoaddcat)\)
    there is an equivalence of (stable) \(\infty\)-categories
    between
    \begin{itemize}
    \item
      chain maps \(\C\to\Fib(f)\) and
    \item
      chain maps \(g\colon \C\to\A\)
      together with a lax null-homotopy \(E\) of \(fg\).
    \end{itemize}
  \end{enumerate}
\end{cor}
\begin{proof}
  We prove the first statement; the second is dual.
  Let \(\terlhenh[f]\) be the terminal \lhenhancement{} of \(f\)
  and recall that we have
  \(\Cof(f)=\oplaxmapcone{\terlhenh[f]}\).
  Observe further, that composition with \(g\)
  sends the \lhenhancedmor{} \(\terlhenh[f]\)
  to \(g(\terlhenh[f])\simeq \terlhenh[(gf)]\),
  which is thus a terminal object of \(\Maplh[gf](\A_\bullet,\C_\bullet)\).
  Therefore the claim follows from \Cref{thm:univ-prop-of-lh-cone}
  after identifying \exactlhenhancement{}s with lax null-homotopies
  (see \Cref{rem:exact-lh-enhancements=null-homotopies}).
\end{proof}

\begin{proof}[Proof of \Cref{thm:univ-prop-of-lh-cone}]
Fix an \lhenhancedmor{}
\(F=(f,h,\epsilon,\eta)\colon
(\A_\bullet,\alpha_\bullet)\lhenhmor(\B_\bullet,\beta_\bullet)\)
and a test chain complex \((\C,\gamma)\) in \(\extwoaddcat\).
We unravel the data encoded in a chain map
\((\oplaxmapcone[\bullet]{F},\delta)\to(\C_\bullet,\gamma)\),
using lax-oplax matrices.
For each \(n\), we have a map
\begin{equation}
  G_n=\Donecovectorop[\mu_{n-1}]{k_{n-1}}{g_{n}}\colon\oplaxmapcone[n]{F}
  =\A_{n-1}\laxplus\B_n\to\C_n
\end{equation}
and an equivalence \(G_{n}\delta_{n+1} \xrightarrow{\simeq} \gamma_{n+1} G_{n+1}\),
which we can expand to
\begin{equation}
  \Donecovectorop{\cof{(k_{n-1}\alpha_{n}\to g_nf_{n})}}{g_n\beta_{n+1}}
  =
  \Donecovectorop{k_{n-1}}{g_n}
  \Donematrixmixed{\alpha_{n}}{0}{f_{n}}{\beta_{n+1}}
  \xrightarrow{\simeq}
  \Donecovectorop{\gamma_{n+1}k_{n}}{\gamma_{n+1}g_{n+1}}
\end{equation}

Therefore, the map
\(G_n\delta_{n+1}\to\gamma_{n+1}G_{n+1}\)
amounts to a cube (read back to front)
\begin{equation}
  \label{eq:map-from-cone-as-cube}
  \begin{tikzcd}
    {k_{n-1}\alpha_n}\ar[rrr]\ar[dr]\ar[d,"\mu_{n-1}\alpha_n"']
    \ar[dd,"\zeta_{n}"',bend right =90]
    &&&{0}\ar[dd]\ar[dr]
    \\
    g_nh_{n-1}\ar[d,"g_n\epsilon_{n}"']\alpha_n\ar[rrru]&
    {0}\ar[from=ul,crossing over]&&
    &{0}\ar[dd]
    \ar[from=lll,crossing over]
    \\
    {g_nf_n}\ar[dr,"\nu_n"']\ar[rrr,"g_n\eta_{n}" near end]&&&{g_n\beta_{n+1}}\ar[rd,"\phi_{n+1}"]
    \\
    &{\gamma_{n+1}k_n}\ar[rrr,"\gamma_{n+1}\mu_{n}"]\ar[from=uu,crossing over]
    &&&{\gamma_{n+1}g_{n+1}}
    \\
    \phantom{1}\ar[r,"{\extwoaddcat(\A_n,\C_n)}", phantom]&\phantom{1}
    &&
    \phantom{1}\ar[r,"{\extwoaddcat(\B_{n+1},\C_n)}", phantom]
    \ar[ll,"{-\circ h_n}"]
    &\phantom{1}
  \end{tikzcd}
\end{equation}
in the contravariant Grothendieck construction of \(-\circ h_n\).
The fact that this map is an equivalence amounts to saying that
the left and right squares of the cube are biCartesian.
In particular we can focus on the the right face
and see an equivalence
\(\phi_{n+1}\colon g_{n}\beta_{n+1}\xrightarrow{\simeq} \gamma_{n+1}g_{n+1}\),
exhibiting \(g_\bullet\colon (\B_\bullet,\beta)\to(\C_\bullet,\gamma)\)
as a chain map.

Consider the functor
\begin{equation}
  \Xi\colon\ZZ^\op\times \ZZ\times \Delta^1\rightarrow\st;
  \quad\quad
  (m,n,-) \mapsto
  \left(\extwoaddcat(\A_n,\B_m)\xrightarrow{g_n\circ -}\extwoaddcat(\A_n,\C_m)\right),
\end{equation}
which is well defined because \(g\colon \B_\bullet\to\C_\bullet\) is a chain map.

By direct comparison with the diagram \eqref{eq:map-from-cone-systematic} below,
one verifies that all of the data
\eqref{eq:map-from-cone-as-cube}
can then be equivalently encoded as \(\Zgeq\times\Delta^1\)-sections
of the mixed
(contravariant, contravariant, covariant) Grothendieck construction
of \(\Xi\)
that satisfy
\begin{itemize}
\item
  the restriction to \(\Zgeq\times{0}\)
  is the original \lhenhancedmor{} \(F\colon\A\lhenhmor\B\),
\item
  the value on each edge \((n,n,1)\to (n,n,0)\) is Cartesian,
\item
  the restriction to \(\Zgeq\times\set{1}\)
  is an \exactlhenhancedmor{} \(E\colon \A\lhenhmor\C\).
\end{itemize}
\begin{equation}
  \label{eq:map-from-cone-systematic}
  \begin{tikzcd}[column sep= small]
    &&\phantom{1}\ar[r,phantom,"{\A_{n+1}}"]&\phantom{1}
    \ar[rr,"\alpha"]
    &&
    \phantom{1}\ar[r,phantom,"{\A_{n}}"]&\phantom{1}
    \ar[rr,"\alpha"]
    &&
    \phantom{1}\ar[r,phantom,"{\A_{n-1}}"]&\phantom{1}
    \\
    &\C_{n+1}\ar[ddd,"\gamma"]&&gf\ar[dl,"!"']&&&
    k\ar[from=ddd,"\nu" near start]
    \ar[lll,"\zeta"']
    \ar[dl,"\mu"']
    &&&
    0\ar[lll]\ar[from=ddd]
    \ar[dl]
    \\
    \B_{n+1}\ar[ur,"g"]\ar[ddd,"\beta"]&&f&&&
    h\ar[lll,"\epsilon"']
    \ar[from=ddd,"\eta"]
    &&&
    0\ar[lll,crossing over]
    &\\
    \\
    &\C_{n}\ar[ddd,"\gamma"]&&&&&gf\ar[dl,"!"']&&&
    k\ar[from=ddd,"\nu"]
    \ar[lll,"\zeta"' near end]
    \ar[dl,"\mu"']
    \\
    \B_{n}\ar[ur,"g"]\ar[ddd,"\beta"]
    &&&&&f&&&
    h\ar[lll,"\epsilon"']
    \ar[uuu,crossing over]
    \\
    \\
    &\C_{n-1}&&&&&&&&gf\ar[dl,"!"']\\
    \B_{n-1}\ar[ur,"g"]&&&&&&&&f\ar[uuu,"\eta"]
  \end{tikzcd}
\end{equation}

In other words,
we have \exactlhenhancedmor{}s \(E\colon \A\lhenhmor\C\)
equipped with a map \(E\to gF\)
which induces an equivalence on the underlying chain maps.
This completes the proof.
\end{proof}

\appendix

\section{Some lemmas from (2-)category theory}
\label{sec:lemmas}

\newcommand\exB{\mathbb{B}}
\newcommand\exC{\mathbb{C}}
\newcommand\Graytensor{\mathbin{\otimes^{\lax}}}

\subsection{About (op)lax limits of \texorpdfstring{\(\infty\)}{infinity}-categories}
We collect here a few useful lemmas regarding various types of 2-categorical
limits of \(\infty\)-categories or stable \(\infty\)-categories.

\begin{con}
  Let \(S\) be an \(\infty\)-category and
  \(\X\colon S\to \Catinfty\) an \(S\)-indexed diagram of
  \(\infty\)-categories.
  Let \(p\colon\covGroth [S] {\X}\to S\)
  be its (covariant) Grothendieck construction.
  Assume that for every arrow \(f\colon s\to t\) in \(S\),
  the functor \(\X_f\) admits a right adjoint.
  In this case, the cocartesian fibration \(p\) is also cartesian;
  it corresponds to the diagram
  \(\ra{\X}\colon S^\op\to \Catinfty\)
  which is obtained form \(\X\) by passing to right adjoints.
  Therefore we obtain a tautological identification
  \begin{equation}
    \label{eq:laxlimX=oplaxlimX^*}
    \laxlim_S\X = \set{\text{sections of p}} = \oplaxlim_{S^\op}\ra{\X}.
  \end{equation}
\end{con}

\begin{cor}
  \label{cor:lax-oplax-adjoint}
  Let \(f\colon \A\to\B\) be a diagram of \(\infty\)-categories
  and assume that \(f\) has a right adjoint \(\ra{f}\).
  Then there is a natural identification
  \begin{equation}
    \laxlima{\A}{f}{\B}
    =
    \oplaxlima{\B}{\ra{f}}{\A}
  \end{equation}
  given by the formula
  \begin{equation}
    (a,b, fa\xrightarrow{u}b)
    \leftrightarrow
    (a,b, a\xrightarrow{\overline{u}} \ra{f} b)
  \end{equation}
\end{cor}

\begin{proof}
  This is just the special case \(S=\Delta^1\) of
  the identification \eqref{eq:laxlimX=oplaxlimX^*}.
\end{proof}

\begin{lem}
  \label{lem:lax-pullback-adjoint}
  Let \(\A\xrightarrow{f}\C\xleftarrow{g}{\B}\)
  be a diagram of \(\infty\)-categories.
  \begin{enumerate}
  \item
    \label{it:lem:lax-pullback-right-adjoint}
    Assume that \(f\) has a right adjoint \(\ra{f}\).
    Then there is a natural equivalence
    \begin{equation}
      \laxpull{\A}{\C}{\B}\simeq
      \oplaxlima{\B}{\ra{f}g}{\A}
    \end{equation}
    given by the formula
    \begin{equation}
      (a, b, fa\xrightarrow{u} gb)
      \leftrightarrow
      (a, b, a\xrightarrow{\overline{u}} \ra{f} gb).
    \end{equation}
  \item
    \label{it:lem:lax-pullback-left-adjoint}
    Assume that \(g\) has a left adjoint \(\la{g}\).
    Then there is a natural equivalence
    \begin{equation}
      \laxpull{\A}{\C}{\B}\simeq
      \laxlima{\A}{\la{g}f}{\B}
    \end{equation}
    given by the formula
    \begin{equation}
      (a, b, fa\xrightarrow{u} gb)
      \leftrightarrow
      (a, b, \la{g}fa\xrightarrow{\overline{u}} b).
    \end{equation}
  \end{enumerate}
\end{lem}

\begin{proof}
  We compute
  \begin{align}
    \laxpull{\A}{\C}{\B}
    &=
    \A\times_{\C^{\set{0}}}\C^{\set{0\to 1}}\times_{\C^{\set{1}}}\B
    \simeq
      \laxlima{\A}{f}{\C} \times_{\C} \B
    \\
    &\simeq
    \oplaxlima{\C}{\ra{f}}{\A} \times_{\C} \B
    \simeq
    \A^{\set{0\to 1}}\times_{\A^\set{1}} \C \times_{\C} \B
    \\
    &\simeq
    \A^{\set{0\to 1}}\times_{\A^\set{1}} \B
    \simeq
    \oplaxlima{\B}{\ra{f}g}{\A}
  \end{align}
  where we have used \Cref{cor:lax-oplax-adjoint} in the third step
  and the explicit construction of the lax/oplax limit
  in steps two, four and six.
  Chasing through the chain of identifications one immediately obtains
  the desired formula.

  The second statement is analogous,
  this time using the description
  \begin{equation}
    \laxpull{\A}{\C}{\B}\simeq \A\times_{\C}\oplaxlima{\B}{f}{\C}
  \end{equation}
  and applying \Cref{cor:lax-oplax-adjoint} in the other direction.
\end{proof}

\begin{cor}
  \label{cor:summary-lemmas-lax-directed}
  Let \(\exC\) be an \((\infty,2)\)-category.
  \begin{enumerate}
  \item
    For each arrow \(f\colon \A\to\B\) in \(\exC\)
    with a right adjoint \(\ra{f}\), we have natural equivalences
    \begin{equation}
      \laxlima{\A}{f}{\B}
      \simeq
      \oplaxlima{\B}{\ra{f}}{\A}
      \quad\text{and}\quad
      \laxcolima{\B}{\ra{f}}{\A}
      \simeq
      \oplaxcolima{\A}{f}{\B}.
    \end{equation}
  \item
    For each diagram
    \(\A\xrightarrow{f}\C\xleftarrow{g}\B\) in \(\exC\),
    we have a natural equivalences
    \begin{equation}
      \laxlima{\A}{\la{g}f}{\B}
      \simeq
      \laxpull{\A}{\C}{\B}
      \simeq
      \oplaxlima{\B}{\ra{f}g}{\A}.
    \end{equation}
    assuming that \(g\) has a left adjoint or \(f\) has a right adjoint.
  \item
    For each diagram
    \(\A\xleftarrow{f}\C\xrightarrow{g}\B\) in \(\exC\),
    we have a natural equivalences
    \begin{equation}
      \laxcolima{\B}{f\ra{g}}{\A}
      \simeq
      \laxpush{\A}{\C}{\B}
      \simeq
      \oplaxcolima{\A}{g\la{f}}{\B}
    \end{equation}
    assuming that \(g\) has a right adjoint or \(f\) has a left adjoint.
  \end{enumerate}
  Each of these equivalences represents (in the case ``\({\amalg}\)'')
  or corepresents (in the case ``\({\times}\)'')
  the corresponding equivalences of \Cref{cor:lax-oplax-adjoint}
  and \Cref{lem:lax-pullback-adjoint}.
\end{cor}

\begin{proof}
  All the relevant objects are characterized either
  by their represented or corepresented functor,
  hence we may reduce to the case of lax limits
  and directed pullbacks in \(\CCatinfty\).
  This case is established in
  \Cref{cor:lax-oplax-adjoint} and \Cref{lem:lax-pullback-adjoint}.
\end{proof}

\subsection{About adjoints in diagram 2-categories}
\label{sec:adjoints-in-diagrams}

Let \(\exB,\exC\) be two \((\infty,2)\)-categories.
\begin{itemize}
\item
  By \(\FUNlax(\exB,\exC)\) and \(\FUNoplax(\exB,\exC)\)
  we denote the \((\infty,2)\)-category
  of functors \(\exB\to\exC\)
  and \emph{lax}/\emph{oplax} natural transformations
  \(\eta\colon F\to G\) between them,
  which assigns to each morphism \(f\colon B\to B'\)
  in \(\exB\) a square
  \begin{equation}
    \label{eq:lax-nat-transf}
    \laxsquare{FB}{FB'}{GB}{GB'}{Ff}{\eta^B}{\eta^{B'}}{Gf}{}
    \quad\text{or}\quad
    \oplaxsquare{FB}{FB'}{GB}{GB'}{Ff}{\eta^B}{\eta^{B'}}{Gf}{}
  \end{equation}
  respectively.
  Formally, the functors \(\FUNlax(\exB,-)\) and \(\FUNoplax(\exB,-)\)
  can be defined as right adjoints
  to the lax and oplax Gray tensor products;
  for example, see \cite[Section~3]{haugseng21} and references therein.
\item
  By \(\FUN(\exB,\exC)\) we denote the standard internal hom
  in the \((\infty,2)\)-category of \((\infty,2)\)-categories;
  it can be identified with the wide, locally full subcategory
  of \(\FUNlax(\exB,\exC)\) and \(\FUNoplax(\exB,\exC)\)
  containing only those 1-morphisms \(\eta\),
  where the squares \eqref{eq:lax-nat-transf}
  contain invertible 2-cells.
\end{itemize}

If each component \(\eta^B\)
of a lax natural transformation \(\eta\colon F\to G\) has a left adjoint
\(\la{\eta^B}\),
then these assemble to an oplax natural transformation
\(\la{\eta}\colon F\to G\)
whose oplax naturality squares
\begin{equation}
  \oplaxsquare{GB}{GB'}{FB}{FB'}{Gf}{\la {\eta^B}}{\la {\eta^{B'}}}{Ff}{}
\end{equation}
are the canonical mates of the squares \eqref{eq:lax-nat-transf}.
Dually, each oplax transformation \(\eta\) has a canonical mate
\(\ra{\eta}\) (which is a lax transformation),
whenever its components have right adjoints.
Finally note that each natural transformation \(\eta\) can be viewed
both as a lax and as an oplax transformation,
thus has both mates \(\la{\eta}\) (oplax) and \(\ra{\eta}\) (lax), provided
that all the required componentwise adjoints exist.

The following result due to Haugseng characterizes the morphisms in
\(\FUN(\exB,\exC)\) which have a adjoints.

\begin{prop}[\cite{haugseng21}, Theorem~4.6]
  Let \(\eta\colon F\to G\colon \exB\to\exC\) be a natural transformation.
  \begin{enumerate}
  \item
    As a morphism in \(\FUNlax(\exB,\exC)\),
    the transformation \(\eta\) has a right adjoint if and only if
    each component \(\eta^B\) has a right adjoint in \(\exC\).
    The right adjoint \(\ra{\eta}\) is its canonical mate,
    where \(\eta\) is viewed as an oplax transformation.
  \item
    As a morphism in \(\FUNoplax(\exB,\exC)\),
    the transformation \(\eta\) has a left adjoint if and only if
    each component \(\eta^B\) has a left adjoint in \(\exC\).
    The left adjoint \(\la{\eta}\) is its canonical mate,
    where \(\eta\) is viewed as a lax transformation.
  \end{enumerate}
\end{prop}

This result also explains our terminology from
\Cref{defi:adjointable-square}.

\begin{cor}
  \label{cor:adjoint-criterion-nat-traf}
  Let \(\eta\colon F\to G\colon \exB\to\exC\) be a natural transformation.
  As a morphism in \(\FUN(\exB,\exC)\) it has
  \begin{enumerate}
  \item
    a right adjoint if and only if the naturality square
    \eqref{eq:lax-nat-transf}
    is vertically right adjointable,
  \item
    a left adjoint if and only if the naturality square
    \eqref{eq:lax-nat-transf}
    is vertically left adjointable,
  \end{enumerate}
  In each case, the left/right adjoint is the corresponding canonical mate.
\end{cor}

\begin{proof}
  Beyond the existence of adjoints,
  the right/left vertical adjointability condition
  states precisely that the 2-cells in the canonical mates
  \(\ra{\eta}\), \(\la{\eta}\) are again invertible,
  thus providing a right/left adjoint in \(\FUN(\exB,\exC)\)
  and not just in \(\FUNlax(\exB,\exC)\)/\(\FUNoplax(\exB,\exC)\).
\end{proof}

\begin{rem}
  \label{rem:lax-(co)limit-adjoints-are-pointwise}
  Let \(S\) be an \(\infty\)-category and
  \(\alpha\colon\X\to\Y\colon S\to\exC\)
  a natural transformation of \(S\)-diagrams in \(\exC\).
  Assume that each component
  \(\alpha_s\) has a left/right adjoint \(\beta_s\)
  and that all naturality squares of \(\alpha\)
  are vertically left/right adjointable.
  \Cref{cor:adjoint-criterion-nat-traf} tells us that in this case
  the components \(\beta_s\) assemble to a natural transformation
  \(\beta\colon \Y\to\X\) which is in the diagram category
  \(\FUN(S,\exC)\) a left/right adjoint to \(\alpha\).
  Assuming that \(\exC\) has lax limits or colimits of shape \(S\),
  we can apply the \(2\)-functors
  \begin{equation}
    \laxcolim\colon\FUN(S,\exC) \to \exC
    \quad\text{and}\quad
    \laxlim\colon\FUN(S,\exC) \to \exC
  \end{equation}
  to get corresponding adjunctions
  \begin{equation}
    \laxcolim_s\alpha_s\colon \laxcolim\X\leftrightarrow
    \laxcolim\Y \noloc \laxcolim_s\beta_s
  \end{equation}
  and
  \begin{equation}
    \laxlim_s\alpha_s\colon \laxlim\X\leftrightarrow
    \laxlim\Y\noloc\laxlim_s\beta_s
  \end{equation}
  which in the lax semiadditive case are identified with each other.
\end{rem}

\section{Recollements of stable $\infty$-categories}

\label{sec:recollements}
For the convenience of the reader,
we quickly summarize the basic theory of recollements of stable \(\infty\)-categories, or equivalently that of semiorthogonal decompositions with a gluing functor, as it is used in \Cref{subsec:chain-maps} and \Cref{sec:lax-mapping-cone}
without further mention.
For comprehensive treatments, see for instance \cite[Appendix~A.8]{lurie:ha} and \cite[Section 2]{DKSS:spherical}.

All \(\infty\)-categories in this section are stable, all functors exact.

\begin{defi}
  \label{defi:recollement}
  A recollement of stable \(\infty\)-categories is a diagram
  \begin{equation}
    \label{eq:recollement}
    \cdrecollementcompl
    \A\C\B
    i p
    j q
    {j'} {q'}
  \end{equation}
  of adjunctions \(q\dashv i\dashv q'\) and \(j\dashv p \dashv j'\),
  such that
  \begin{itemize}
  \item
    the functor \(i\) is fully faithful and exhibits \(\A\) as the kernel of \(p\);
  \item
    the functors \(j\) and \(j'\) are fully faithful and exhibit
    \(\B\) as the kernel of \(q\) and of \(q'\), respectively.
  \end{itemize}
\end{defi}

\newcommand\lorth[1]{{^{\bot}{#1}}}
\newcommand\rorth[1]{{{#1}^{\bot}}}

In such a recollement, one often views \(\A\) as a full subcategory of \(\C\)
via \(i\);
then \(\B\) is viewed either as its left or its right orthogonal complement
\begin{equation}
  \lorth{\A}\coloneqq \{c : \C \mid \C(c,i(-))=0\}=\ker(q)\simeq \B
  \quad
  \text{and}
  \quad
  \rorth{\A}\coloneqq \{c : \C \mid \C(i(-),c)=0\}=\ker(q')\simeq \B
\end{equation}
depending on whether we view \(\B\) as embedded in \(\C\)
via \(j\) or via \(j'\).
Note that while these two complements are both identified with \(\B\),
they are \emph{not} the same subcategory of \(\C\),
unless the recollement is trivial, i.e., \(\C=\A\times\B\).

\begin{rem}
  We prefer the formulation of \Cref{defi:recollement}
  because it presents the full datum of a recollement in a way which
  is nicely symmetric between left and right adjoints.
  In this way our definition seemingly differs from Lurie's
  (see \cite[Definition~A.8.1]{lurie:ha}),
  who in a more general left exact setting defines
  recollements asymmetrically only in terms of
  two full subcategories \(\C_0\triangleq\A\) and \(\C_1\triangleq \rorth{\A}\) of \(\C\),
  whose inclusions \(i\) and \(j'\) admit left adjoints \(L_0\triangleq q\) and \(L_1\triangleq p\),
  respectively.
  In the stable setting one then automatically has the other two adjoints \(q'\) and \(j\), see \cite[Remark~A.8.19]{lurie:ha} (this can also be deduced from \cite[Proposition 2.3.2]{DKSS:spherical}),
  justifying our more redundant definition.
\end{rem}

One has the following computation.

\begin{lem}
  \label{lem:gluing-functor-char}
  Given a recollement \eqref{eq:recollement},
  the units and counits of the various adjunctions
  yield canonical identifications
  \begin{equation}
    q' j \cong \fib(j\to j') \cong q j'[-1]
  \end{equation}
  of functors \(\B\to \A\),
  where the canonical map \(j\to j'\) can be obtained either
  by transposing the counit \(pj\xrightarrow{\cong} \id_\B\) along the adjunction \(p\dashv j'\)
  or, equivalently, by transposing the unit \(\id_\B\xrightarrow{\cong} pj'\)
  along the adjunction \(j\dashv p\).
  Note that in the middle there is the implicit claim that
  \(\fib(j\to j')\colon \B\to\C\)
  factors through \(i\), so that we can view it as a functor \(\B\to \A\).
\end{lem}

\begin{defi}
  Any of the equivalent functors \(F\colon\B\to\A\) of \Cref{lem:gluing-functor-char}
  is called the \emph{gluing functor} of the recollement \eqref{eq:recollement}.
\end{defi}

\begin{con}
  \label{con:laxlim-recollement}
  Conversely given a \(F\colon \B\to \A\),
  one can construct a canonical recollement
  \begin{equation}
    \label{eq:recollement-cocart}
    \cdrecollementcompl
    \A{\laxlima\B F \A}\B
    {(0\to a)} {\ev_b}
    {(b \xrightarrow{!} Fb)} {\cof}
    {(b\to 0)} {\ev_a}
  \end{equation}
  where in the middle we have the lax limit
  \begin{equation}
    \laxlima \B F \A = \{(b\to a) = (b : B, a : A, Fb \to a)\},
  \end{equation}
  i.e., the category of sections of the Grothendieck construction
  for the functor \(\Delta^1 \to \St\) classifying \(F\).
  Manifestly, the gluing functor of the recollement 
  \eqref{eq:recollement-cocart} is the original functor \(F\colon \A\to \B\).
\end{con}

The main structural result of the theory is that one can construct recollements
starting with only very minimal amount of data.

\begin{thm}
  A recollement \eqref{eq:recollement} can be uniquely recovered/constructed from
  any one of the following pieces of data:
  \begin{enumerate}
  \item
    \label{it:recoll-from-A}
    A fully faithful functor \(i\colon \A\hookrightarrow \C\)
    which admits a left and a right adjoint.
    In this case \(\B\) is determined as the Verdier quotient \(\C/\A\).
  \item
    \label{it:recoll-from-B}
    A functor \(p\colon \C\to \B\)
    which admits a left and a right adjoint, both fully faithful.
    In this case \(\A\) is determined as the kernel of \(p\).
  \item
    \label{it:recoll-from-F}
    An arbitrary functor \(F\colon \B\to \A\).
    In this case the recollement is determined by 
    \Cref{con:laxlim-recollement}.
  \end{enumerate}
\end{thm}

\begin{proof}
Part \ref{it:recoll-from-A} is the statement of \cite[Proposition~2.3.3]{DKSS:spherical} or \cite[Proposition~A.8.20]{lurie:ha}. Part \ref{it:recoll-from-F} follows from \cite[Proposition~2.2.11]{DKSS:spherical} or \cite[Remark~A.8.18]{lurie:ha}. The only thing to note there is that Lurie's convention is dual to ours:
  instead of reconstructing the recollement from the composite
  \(F=q'j\colon \B\to \A\), he uses the composite \(F'=qj'\colon \B \to \A\)
  (which is written \(L_0|\C_1\) in his notation);
  and instead of considering sections of the coCartesian Grothendieck construction,
  he uses the Cartesian one.
  Passing from one convention to the other is just a matter
  of replacing each stable \(\infty\)-category with its opposite.

  Part \ref{it:recoll-from-B} follows from \cite[Proposition~2.3.2]{DKSS:spherical}: since the inclusion $j\colon \B\to \C$ admits a right adjoint, there is a semiorthogonal decomposition $(\A,\B)$ of $\C$ with $\A=\ker(p)$, which in turn implies that $i\colon \A\to \C$ admits a left adjoint. A similar argument shows that $i$ admits a right adjoint. 
\end{proof}


\bibliographystyle{alpha} 
\bibliography{refs} 

\end{document}